\documentclass[11 pt,fullpage]{article}
\usepackage{amsfonts,amssymb}
\usepackage{ graphicx, mathrsfs}
\usepackage{amsmath,amsthm}
\usepackage{color} 
\usepackage{hyperref}
\usepackage{mathrsfs,dsfont}
\usepackage[right = 2.5cm, left=2.5cm, top = 2.5cm, bottom =2.5cm]{geometry}
\theoremstyle{plain}
\newtheorem{theorem}{Theorem}

\newtheorem{proposition}{Proposition}[section]
\newtheorem{lemma}[proposition]{Lemma}

\newtheorem{definition}{Definition}[section]

\theoremstyle{definition}
\newtheorem{remark}{Remark}

\allowdisplaybreaks

\numberwithin{equation}{section}

\newcommand\R{{\mathbb R}}
\newcommand\T{{\mathbb T}}
\newcommand\Torus{{\mathbb T}}

\newcommand\Z{{\mathbb Z}}

\newcommand{\cF}{\mathcal F}
\newcommand{\cG}{\mathcal G}

\newcommand{\cL}{\mathcal L}

\newcommand{\Real}{\mathbb R}
\newcommand{\Complex}{\mathbb C}
\newcommand{\Integer}{\mathbb Z}
\newcommand{\norm}[1]{\left\lVert#1\right\rVert}
\newcommand{\abs}[1]{\left\vert#1\right\vert}
\newcommand{\set}[1]{\left\{#1\right\}}

\newcommand{\grad}{\nabla}

\newcommand{\Naturals}{\mathbb N}
\newcommand{\jap}[1]{\langle #1 \rangle} 
\newcommand{\ddz}{\, dz}
\newcommand{\ddv}{\, dv}

\newcommand{\dss}{\displaystyle}
\newcommand{\vsp}{\vspace{0.2cm}}

 \newcommand{\dd}{{\, \mathrm d}}

\newcommand{\mk}{\medskip}

\begin{document}
 
\title{Landau damping: paraproducts and Gevrey regularity}
\author{Jacob Bedrossian\footnote{\textit{jacob@cims.nyu.edu}, Courant Institute of Mathematical Sciences. Partially supported by NSF Postdoctoral Fellowship in Mathematical Sciences, DMS-1103765}, \, Nader Masmoudi\footnote{\textit{masmoudi@cims.nyu.edu}, Courant Institute of Mathematical Sciences. Partially supported by NSF  grant DMS-1211806} \, and Cl\'ement Mouhot\footnote{\textit{c.mouhot@dpmms.cam.ac.uk}, Centre for Math. Sc., Univ. of Cambridge. Partially
      funded by ERC grant MATKIT}}

\date{\today}
\maketitle

\begin{abstract}
We give a new, simpler, proof of nonlinear Landau damping on $\Torus^d$ in Gevrey$-\frac{1}{s}$ regularity ($s > 1/3$)
 which matches the regularity requirement predicted by the formal analysis of Mouhot and Villani \cite{MouhotVillani11}.
Our proof combines in a novel way ideas from the original proof of Landau damping \cite{MouhotVillani11} and the proof of inviscid damping in 2D Euler \cite{BM13}. 
As in \cite{BM13}, we use paraproduct decompositions and controlled regularity loss to replace the Newton iteration scheme of \cite{MouhotVillani11}. 
We perform time-response estimates adapted from \cite{MouhotVillani11} to control the plasma echoes and couple them to energy estimates 
on the distribution function in the style of the work \cite{BM13}. 
\end{abstract}


\setcounter{tocdepth}{1}
{\small\tableofcontents}

\section{Introduction} \label{sec:vlas-poiss-equat}

The collisionless Vlasov equations are a fundamental model of plasma
physics and galactic dynamics (see
e.g. \cite{BoydSanderson,Binney-Tremaine,Krall-Trivelpiece,Ryutov99,villani2010}),
and it writes in the periodic box $x \in \T^d_L := [-L,L]^d$ with size $L>0$: 
\begin{equation}\label{def:VPE0}
  \left\{
\begin{array}{l} \dss 
\partial_t f + v\cdot \grad_x f + F(t,x)\cdot \grad_v f  = 0, \\
F(t,x) = -\grad_x W \ast_{x} \left(\rho_f(t,x) - L^{-d}\int_y \rho_f(t,y) dy\right), \\
\rho_f(t,x) = \int_{\R^d} f(t,x,v) \, dv, \\ 
f(t=0,x,v) = f_{in}(x,v),
\end{array}
\right.
\end{equation}
with $f(t,x,v) :\Torus^d _L \times \Real^d \rightarrow [0,\infty)$, the distribution function in phase space. 
We are interested in solutions of the form $f(t,x,v) = f^0(v) +
h(t,x,v)$, where $f^0(v)$ is a spatially homogeneous background
distribution and $h$ is a mean-zero perturbation. 
If we denote simply the (perturbation) density $\rho(t,x)$, then the
Vlasov equations can be written as
\begin{equation} \label{def:VPE}
\left\{
\begin{array}{l} \dss 
\partial_t h + v\cdot \grad_x h + F(t,x)\cdot \grad_v h + F(t,x)\cdot
\grad_vf^0 = 0, \\ 
F(t,x) = -\grad_x W \ast_{x} \rho(t,x), \\ 
\dss \rho(t,x) = \int_{\R^d} h(t,x,v)
\ddv, \\  
h(t=0,x,v) = h_{in}(x,v). 
\end{array}
\right.
\end{equation}
The potential $W(x)$ describes the mean-field interaction between
particles; the cases of most physical interest are (1) Coulomb
repulsive interactions $F = e \nabla_x \Delta_x ^{-1} \rho_f$ between
electrons in plasmas (where $e >0$ is the electron charge-to-mass ratio)
and Newtonian attractive interactions $F = -m \cG \nabla_x \Delta_x
^{-1} \rho_f$ between stars in galaxies (where $m>0$ is the mass of
the identical stars and $\cG$ is the gravitational constant). In
Fourier variables (see later for the notation) these two cases
correspond respectively to $\widehat{W}(k) = (2\pi)^{-2} e L^2 \abs{k}^{-2}$ and
$\widehat{W}(k) = - (2\pi)^{-2} m\cG L^2 \abs{k}^{-2}$.
The former arises in plasma physics where \eqref{def:VPE0} describes
the distribution of electrons in a plasma interacting with a
background of ions ensuring global electrical neutrality, after
neglecting magnetic effects and ion acceleration. 
The latter arises in
galactic dynamics where \eqref{def:VPE0} describes a distribution of stars interacting via Newtonian gravitation, neglecting
smaller planetary objects as well as relativistic effects, and
assuming Jean's swindle (see \cite{MouhotVillani11} and the references therein).

By re-scaling $t,x$ and $W$, we may normalize the size of the box to $L=2\pi$ without
loss of generality and write $\T^d = \T^d_{2\pi}$ (see Remark \ref{rmk:confinement}). 
 For simplicity of notation and mathematical generality, we consider a general class
of potentials with Coulomb/Newton representing the most
singular examples. Specifically, we only require that there exists
$C_W < \infty$ and $\gamma \geq 1$ such that
\begin{align} 
|\widehat{W}(k)| \leq C_W \abs{k}^{-1-\gamma}. \label{ineq:Wbd}
\end{align}

This paper is concerned with the mathematically rigorous treatment of
Landau damping for the full nonlinear mean-field dynamics, as
initiated in \cite{CagliotiMaffei98,HwangVelazquez09,MouhotVillani11}. 
We will not provide a historical
background for Landau damping as the topic was discussed at length
in \cite{Ryutov99,MouhotVillani11,villani2010}; see references therein. 
Let us just briefly recall that Landau damping is a mechanism discovered by Landau
\cite{Landau46} (after preliminary works of
Vlasov~\cite{Vlasov-damping}) predicting the decay of spatial
oscillations in a plasma when perturbed around certain stable
spatially homogeneous distributions. 
This effect is now considered fundamental to modern plasma physics (see e.g. \cite{Ryutov99,BoydSanderson,Stix,MouhotVillani11}). 
It was later ``exported'' to galactic dynamics
by Lynden-Bell \cite{Lynden,LyndenBell67} where it is thought to play a key role in the stability of galaxies.  
Landau damping is one example of a more general effect usually referred to as ``phase mixing'', which 
arises in many physical settings; see \cite{MouhotVillani11,BM13,BMT13} and the references therein.
See also \cite{Degond86,CagliotiMaffei98,MouhotVillani11,BM13} for a discussion about the differences and similarities with dispersive phenomena.
 
The original works in physics neglected nonlinear effects, which lead to some speculation (see \cite{MouhotVillani11,LZ11b} and the references therein). 
The mathematically rigorous theory of the linear
damping was pioneered by Backus \cite{Backus} and Penrose
\cite{Penrose}, and further clarified by many mathematicians, see
e.g. \cite{Maslov,Degond86}.  
The first limited nonlinear results were obtained by \cite{CagliotiMaffei98,HwangVelazquez09} which showed that Landau damping was at least possible in \eqref{def:VPE} for analytic data (see also \cite{LZ11b} for a negative result). 
However, Landau damping was not shown to hold for all small data until \cite{MouhotVillani11}. 
The results therein hold in analytic spaces or in Gevrey spaces ``close to analytic'' \cite{Gevrey18} and the authors made heuristic conjectures about the
minimal regularity required. 
The proof involved an intricate use of Eulerian and Lagrangian coordinates combined with a global-in-time Newton
iteration reminiscent of the proof of the KAM theorem. 
A key step in the analysis of \cite{MouhotVillani11} was controlling the potentially destabilizing influence of \emph{plasma echoes}, a weakly nonlinear effect discovered by Malmberg et. al. in \cite{MalmbergWharton68}. See \cite{MouhotVillani11} for a detailed discussion of the role this effect plays in the nonlinear theory. 

In this work, we provide a new, simpler and shorter proof of Landau damping for \eqref{def:VPE} that nearly obtains the ``critical'' Gevrey regularity predicted in \cite{MouhotVillani11}. 
From a physical point of view, our proof has many of the same ingredients as the proof  in \cite{MouhotVillani11}:
1) we use the same abstract linear stability condition, 2) the interplay between regularity and decay is similar and 3) the isolation and control of the plasma echoes is still the main challenge. 
However, on a mathematical level,  the two proofs are quite different; see \S\ref{sec:Outline} for a full discussion.  
In short, our proof combines the viewpoint in the original work \cite{MouhotVillani11} with the recent work on the 2D Euler equations in \cite{BM13} (see also the expository note \cite{BM13_AMRX}). 
This latter work proves the asymptotic stability (in a suitable sense) of sufficiently smooth shear flows near Couette flow in $\Torus \times \Real$ via `inviscid damping' (the hydrodynamic analogue of Landau damping). 
The proof in \cite{BM13} uses a number of ideas specific to the structure of 2D Euler, however, some aspects of the viewpoint taken therein will be useful here as well (when suitably combined with ideas from \cite{MouhotVillani11}). 
  
One of the main ingredients of our proof, employed also in \cite{BM13}, is the use of the paradifferential calculus to split nonlinear terms into one that carries the transport structure and another that is analogous to the nonlocal interaction term referred to as `reaction' in \cite{MouhotVillani11}. 
It has been long believed that Nash-Moser or similar Newton iterations (see the classical work \cite{Moser1966}) can generally be replaced by a more standard  fixed point argument if one uses better all of the structure in the
equation. This has been the case in most examples in the literature (e.g. Nash's isometric embedding theorem \cite{MR0075639} by G\"unther \cite{MR1029846,MR1037168,MR1159298}), with maybe the only exception being KAM theory. 
Other examples can be found in \cite{Hormander1990}, where H\"ormander specifically points out paradifferential calculus as 
a useful tool in this context, stating: 

``The Nash-Moser method contains constructions which very much resemble the dyadic
decompositions which are central in the paradifferential calculus of Bony. In fact,
we show that the Nash-Moser technique can often be replaced by elementary nonlinear
functional analysis combined with the paradifferential calculus of Bony.''

It is well known that one of the main physical barriers to Landau damping 
is nonlinear particle trapping, whereby particles are trapped in the potential wells of (say) electrostatic waves.
Exact traveling wave solutions of this type exist in plasma physics and are known as BGK waves \cite{BGK57}. 
They were used by Lin and Zeng in \cite{LZ11b} to show that
 one needs at least $H^\sigma$ with $\sigma > 3/2$ regularity on the distribution function to expect Landau damping in \eqref{def:VPE} in the neighborhood of Landau-stable stationary solutions.
The plasma echoes provide a natural nonlinear bootstrap mechanism by which the electrostatic waves can persist long enough to trap particles.
After modding out by particle free streaming, the bootstrap appears as a cascade of information to high frequencies and
the regularity requirement of Gevrey-$(2+\gamma)^{-1}$ comes from formal `worst-case' calculations carried out in \cite{MouhotVillani11}.
Mathematically, the same requirement arises here in \S\ref{sec:Proofa}.
Lower regularity is an open question: it seems plausible that Theorem \ref{thm:Main} is false for all $s < (2+\gamma)^{-1}$, 
however there might be additional cancellations 
that could allow it to hold in lower regularities.
Finally, in weakly collisional plasmas, the requirement could perhaps be relaxed in some suitable sense (e.g. permitting data which is Gevrey plus a smaller rough contribution that will be instantly regularized).

It is well known that many areas in physics present striking analogies with each other and
  many  important developments  have come 
 from a good
 understanding of these analogies. 
 This work is  an example where
 the analogy   between  2D incompressible  Euler  and Vlasov-Poisson
  proved fruitful.  
The connection between inviscid damping and Landau damping has been acknowledged by a number of authors, for example \cite{BouchetMorita10,SchecterEtAl00,Briggs70,BM95,LinZeng11,GSV13}.
Both have similar weakly nonlinear echoes \cite{YuDriscoll02,YuDriscollONeil,VMW98,Vanneste02}, moreover, the work of Lin and Zeng \cite{LinZeng11,LZ11b} and \cite{BMT13} show that particle trapping and vortex roll-up may in fact be essentially the same `universal' over-turning instability.  
On a more general level, both systems are conservative transport equations governed by a single scalar quantity (the vorticity in 2D Euler and the distribution function in Vlasov).  
 Both equations have a large set of stationary 
solutions: the shear flows for Euler and spatially homogeneous distributions for Vlasov being the simplest.     
Each can be viewed as Hamiltonian systems and variational methods have
been used for both to provide nonlinear stability results in
low-regularity spaces (i.e. functional spaces invariant under the
free-streaming operator): we refer for
instance, among a huge literature, to \cite{AK98} in the context of
the Euler equation, and to the recent remarkable results
\cite{Guo-Lin,LMR12} in the context of the gravitational
Vlasov-Poisson equation (see also the review paper \cite{bourbaki-vp}
and the references therein).
One can also derive  the incompressible Euler system from 
the {V}lasov-{P}oisson system in the quasi-neutral regime for 
cold electrons  (vanishing initial temperature) \cite{Brenier00,Masmoudi01}. 

\subsection{The main result}
\label{sec:main-result}
In what follows, we denote the Gevrey-$\frac{1}{\nu}$ norms, $\nu \in (0,1]$, with Sobolev corrections $\sigma \in \Real$
\begin{subequations} \label{def:norm}
\begin{align} 
\norm{h}^2_{\cG^{\lambda,\sigma;\nu}} = \sum_{k \in \Integer^d}\int_\eta \abs{\hat{h}_k(\eta)}^2 \jap{k,\eta}^\sigma e^{2\lambda\jap{k,\eta}^{\nu}} d\eta \\ 
\norm{\rho(t)}^2_{\cF^{\lambda,\sigma;\nu}} = \sum_{k \in \Integer^d} \abs{\hat{\rho}_k(t)}^2 \jap{k,kt}^\sigma e^{2\lambda\jap{k,kt}^{\nu}}. 
\end{align} 
\end{subequations}
When $\sigma = 0$ or $\nu = s$, (defined below) these parameters are usually omitted.
See \S\ref{sec:Notation} for other notation conventions and \S\ref{sec:toolbox} for Fourier analysis conventions.  

Recall the sufficient linear stability condition \textbf{(L)} introduced in \cite{MouhotVillani11} for analytic background distributions, which we slightly adapt to the norms we are using. 
In \cite{MouhotVillani11}, it is shown that \textbf{(L)} is practically equivalent to several well known stability criteria in plasma and galactic dynamics (see \S\ref{sec:penrose-criterion} for a completion of the proof that the Penrose condition \cite{Penrose} implies \textbf{(L)}). 

\begin{definition} \label{def:L} 
Given a homogeneous analytic distribution $f^0(v)$ we say it satisfies \emph{stability condition \textbf{(L)}} if 
there exists constants $C_0,\bar \lambda,\kappa > 0$ and an integer $M > d/2$ such that 
\begin{align} 
\sum_{\alpha \in \Naturals^d : \abs{\alpha} \leq M}\norm{v^{\alpha}f^0}^2_{\cG^{\bar\lambda;1}} \leq C_0, \label{ineq:f0Loc}
\end{align} 
and for all $\xi \in \Complex$ with $\textup{Re}\, \xi < \bar\lambda$, 
\begin{align} 
\inf_{k \in \Integer^d}\abs{\mathcal{L}(\xi,k) - 1} \geq \kappa,  
\end{align} 
where $\mathcal{L}$ is defined by the following (where $\bar{\xi}$ denotes the complex conjugate of $\xi$), 
\begin{align} 
\mathcal{L}(\xi,k) = -\int_0^\infty e^{\bar{\xi}\abs{k}t} \widehat {f^0}\left(kt\right) \widehat{W}(k)\abs{k}^2t dt. \label{ineq:mL}
\end{align}
\end{definition}
\begin{remark} 
Note \eqref{ineq:f0Loc} implies the integral in \eqref{ineq:mL} is absolutely convergent by the $H^{d/2+}\hookrightarrow C^0$ embedding theorem. 
\end{remark}

We prove the following nonlinear Landau damping result, which for Coulomb/Newton interaction nearly obtains the Gevrey-3 regularity predicted heuristically in \cite{MouhotVillani11}.
In \S\ref{sec:Outline} we give the outline of the proof and discuss the relationship with the original proof in \cite{MouhotVillani11} and the proof of inviscid damping in 2D Euler \cite{BM13}. 

\begin{theorem}\label{thm:Main} 
Let $f^0$ be given which satisfies stability condition \textbf{(L)} with constants $M$,$\bar{\lambda}$, $C_0$ and $\kappa$. 
Let $\frac{1}{(2+\gamma)} < s \leq 1$ and $\lambda_0 > \lambda^\prime > 0$ be arbitrary (if $s = 1$ we require $\bar\lambda > \lambda_0$).
Then there exists an $\epsilon_0 = \epsilon_0(d,M,\bar{\lambda},C_0,\kappa,\lambda_0,\lambda^\prime,s)$ such that if $h_{in}$ is mean zero and
\begin{align*} 
\sum_{\alpha \in \Naturals^d : \abs{\alpha} \leq M}\norm{v^{\alpha}h_{in}}^2_{\cG^{\lambda_0;s}} < \epsilon^2 \leq \epsilon^2_0, 
\end{align*} 
then there exists a mean zero $h_\infty \in \cG^{\lambda';s}$ satisfying for all $t \ge0$,
\begin{subequations} \label{ineq:damping}
\begin{align} 
\norm{h(t,x+vt,v) - h_\infty(x,v)}_{\cG^{\lambda^{\prime};s}} & \lesssim \epsilon e^{-\frac{1}{2}(\lambda_0 - \lambda^\prime)t^{s}}, \label{ineq:glidingconverg}\\
\norm{\rho(t)}_{\cF^{\lambda^\prime;s}} & \lesssim  \epsilon e^{-\frac{1}{2}(\lambda_0 - \lambda^\prime)t^{s}}. \label{ineq:densitydecay}
\end{align}  
\end{subequations}
\end{theorem} 

\begin{remark} \label{rmk:confinement}
Through the rescaling on $W$, our estimate of $\epsilon_0$ in Theorem \ref{thm:Main} is a decreasing function of the side-length of the original torus, $L$. That is, the restriction for nonlinear stability becomes more stringent as the confinement is removed. Moreover, through the rescaling on time, Theorem \ref{thm:Main} predicts damping on a characteristic time-scale of $O(L)$. See \cite{glassey94,glassey95,MouhotVillani11,villani2010} for more discussion about what can happen without confinement. 
\end{remark}

\begin{remark} It is immediate to deduce estimates on the complete
  distribution $f(t,x,v) = f^0(v) + h(t,x,v)$ of the original
  equation~\eqref{def:VPE0}. In particular, Theorem \ref{thm:Main} shows that all solutions to \eqref{def:VPE0} with non-trivial spatial dependence close to $f^0$ satisfy $\norm{f(t)}_{H^N} \approx \jap{t}^N$ (the same as free transport).  
\end{remark}

\begin{remark} 
From \eqref{ineq:glidingconverg} we have the homogenization $h(t,x,v) \rightharpoonup <f_\infty(\cdot,v)>_x := (2\pi)^{-d}\int f(x,v) dx$ and the  
  exponential decay of the electrical or gravitational field
  $F(t,x)$. 
Note that the estimates on $h$ and
  $\rho$ are more precise as they show decay rates which increase with the spatial frequency. 
\end{remark}

\begin{remark} 
The theorem also holds backwards in time for some $h_{-\infty} \in \cG^{\lambda';s}$. 
\end{remark} 
 
\begin{remark} \label{rmk:LinearTh}
The asymptotic distribution function $f_\infty(x,v) := f^0(v) + h_\infty(x,v)$ depends on the
  entire nonlinear evolution, however, 
  at least one can show that $f_\infty$ is within $O(\varepsilon^2)$ of the
  distribution predicted by the linear theory 
in $\cG^{\lambda^{\prime\prime};s}$ for any $\lambda^{\prime\prime} < \lambda^\prime$ (similar to
  \cite{MouhotVillani11}). See
  \S\ref{sec:ConcProof} for a sketch.
\end{remark}

\begin{remark} 
Though $f^0$ is analytic, the statement shows the asymptotic stability of homogeneous distributions within a small neighborhood of $f^0$ in Gevrey-$\frac{1}{s}$ since we are really making a perturbation of $f^0(v) + <h_{in}(\cdot,v)>_x$.
However, the size of that neighborhood depends on the parameters in Definition \ref{def:L}, so it still must be close to an analytic distribution satisfying \textbf{(L)}.   
\end{remark} 

\begin{remark} 
Requiring $h_{in}$ to be average zero does not lose any generality. Indeed, if $h_{in}$ were not mean zero we may apply Theorem \ref{thm:Main} to $\tilde{f}^0(v) = f^0(v) +\left( <h_{in}(\cdot,v)>_x\right)_{<1}$ and $\tilde h_{in} = h_{in} -\left( <h_{in}(\cdot,v)>_x\right)_{<1}$, where $g_{<1}$ denotes projection onto frequencies less than one. 
Since \textbf{(L)} is open, for $\epsilon_0$ sufficiently small, $\tilde{f}^0$ still satisfies \textbf{(L)} with slightly adjusted parameters $C_0$ and $\kappa$.
\end{remark} 

\begin{remark} 
If $W$ is in the Schwartz space, then Theorem \ref{thm:Main} holds for all $0 < s < 1$ (although $\epsilon_0$ goes to zero as $s \searrow 0$). 
The heuristics of \cite{MouhotVillani11} suggest that even in the case of analytic $W$, we cannot hope to work in the Sobolev scale without major new ideas, if such a result holds at all. 
\end{remark}

\subsubsection{Notation and conventions} \label{sec:Notation}
We denote $\Naturals = \set{0,1,2,...}$ (including zero) and $\Integer_\ast = \Integer \setminus \set{0}$. 
For $\xi \in \Complex$ we use $\bar{\xi}$ to denote the complex conjugate. 
We will denote the $\ell^1$ vector norm $\abs{k,\eta} = \abs{k} + \abs{\eta}$, which by convention is the norm taken in our work. 
We denote 
\begin{align*} 
\jap{v} = \left( 1 + \abs{v}^2 \right)^{1/2}. 
\end{align*} 
We use the multi-index notation: given $\alpha = (\alpha_1,...,\alpha_d) \in \Naturals^d$ and $v = (v_1,...,v_d) \in \Real^d$ then
\begin{align*} 
v^\alpha  = v^{\alpha_1}_1v^{\alpha_2}_2...v^{\alpha_d}_d, \quad\quad\quad D_\eta^\alpha  = (i\partial_{\eta_1})^{\alpha_1} ... (i\partial_{\eta_d})^{\alpha_d}. 
\end{align*}  
We denote spatial averages $<f(\cdot,v>_x = (2\pi)^{-d}\int_{\T^d} f(x,v) dx$ and denote Lebesgue norms for $p,q \in [1,\infty]$ and $a,b$ either in $\Real^d$, $\Integer^d$ or $\Torus^d$ as
\begin{align*}  
\norm{f}_{L_a^p L_b^q} = \left(\int_a \left(\int_b \abs{f(a,b)}^q db \right)^{p/q} da\right)^{1/p}
\end{align*}
and Sobolev norms (usually applied to Fourier transforms) as
\begin{align*} 
\norm{\hat{f}}^2_{H^M_\eta} = \sum_{\alpha \in \Naturals^d; \abs{\alpha} \leq M} \norm{D_\eta^\alpha \hat{f}}_{L^2_\eta}^2. 
\end{align*}
We will often use the short-hand $\norm{\cdot}_2$ for
  $\norm{\cdot}_{L^2_{z,v}}$ or $\norm{\cdot}_{L^2_v}$ depending on the context.

See \S\ref{Apx:LPProduct} for the Fourier analysis conventions we are taking.
A convention we generally use is to denote the discrete $x$ (or $z$) frequencies as subscripts.   
We use Greek letters such as $\eta$ and $\xi$ to denote frequencies in $\Real^d$ and lowercase Latin characters such as $k$ and $\ell$ to denote frequencies in $\Integer^d$. 
Another convention we use is to denote $N,N^\prime$ as dyadic integers $N \in \mathbb{D}$ where  
\begin{align*} 
\mathbb{D} = \set{\frac{1}{2},1,2,...,2^j,...}. 
\end{align*}
When a sum is written with indices $N$ or $N^\prime$ it will always be over a subset of $\mathbb D$. 
This will be useful when defining Littlewood-Paley projections and paraproduct decompositions, see \S\ref{Apx:LPProduct}. 
Given a function $m \in L^\infty$, we define the Fourier multiplier $m(\grad_{z,v}) f$ by 
\begin{align*} 
(\widehat{m(\grad_{z,v})f})_k(\eta) =  m( (ik,i\eta) ) \hat{f}_k(\eta). 
\end{align*}   
For any multiplier $m$ we define the following notation for $m(\grad_{z,v})\rho(t,x)$ (also for similar functions only of $x$) which will appear more natural below: 
\begin{align}  
\widehat{(m(\grad_{z,v})\rho)}_k(t) & = m((ik,ikt))\hat{\rho}_k(t). \label{def:Fmrho} 
\end{align}
We use the notation $f \lesssim g$ when there exists a constant $C > 0$ independent of the parameters of interest 
such that $f \leq Cg$ (we analogously define $f \gtrsim g$).
Similarly, we use the notation $f \approx g$ when there exists $C > 0$ such that $C^{-1}g \leq f \leq Cg$. 
We sometimes use the notation $f \lesssim_{\alpha} g$ if we want to emphasize that the implicit constant depends on some parameter $\alpha$.

\section{Outline of the proof} \label{sec:Outline}
\subsection{Summary and comparison with original proof \cite{MouhotVillani11}} \label{sec:main-ingr-proof} 
Landau damping predicts that the solution evolves by kinetic free transport as $t \rightarrow \infty$: 
\begin{align*} 
h(t,x,v) \sim h_\infty(x-vt,v).
\end{align*} 
We `mod out' by the characteristics of free transport and work in the coordinates 
$z = x-vt$ with $f(t,z,v) = h(t,x,v)$ (not to be confused with the notation in \eqref{def:VPE0}). 
Then \eqref{ineq:glidingconverg} becomes equivalent to $f(t) \rightarrow h_\infty$ strongly in Gevrey$-\frac{1}{s}$. 
This coordinate shift was used in \cite{CagliotiMaffei98,HwangVelazquez09} and is related to the notion of `gliding regularity' used in \cite{MouhotVillani11} (although we will not have to use the time-shifting tricks employed therein).  
Moreover, it is also closely related to the notion of `profile' used in nonlinear dispersive equations (see e.g. \cite{GMS12}).
A related coordinate change is used in \cite{BM13}, although there it is complicated by the fact that it depends on the solution.

A straightforward computation gives the evolution equation:
\begin{equation}\label{def:VPEgliding}
\left\{
\begin{array}{l}
\partial_t f + F(t,z+vt) \cdot (\grad_v -t\grad_z)f + F(t,z+vt)\cdot \grad_vf^0 = 0, \\ 
f(t = 0,z,v) = h_{in}(z,v). 
\end{array}
\right.
\end{equation}
Note that the density satisfies $\hat{\rho}_k(t) =  \hat f_k(t,kt)$; by the $H^{d/2+} \hookrightarrow C^0$ embedding theorem and the requirement $M > d/2$, this formula at least makes sense pointwise in time provided $\sum_{\alpha \leq M}\norm{v^\alpha f}_2$ is finite.  
Moreover, from this formula we see that a uniform bound on the regularity of $f$ translates directly into decay of the density: this is precisely the phase mixing mechanism. 
Hence, to prove \eqref{ineq:densitydecay}, we are aiming for a uniform-in-time bound on the regularity in a velocity polynomial-weighted space (so that we may restrict the Fourier transform using $H^{d/2+} \hookrightarrow C^0$). 

All of our analysis is on the Fourier side; taking Fourier transforms and using $\widehat{F(t,z+vt)}(t,k,\eta) = -(2\pi)^dik\widehat{W}(k)\hat{\rho}_k(t) \delta_{\eta = kt}$, \eqref{def:VPEgliding} becomes
\begin{align} 
\partial_t  \hat f_k(t,\eta) + \hat{\rho}_k(t) \widehat{W}(k) k \cdot (\eta - tk) \hat{f}^0(\eta - kt) + \sum_{\ell \in \Integer_\ast^d} \hat{\rho}_\ell(t)\widehat{W}(\ell)\ell\cdot \left[ \eta  - tk\right] \hat f_{k-\ell}(t,\eta - t\ell) = 0. \label{eq:feqn}
\end{align}
Note in the summation that $\eta - tk = (\eta-t\ell) - t(k-\ell)$.
By the formula $\hat{\rho}_k(t) = \hat{f}_k(t,kt)$, \eqref{eq:feqn} is a closed equation for $f$, however, we will derive a separate integral equation for $\rho(t)$ and look at $(f,\rho)$ as a coupled system (similar to \cite{MouhotVillani11}). 
Integrating \eqref{eq:feqn} in time and evaluating at $\eta = kt$ gives 
\begin{align} 
\hat{\rho}_k(t) &= \hat h_{in}(k,kt)   - \int_0^t\hat{\rho}_k(\tau) \widehat{W}(k)k \cdot k(t - \tau) f^0\left(k(t-\tau)\right) d\tau \nonumber \\ & \quad
 - \int_0^t \sum_{\ell \in \Integer_\ast^d} \hat{\rho}_\ell(\tau)\widehat{W}(\ell)\ell \cdot k\left(t-\tau\right) \hat f_{k-\ell}(\tau,kt - \tau \ell) d\tau. \label{eq:rhointegral}
\end{align} 

As in \cite{MouhotVillani11}, our goal now is to use the system \eqref{eq:feqn},\eqref{eq:rhointegral} to derive a uniform control on the regularity of $f$ in the moving frame (referred to as `gliding regularity' in \cite{MouhotVillani11}).  
The linear term in \eqref{eq:rhointegral} is handled with the help of the abstract stability condition \textbf{(L)}; the difference here with \cite{MouhotVillani11} is that we must adapt this to the Gevrey norms we are using, which is done using a slightly technical decomposition technique similar to one which appeared in \cite{MouhotVillani11} to treat Gevrey data (carried out below in \S\ref{sec:LinDamp}). 
The main point of departure from \cite{MouhotVillani11} is our treatment of the nonlinear terms in \eqref{eq:rhointegral} and \eqref{eq:feqn}. 
There are schematically two mechanisms of potential loss of regularity in \eqref{eq:feqn} and \eqref{eq:rhointegral} and one potential loss of localization in velocity space: 
\begin{enumerate}
\item Equation \eqref{eq:feqn} describes a transport structure in $\Torus^d \times \Real^d$, and hence we can expect this to induce the loss of regularity usually associated with transport by controlled coefficients. 
A different loss of regularity occurs in \eqref{eq:rhointegral} due to the $ k(t-\tau)\hat{f}_{k-\ell}(\tau,kt-\ell \tau)$: here 
there is a derivative of $f$ appearing but no transport structure to take advantage of. 
However, we still refer to this effect as `transport' and remark that it seems related to the beam instability \cite{BoydSanderson}. 

\item Following \cite{MouhotVillani11}, we can see from \eqref{eq:rhointegral} that $\hat{\rho}_\ell(\tau)$ has a strong effect on $\hat{\rho}_k(t)$ when $kt \sim \ell \tau$, which is referred to as \emph{reaction}. 
  These nonlinear resonances are exactly the plasma echoes of \cite{MalmbergWharton68}, and arise due to interaction with the oscillations in the velocity variable. 
  These effects can potentially amplify high frequencies (e.g. \emph{costing regularity}) at localized times
  of strong interaction. 
  A related effect will appear in \eqref{eq:feqn} when $\rho$ forces $f$ via interaction with lower frequencies, and we will refer to this as reaction as well. 
  \item The density $\rho$ is a restriction of the Fourier transform of $f$ and the nonlinear terms in \eqref{eq:feqn} and \eqref{eq:rhointegral} each involve Fourier restrictions.  
      This issue was treated in \cite{MouhotVillani11} by adapting
  carefully the norms used in order to keep under control some
  $L^1$-based norms of regularity.
\end{enumerate}

The proof of \cite{MouhotVillani11} employs a global-in-time Newton scheme which loses a decreasing amount of analytic regularity at each step.  
The linearization of the Newton scheme provided a natural way to isolate 
transport effects from reaction effects.
The transport was treated by Lagrangian methods to estimate analytic regularity along the characteristic trajectories.
The reaction effects were treated via time-integrated estimates on \eqref{eq:rhointegral} that account carefully for the localized, time-delayed effects of the plasma echoes. 
The proof treats \eqref{eq:rhointegral} and \eqref{eq:feqn} as a coupled system in the sense that the main estimates are two coupled but separate controls, one on $f$ and one on $\rho$.
To extend the results to Gevrey regularity, a frequency decomposition of the initial data was employed so that at each step in the scheme everything was analytic.  

Here the two different mechanisms of \emph{transport} and
\emph{reaction} are recovered by a rather different approach, employed recently in the proof of inviscid damping \cite{BM13}. We use a
paraproduct decomposition in order to split the bilinear terms as 
\begin{equation*}
  G_1 G_2 = (G_1)_{\mbox{{\scriptsize lower}}}
  (G_2)_{\mbox{{\scriptsize higher}}} +
  (G_1)_{\mbox{{\scriptsize higher}}} 
(G_2)_{\mbox{{\scriptsize lower}}} + (G_1
  G_2)_{\mbox{{\scriptsize similar frequencies}}}. 
\end{equation*}
In a general sense, one of the first two terms on the RHS will capture the transport effects and the other will capture the reaction effects. 
The last term is a remainder which roughly corresponds to the quadratic error term in the Newton iteration.  
Indeed, the paraproduct decomposition can be thought of formally linearizing the evolution of higher frequencies around lower frequencies.  

As in \cite{BM13}, the transport terms in \eqref{eq:feqn} are treated via an adaptation of the Gevrey energy methods of
\cite{FoiasTemam89,LevermoreOliver97}.  
The essential content is a \emph{commutator estimate}
to take advantage of the cancellations inherent in the natural transport structure. 
For this step to work we need to use $L^2$ based norms, and so to deal with the Fourier restrictions, we use norms with polynomial weights in velocity. 
 
The reaction effects in \eqref{eq:rhointegral} are treated here by making use of a 
refinement of the integral estimates on $\rho$ of \cite{MouhotVillani11} but with some important
conceptual changes inspired from \cite{BM13}: the loss of Gevrey
regularity can occur along time rather than iteratively in a Newton scheme. 
Together with the paraproduct decompositions, this allows for avoiding the Newton iteration altogether and a significant shortening
of the proof. 

 Since we do not use a Newton scheme, we have an additional new constraint: we are not allowed to lose any
derivatives in our coupled estimates on $\rho$ and $f$, a problem due to the derivative of $f$ appearing in \eqref{eq:rhointegral}.
However, since regularity can be traded for decay, we solve this problem by propagating controls on both
``high'' and ``low'' norms of regularity -- the ``high'' ones being
mildly growing in time, a general scheme which is common in the study of wave and dispersive equations (see e.g. \cite{Klainerman,LR05,GMS12} and the references therein). 

\subsection{Gevrey functional setting} 
As discussed in the previous section, our goal now is to use the system \eqref{eq:feqn},\eqref{eq:rhointegral} to derive a uniform control on the regularity of $f$. 
Unlike the norms used in \cite{BM13} and \cite{MouhotVillani11} we will only need the standard norm $\cG^{\lambda,\sigma;s}$ (and the variant $\cF^{\lambda,\sigma;s}$) defined in \eqref{def:norm} with time-dependent $\lambda(t)$.  
For future notational convenience, we define the multiplier $A(t)$ such that $\norm{A(t)f}_2 = \norm{f}_{\cG^{\lambda(t),\sigma;s}}$: 
\begin{align} 
A_k(t,\eta) & = e^{\lambda(t)\jap{k,\eta}^s} \jap{k,\eta}^\sigma, \label{eq:Akteta}
\end{align}
where $\sigma > d+8$ is fixed and $\lambda(t)$ is an index (or `radius') of regularity which is decreasing in time.
In the sequel, we also denote for $\sigma^\prime \in \Real$:
\begin{align*} 
A^{(\sigma^\prime)}_k(t,\eta) & = e^{\lambda(t)\jap{k,\eta}^s} \jap{k,\eta}^{\sigma + \sigma^\prime}.  
\end{align*}   

We will use our choices of $\lambda$ and $s$ to absorb the potential loss of regularity due to the plasma echoes. 
In particular, the choice $s > (2+\gamma)^{-1}$ will be necessary to ensure that the nonlinear plasma echoes do not destabilize the phase mixing mechanism. 
This restriction is used only in equations \eqref{ineq:IRhard} and \eqref{ineq:IRHardDual} of \S\ref{sec:Proofa}. 
Additionally, in order to absorb the loss of regularity from the echoes, we will need to choose $\lambda(t)$ to decay over a very long time; this restriction is also only used in \eqref{ineq:IRhard} and \eqref{ineq:IRHardDual} of \S\ref{sec:Proofa}. 
It will suffice to make the following choices: 
\begin{align*} 
s > (2+\gamma)^{-1}, \quad\quad \alpha_0 = \frac{\lambda_0}{2} + \frac{\lambda^\prime}{2}, \quad\quad a = \frac{(2+\gamma)s-1}{(1+\gamma)},  
\end{align*}
(note $0 < a < s$ if $s < 1$) and  
\begin{align} 
\lambda(t) = \frac{1}{8}\left(\lambda_0 - \lambda^\prime\right)\left(1-t\right)_+ + \alpha_0 + \frac{1}{4}\left(\lambda_0 - \lambda^\prime\right)\min\left(1,\frac{1}{t^a}\right).  \label{def:lambda}
\end{align} 
Then $\alpha_0 \leq \lambda(t) \leq \frac{7}{8}\lambda_0 + \frac{1}{8}\lambda^\prime$ and the derivative never vanishes: 
\begin{align} 
\dot{\lambda}(t) \lesssim -a(\lambda_0-\lambda^\prime)\jap{t}^{-1-a}. \label{ineq:dotlambda}
\end{align}
Norms such as $\mathcal{G}^{\lambda(t),\sigma;s}$ are common when dealing with analytic or Gevrey regularity, for example, see the works \cite{FoiasTemam89,LevermoreOliver97,Chemin04,KukavicaVicol09,CGP11,GM13,BM13,MouhotVillani11}. The Sobolev correction $\sigma$ is included mostly for technical convenience. 
In order to study the analytic case, $s = 1$, we would need to add an additional Gevrey-$\frac{1}{\zeta}$ correction with $(2+\gamma)^{-1} < \zeta < 1$ as an intermediate regularity (and define $a$ in terms of $\zeta$) so that we may take advantage of beneficial properties of Gevrey norms (see, for example, Lemma \ref{lem:GevProdAlg}). 
For the duration of the proof we ignore this minor technical detail and assume $s < 1$. 

The reason for using $\jap{\cdot}$ as opposed to $\abs{\cdot}$ is so that for all $\alpha \in \Naturals^d$ and $\sigma^\prime \in \Real$,  
\begin{align}  \label{ineq:MomentDerivs}
\abs{D_\eta^\alpha A^{(\sigma^\prime)}_k(t,\eta)} & \lesssim_{\abs{\alpha},\lambda_0,\sigma^\prime} \frac{1}{\jap{k,\eta}^{\abs{\alpha}(1-s)}}A^{(\sigma^\prime)}_k(t,\eta), 
\end{align} 
which is useful when estimating velocity moments.

\subsection{Uniform in time regularity estimates} \label{sec:cont}
In this section we set up the continuity argument we use to derive a uniform bound on $(\rho,f)$. 
In order to ensure the formal computations are rigorous, we first regularize the initial data to be analytic. 
The following standard lemma provides local existence of a unique analytic solution which remains analytic as long as 
a suitable Sobolev norm remains finite. 
The local existence of analytic solutions can be proved with an abstract Cauchy-Kovalevskaya theorem, see for example \cite{Nirenberg72,Nishida77}. 
The propagation of analyticity by Sobolev regularity can be proved by a variant of the arguments in \cite{LevermoreOliver97} along with the inequality $\norm{B\rho}_2 \lesssim \sum_{\alpha \leq M} \norm{v^\alpha B f}_2$ for all Fourier multipliers $B$ (with our notation \eqref{def:Fmrho}) and all integers $M > d/2$. 
We omit the proof of Lemma \ref{lem:loctheory} for brevity. 

\begin{lemma}[Local existence and propagation of analyticity] \label{lem:loctheory}
Let $M > d/2$ be an integer and $\tilde{\lambda}>0$. Suppose $h_{in}$ is analytic such that 
\begin{align*} 
\sum_{\alpha \in \Naturals^d: \abs{\alpha} \leq M}\norm{v^\alpha h_{in}}_{\cG^{\tilde{\lambda};1}} < \infty.  
\end{align*} 
Then there exists some $T_0 > 0$ such that for all $T < T_0$ there exists a unique analytic solution $f(t)$ to \eqref{def:VPEgliding} on $[0,T]$ such that for some index $\tilde{\lambda}(t)$ with $\inf_{t \in [0,T]} \tilde\lambda(t) > 0$ we have,  
\begin{align*} 
\sup_{t \in [0,T]}\sum_{\alpha \in \Naturals^d: \abs{\alpha} \leq M}\norm{v^\alpha f(t)}_{\cG^{\tilde{\lambda}(t);1}} < \infty.  
\end{align*} 
Moreover, if for some $T \leq T_0$ and $\sigma > d/2$, $\limsup_{t \nearrow T} \sum_{\alpha \in \Naturals^d:\abs{\alpha} \leq M}\norm{v^\alpha f(t)}_{H^{\sigma}_{x,v}} < \infty$, then $T < T_0$. 
\end{lemma}
\begin{remark} 
If $d \leq 3$ and the solution has finite kinetic energy, then global analytic solutions to \eqref{def:VPE} exist even for large data. 
See the classical results \cite{Pfaffelmoser92,Schaeffer91,LP91,Horst93,Benachour89} for the global existence of strong solutions (from which analyticity can be propagated by a variant of e.g. \cite{LevermoreOliver97}).
\end{remark}
\begin{remark} 
To treat the case $s = 1$ we would need to be slightly more careful in applying Lemma \ref{lem:loctheory}. 
In this case, we may regularize the data to a larger radius of analyticity $\tilde{\lambda} > \lambda_0$, perform our a priori estimates until $\tilde{\lambda}(t) = \lambda(t)$, at which point we may stop, re-regularize and restart iteratively. 
\end{remark} 

Lemma \ref{lem:loctheory} implies that as long as we retain control on the moments and regularity of the regularized solutions, they exist and remain analytic. We perform the a priori estimates on these solutions, for which the computations are rigorous, and then we may pass to the limit to show that the original solutions exist globally and satisfy the same estimates as the regularized solutions. For the remainder of the paper, we omit these details and discuss only the a priori estimates.  

For constants $K_i$, $1 \leq i \leq 3$ fixed in the proof depending only on $C_0,\bar{\lambda},\kappa$, $s$, $d$, $\lambda_0$ and $\lambda^\prime$, let $I \subset \Real_+$ be the largest interval of times such that $0 \in I$ and the following controls hold for all  $t \in I$: 
\begin{subequations} \label{ctrl:Boot}
\begin{align}
\sum_{\alpha \in \Naturals^d:\abs{\alpha} \leq M}\norm{\jap{\grad_{z,v}}A(v^\alpha f)(t)}^2_{2} = \sum_{\alpha \in \Naturals^d:\abs{\alpha} \leq M}\norm{A^{(1)}(v^\alpha f)(t)}^2_{2} & \leq 4K_1\jap{t}^7\epsilon^2 \label{ctrl:HiLocalized} \\
\sum_{\alpha \in \Naturals^d:\abs{\alpha} \leq M}\norm{\jap{\grad_{z,v}}^{-\beta}A(v^\alpha f)(t)}^2_{2} = \sum_{\alpha \in \Naturals^d:\abs{\alpha} \leq M}\norm{A^{(-\beta)}(v^\alpha f)(t)}^2_{2} & \leq 4K_2\epsilon^2 \label{ctrl:LowCommLoc} \\
\int_0^t \norm{A\rho(\tau)}_2^2 d\tau & \leq 4K_3\epsilon^2,  \label{ctrl:MidBoot}
\end{align}
\end{subequations}
where we may fix $\beta > 2$ arbitrary. Recall the definition of $A$ in \eqref{eq:Akteta}.
It is clear from the assumptions that if $K_i \geq 1$ then $0 \in I$. 
The primary step in the proof of Theorem \ref{thm:Main} is to show that $I = [0,\infty)$.

For the regularized solutions it will be clear from the ensuing arguments that the quantities on the LHS of \eqref{ctrl:Boot} take values continuously in time, from which it follows that $I$ is relatively closed in $\Real_+$. Hence define $T^\star \leq \infty$ such that $I = [0,T^\star]$.
 In order to prove that $T^\star = \infty$ it suffices to establish that $I$ is also relatively open, which is implied by the following bootstrap. 
\begin{proposition}[Bootstrap] \label{lem:Boot}
There exists $\epsilon_0, K_i$  depending only on $d,M,\bar{\lambda},C_0,\kappa,\lambda_0,\lambda^\prime$ and $s$ 
 such that if \eqref{ctrl:Boot} holds on some time interval $[0,T^\star)$ and $\epsilon < \epsilon_0$,  
then for $t \leq T^\star$, 
\begin{subequations} \label{ctrl:BootRes}
\begin{align}
\sum_{\alpha \in \Naturals^d:\abs{\alpha} \leq M}\norm{\jap{\grad_{z,v}}A(v^\alpha f)(t)}^2_{2} = \sum_{\alpha \in \Naturals^d:\abs{\alpha} \leq M}\norm{A^{(1)} v^\alpha f(t)}^2_{2} & < 2K_1\jap{t}^7 \epsilon^2 \label{ctrl:HiLocalizedB} \\
\sum_{\alpha \in \Naturals^d:\abs{\alpha} \leq M}\norm{\jap{\grad_{z,v}}^{-\beta}A(v^\alpha f)(t)}^2_{2} =\sum_{\alpha \in \Naturals^d:\abs{\alpha} \leq M}\norm{A^{(-\beta)} v^\alpha f(t)}^2_{2} & < 2K_2\epsilon^2 \label{ctrl:LowCommLocB} \\
\int_0^{t} \norm{A\rho(\tau)}_2^2 d\tau & < 2K_3\epsilon^2, \label{ctrl:Mid} 
\end{align}
\end{subequations}
from which it follows that $T^\star = \infty$. 
\end{proposition} 
Once Proposition \ref{lem:Boot} is deduced, Theorem \ref{thm:Main} follows quickly. 
This is carried out in \S\ref{sec:ConcProof}.

\begin{remark} 
In order to close the bootstrap in Proposition \ref{lem:Boot} we need to understand how, or if, the constants $K_i$ depend on each other so that we can be sure that they can be chosen self-consistently. 
In fact, $K_3$ is determined by the linearized evolution (from $C_{LD}$ in Lemma \ref{lem:LinearL2Damping}) then $K_1$ is fixed in \S\ref{subsec:HiNorm} depending on $K_3$ (and $s$,$d$,$\lambda_0$,$\lambda^\prime$)
 and $K_2$ is analogously fixed in \S\ref{subsec:LoNorm} depending on $K_3$ (and $s$,$d$,$\lambda_0$,$\lambda^\prime$ but not directly $K_1$). 
Finally, $\epsilon_0$ is chosen small with respect to everything. 
\end{remark}

\begin{remark} \label{rmk:unbalance} 
The unbalance of a whole derivative between \eqref{ctrl:Mid} and \eqref{ctrl:HiLocalizedB} uses the regularization of the interaction potential and is the only aspect of the proof which requires $\gamma \geq 1$.
\end{remark}

The purpose of the weights in velocity is to control derivatives of the Fourier transform so that the trace Lemma \ref{lem:SobTrace} and the $H^{d/2+} \hookrightarrow C^0$ embedding theorem can be applied to restrict the Fourier transform along rays and pointwise.
Both are necessary to deduce the controls on the density in \S\ref{sec:L2I} and \S\ref{sec:PtwiseRho}. 
Moreover, from \eqref{ctrl:LowCommLoc} and \eqref{ineq:MomentDerivs} we have, 
\begin{align}
\norm{A^{(-\beta)}\hat{f}}_{L_k^2H^M_\eta} \lesssim \sum_{\alpha \in \Naturals^d: \abs{\alpha} \leq M}\norm{D^\alpha_\eta A^{(-\beta)}\hat{f}}_{L_k^2 L_\eta^2} & \leq \sum_{\alpha \in \Naturals^d: \abs{\alpha} \leq M}\sum_{j \leq \alpha} \frac{\alpha!}{j!(\alpha-j)!}\norm{(D_\eta^{\alpha-j}A^{(-\beta)})(D_\eta^{j}\hat f)}_{L_k^2 L_\eta^2} \nonumber \\ 
& \lesssim_{M} \sqrt{K_2} \epsilon. \label{ineq:OuterMctrl}
\end{align} 
Similarly, \eqref{ctrl:HiLocalized} implies
\begin{align} 
\norm{A^{(1)}\hat{f}}_{L_k^2 H^M_\eta} \lesssim_{M}  \sqrt{K_1} \epsilon \jap{t}^{7/2}. \label{ineq:OuterActrl}
\end{align} 

Let us briefly summarize how Proposition \ref{lem:Boot} is proved.  
The main step is the $L^2_t$ estimate \eqref{ctrl:Mid}, deduced in \S\ref{sec:L2I}. 
This is done by analyzing \eqref{eq:rhointegral}. 
The linear term in \eqref{eq:rhointegral} is treated using a Fourier-Laplace transform and \textbf{(L)} as in \cite{MouhotVillani11}
with a technical time decomposition in order to get an estimate using $A(t)$. 
The nonlinear term \eqref{eq:rhointegral} is decomposed using a paraproduct into \emph{reaction}, \emph{transport} and \emph{remainder} terms.  
As discussed above in \S\ref{sec:main-ingr-proof}, the reaction term in \eqref{eq:rhointegral} is connected to the plasma echoes. Once the paraproduct decomposition and \eqref{ctrl:LowCommLoc} have allowed us to isolate this effect, it is treated with an adaptation of \S7 in \cite{MouhotVillani11}, carried out in \S\ref{sec:Proofa} (our treatment is in the same spirit but not quite the same). 
The transport term describes the interaction of $\rho$ with `higher' frequencies of $f$; this effect is controlled using \eqref{ctrl:HiLocalized}. 
Once the $L_t^2$ estimate has been established, we also derive a relatively straightforward pointwise-in-time control in \S\ref{sec:PtwiseRho}.  

 The estimate \eqref{ctrl:HiLocalizedB} is deduced in \S\ref{subsec:HiNorm} via an energy estimate in the spirit of \cite{BM13}. 
A paraproduct is again used to decompose the nonlinearity into \emph{reaction}, \emph{transport} and \emph{remainder} terms. 
As in \cite{BM13}, the transport term is treated using an adaptation of the methods of \cite{FoiasTemam89,LevermoreOliver97,KukavicaVicol09}. 
However, perhaps more like \cite{MouhotVillani11}, the reaction term is treated using \eqref{ctrl:MidBoot}.  
The time growth is due to the fact that there is no regularity available to transfer to decay; that the estimate is closable at all requires the regularization from $\gamma \geq 1$. 
The low norm estimate \eqref{ctrl:LowCommLocB} is proved in \S\ref{subsec:LoNorm} in a fashion similar to that of \eqref{ctrl:HiLocalizedB} (the uniform bound is possible due to the regularity gap of $\beta > 2$ derivatives between \eqref{ctrl:MidBoot} and \eqref{ctrl:LowCommLocB}).

\section{Toolbox}\label{sec:toolbox}
\label{Apx:LPProduct}
In this section we review some of the technical tools used in the proof of Theorem \ref{thm:Main}: the Littlewood-Paley dyadic decomposition, paraproducts and useful inequalities for working in Gevrey regularity. 

\subsection{Fourier analysis conventions}
 For a function $f=f(z,v)$ in the Schwartz space, we define its
 Fourier transform $\hat{f}_k(\eta)$ where $(k,\eta) \in \Integer^d
 \times \Real^d$ by
\begin{align*} 
\hat{f}_k(\eta) = \frac{1}{(2\pi)^{d}}\int_{\Torus^d \times \Real^d} e^{-i z k - iv\eta} f(z,v) \ddz \ddv. 
\end{align*}
Similarly we have the Fourier inversion formula, 
\begin{align*} 
f(z,v) = \frac{1}{(2\pi)^{d}}\sum_{k \in \Integer^d} \int_{\Real^d}
e^{i z k + iv\eta} \hat{f}_k(\eta) \, d\eta. 
\end{align*} 
As usual, the Fourier transform and its inverse are extended to
$L^2(\T^d \times \R^d)$ via duality.  We also need to apply the
Fourier transform to functions of $x$ or $v$ alone, for which we use
analogous conventions.  

With these conventions we have the following relations:
\[
\begin{cases}\dss
\int_{\T^d \times \R^d} f(z,v) \overline{g}(z,v) \ddz\ddv  \dss =
\sum_{k \in \Z^d}\int_{\R^d}
\hat{f}_k(\eta) \overline{\hat{g}}_{k}(\eta) \, d\eta, \vsp \\ \dss
\widehat{fg}  = \dss \frac{1}{(2\pi)^{d}}\hat{f} \ast \hat{g},
\vsp \\ \dss
(\widehat{\grad f})_k(\eta) = (ik,i\eta)\hat{f}_k(\eta),
\vsp  \\ \dss
(\widehat{v^\alpha f})_k(\eta) = D_\eta^\alpha \hat{f}_k(\eta). 
\end{cases}
\]

The following versions of Young's inequality occur frequently in the proof. 
\begin{lemma}  
\begin{itemize} 
\item[(a)] Let $f_k(\eta),g_k(\eta) \in L^2(\Integer^d \times \Real^d)$ and $\jap{k}^\sigma h_k(t) \in L^2(\Integer^d)$ for $\sigma > d/2$. Then, for any $t \in \Real$, 
\begin{align} 
\abs{\sum_{k,\ell} \int_\eta  f_k(\eta) h_\ell(t) g_{k-\ell}(\eta-t\ell) \, d\eta} \lesssim_{d,\sigma} \norm{f}_{L^2_{k,\eta} }\norm{g}_{L^2_{k,\eta} } \norm{\jap{k}^\sigma h(t)}_{L_k^2}. \label{ineq:L2L2L1}      
\end{align} 
\item[(b)] Let $f_k(\eta), \jap{k}^\sigma g_k(\eta) \in L^2(\Integer^d \times \Real^d)$ and $h_k \in L^2(\Integer^d)$ for $\sigma > d/2$. Then, for any $t \in \Real$, 
\begin{align} 
\abs{\sum_{k,\ell} \int_\eta  f_k(\eta) h_\ell(t) g_{k-\ell}(\eta-t\ell) \, d\eta} \lesssim_{d,\sigma} \norm{f}_{L^2_{k,\eta} }\norm{\jap{k}^{\sigma}g }_{L^2_{k,\eta} } \norm{h(t)}_{L_k^2}.    \label{ineq:L2L1L2}      
\end{align} 
\end{itemize}
\end{lemma} 
\begin{proof} 
To prove (a): 
\begin{align*} 
\abs{\sum_{k,\ell} \int_\eta  f_k(\eta) h_\ell g_{k-\ell}(\eta-t\ell) \, d\eta} & \leq \sum_{k} \left( \int_{\eta} \abs{f_k(\eta)}^2 \, d\eta\right)^{1/2} \sum_{\ell} h_\ell \left(\int_{\eta} \abs{g_{k-\ell}(\eta-\ell t)}^2 \, d\eta\right)^{1/2} \\ 
& = \sum_{k} \left( \int_{\eta} \abs{f_k(\eta)}^2 \, d\eta\right)^{1/2} \sum_{\ell} h_\ell \left(\int_{\eta} \abs{g_{k-\ell}(\eta)}^2 \, d\eta\right)^{1/2} \\ 
& \leq  \left( \sum_{k}\int_{\eta}\abs{f_k(\eta)}^2 \, d\eta \right)^{1/2} \left[\sum_{k} \left(\sum_{\ell} h_\ell \left(\int_{\eta} \abs{g_{k-\ell}(\eta)}^2 \, d\eta\right)^{1/2}\right)^{2}\right]^{1/2} \\ 
& \leq \norm{f}_{L^2_{k,\eta}} \norm{g}_{L^2_{k,\eta}} \sum_{\ell} \abs{h_{\ell}} \\ 
& \lesssim_{d,\sigma} \norm{f}_{L^2_{k,\eta}} \norm{g}_{L^2_{k,\eta}} \norm{\jap{k}^\sigma h_k}_{L_k^2}, 
\end{align*} 
where the penultimate line followed from the $L^2 \times L^1 \mapsto L^{2}$ Young's inequality. 
The proof of (b) is analogous, simply putting $g$ rather than $h$ in $L^1_k$. 
\end{proof} 

\subsection{Littlewood-Paley decomposition}
\label{sec:littl-paley-decomp}

This work makes use of the Littlewood-Paley dyadic decomposition.
Here we fix conventions and review the basic properties of this
classical theory (see e.g. \cite{BCD11} for more details).

Define the joint variable $\Xi := (k,\eta) \in \Z^d \times \R^d$. 
Let $\psi \in C_0^\infty(\Z^d \times \Real^d;\Real)$ be a radially symmetric
non-negative function such that $\psi(\Xi) = 1$ for $\abs{\Xi} \leq
1/2$ and $\psi(\Xi) = 0$ for $\abs{\Xi} \geq 3/4$. 
Then we define $\phi(\Xi) := \psi(\Xi/2) - \psi(\Xi)$, a 
non-negative, radially symmetric function supported in the annulus $\{ 1/2 \le |\Xi| \le 3/2 \}$. 
Then we define the rescaled functions $\phi_N(\Xi) =
\phi(N^{-1}\Xi)$, which satisfy 
\begin{equation*} 
\textup{supp}\, \phi_N(\Xi) = \{ N/2 \leq \abs{\Xi} \leq
3N/2 \}
\end{equation*}
and we have classically the partition of unity,
\begin{equation*} 
1 = \psi(\Xi) + \sum_{N \in 2^\Naturals} \phi_N(\Xi), 
\end{equation*}  
(observe that there is always at most two non-zero terms in this sum),
where we mean that the sum runs over the dyadic numbers $M =
1,2,4,8,\dots,2^{j},\dots$.

For $f \in L^2(\Torus^d \times \Real^d)$ we define
\begin{equation*}
  \begin{cases} \dss
    f_{N} & := \phi_N\left(\grad_{z,v}\right)f, \vsp \\ 
    f_{\frac{1}{2}} & :=  \psi\left(\grad_{z,v}\right)f, \vsp \\ 
    f_{< N} & := \dss f_{\frac{1}{2}} + \sum_{N' \in 2^{\Naturals} \,
      : \, N' < N}
    f_{N'}, 
 \end{cases}
\end{equation*}
and we have the natural decomposition  (recall $\mathbb D$ is defined in \S\ref{sec:Notation}), 
\begin{equation*}
  f = \sum_{N \in \mathbb D} f_N = f_{\frac{1}{2}} + \sum_{N \in 2^\Naturals} f_N. 
\end{equation*}
Normally one would use $f_0$ or $f_{-1}$ rather than the slightly
inconsistent $f_{1/2}$, however $f_0$ is reserved for the much
more commonly used projection onto the zero mode only in $z$ or $x$. 

There holds the almost orthogonality and the approximate projection
property
\begin{equation} \label{ineq:LPOrthoProject}
\begin{cases} \dss
\sum_{N \in \mathbb D} \norm{f_N}_2^2 \le \norm{f}^2_2 \le 2 \sum_{N \in \mathbb D} \norm{f_N}_2^2 \vsp \\ \dss
\norm{(f_{N})_N}_2 \le \norm{f_N}_2. 
\end{cases}
\end{equation}
Similar to \eqref{ineq:LPOrthoProject} but more generally, if $f =
\sum_{\mathbb D} g_N$ with $\mbox{supp} \, g_N \subset \{
C^{-1}N \le |\Xi| \le C N\}$ for $N \geq 1$ and $\mbox{supp} \, g_{\frac{1}{2}} \subset \set{\abs{\Xi} \leq C}$ then we have
\begin{align} 
  \norm{f}^2_2 \approx_C \sum_{N \in \mathbb D}
  \norm{g_N}_2^2. \label{ineq:GeneralOrtho}
\end{align}
Moreover, the dyadic decomposition behaves nicely with respect to differentiation:
\begin{align} 
\norm{\jap{\grad_{z,v}}f_N}_2 \approx N \norm{f_N}_2. \label{ineq:dyadicDiff}
\end{align}

We also define the notation
\begin{align*} 
f_{\sim N} = \sum_{N' \in \mathbb{D} \, : \, C^{-1} N \leq N' \leq
  C N} f_{N'}, 
\end{align*}
for some constant $C$ which is independent of $N$.  Generally the
exact value of $C$ which is being used is not important; what is
important is that it is finite and independent of $N$.

During some steps of the proof we will apply the Littlewood-Paley
decomposition to the spatial density $\rho(t,x)$.  In this case it will be
more convenient to use the following definition that uses $kt$ in place of the frequency in $v$, a natural convention when one recalls that $\hat{\rho}_k(t) = \hat{f}_k(t,kt)$: 
 \begin{equation*}
\begin{cases} 
  \widehat{\rho(t)}_N & = \phi_N(\abs{k,kt})\hat{\rho}_k(t), \vsp \\
  \widehat{\rho(t)}_{\frac{1}{2}} & = \psi(\abs{k,kt})\hat{\rho}_k(t),
  \vsp \\
  \widehat{\rho(t)}_{< N} & = \dss \rho(t)_{\frac{1}{2}} + \sum_{N' \in
    2^{\Naturals}: N' < N} \rho(t)_{N'}.
\end{cases}
\end{equation*}

\begin{remark}
  We have opted to use the compact notation above, rather
  than the commonly used alternatives $\Delta_{j}f = f_{2^j}$ and
  $S_jf = f_{<2^j}$, in order to reduce the number of characters in
  long formulas.
\end{remark}

\subsection{The paraproduct decomposition}
\label{sec:parapr-decomp}
Another key Fourier analysis tool employed in this work is the
\emph{paraproduct decomposition}, introduced by Bony \cite{Bony81} (see also
\cite{BCD11}).  Given suitable functions $f,g$ we may define the
paraproduct decomposition (in either $(z,v)$ or just $v$),
\begin{equation}\label{eq:paraproduct} 
fg = \sum_{N \geq 8} f_{<N/8}g_N + \sum_{N \geq 8} f_N  g_{<N/8} +
\sum_{N \in \mathbb{D}}\sum_{N/8 \leq N^\prime \leq 8N}
f_{N} g_{N^\prime} := T_fg + T_gf + \mathcal{R}(f,g)
\end{equation}
where all the sums are understood to run over $\mathbb D$. 

The advantage of this decomposition in the energy estimates is that
the first term $T_f g$ contains the highest derivatives on $g$ but
allows to take advantage of the frequency cutoff on the first function
$f$, whereas the second term $T_g f$ contains the highest derivatives
on $f$ but allows to take advantage of the frequency cutoff on the
second function $g$. 
Finally the last ``remainder'' term contains the contribution from comparable frequencies which allows to split regularity evenly between the factors. 

\subsection{Elementary inequalities for Gevrey regularity} \label{apx:Gev}
In this section we discuss a set of elementary, but crucial, inequalities for working in Gevrey regularity spaces. 
First we point out that Gevrey and Sobolev regularities can be related with the following two inequalities. 
\begin{itemize}
\item[(i)] For $x \geq 0$, $\alpha > \beta \geq 0$, $C,\delta > 0$, 
\begin{align} 
e^{Cx^{\beta}} \leq e^{C\left(\frac{C}{\delta}\right)^{\frac{\beta}{\alpha - \beta}}} e^{\delta x^{\alpha}};  \label{ineq:IncExp}
\end{align}
\item[(ii)] For $x \geq 0$, $\alpha > 0$, $\sigma,\delta > 0$, 
\begin{align} 
e^{-\delta x^{\alpha}} \lesssim \frac{1}{\delta^{\frac{\sigma}{\alpha}} \jap{x}^{\sigma}}. \label{ineq:SobExp}
\end{align}
\end{itemize}

Next, we state several useful inequalities regarding the weight $\jap{x} = (1 + \abs{x}^2)^{1/2}$. 
In particular, the improvements to the triangle inequality for $s < 1$ given in \eqref{ineq:TrivDiff2}, \eqref{lem:scon2} and \eqref{lem:strivial2} are important for getting useful bilinear (and trilinear) estimates. The proof is straightforward and is omitted. 
\begin{lemma} \label{lem:jap}
Let $0 < s < 1$ and $x,y \geq 0$.
\begin{itemize}
\item[(i)] We have the triangle inequalities (which hold also for $s = 1$), 
\begin{subequations} \label{ineq:JapTri}
\begin{align} 
\jap{x + y}^s & \leq \jap{x}^s + \jap{y}^s \label{ineq:JapTri1} \\ 
\abs{\jap{x}^s - \jap{y}^s} & \leq \jap{x-y}^s \label{ineq:JapTri2} \\ 
C_s(\jap{x}^s + \jap{y}^s) & \leq \jap{x+y}^s,  \label{ineq:Comps}
\end{align} 
for some $C_s > 0$ depending only on $s$. 
\end{subequations} 
\item[(ii)] In general, 
\begin{align} 
\abs{\jap{x}^s - \jap{y}^s} \lesssim_s \frac{1}{\jap{x}^{1-s} + \jap{y}^{1-s}}\jap{x-y}. \label{ineq:TrivDiff2}
\end{align}
\item[(iii)] If $\abs{x-y} \leq x/K$ for some $K > 1$, then we have the
  improved triangle inequality
\begin{align} 
\abs{\jap{x}^s - \jap{y}^s} \leq \frac{s}{(K-1)^{1-s}}\jap{x-y}^s. \label{lem:scon2}
\end{align} 
\item[(iv)] We have the improved triangle inequality for $x \geq y$, 
\begin{align} 
\jap{x + y}^s \leq \left(\frac{\jap{x}}{\jap{x+y}}\right)^{1-s}\left(\jap{x}^s + \jap{y}^s\right). \label{lem:strivial2}
\end{align} 
\end{itemize}
\end{lemma} 
The following product lemma will be used several times in the sequel.
Notice that $\tilde c < 1$ for $s < 1$, which shows that we gain something by working in Gevrey spaces with $s < 1$; indeed this important gain is used many times in the nonlinear estimates. 
We sketch the proof as it provides a representative example of arguments used several times in \S\ref{sec:Energy}.
\begin{lemma}[Product Lemma] \label{lem:GevProdAlg}
For all $0<s<1$, and $\sigma \geq 0$ there exists a $\tilde c = \tilde{c}(s,\sigma) \in (0,1)$ such that the following holds 
for all $\lambda> 0$, $f,\tilde{f} \in \mathcal{G}^{\lambda,\sigma;s}(\T^d \times \Real^d)$ and $h(t) \in \mathcal{F}^{\lambda,\sigma;s}(\T^d)$:  
\begin{subequations}
\begin{align} 
\sum_{k \in \Integer^d} \sum_{\ell \in \Integer^d} \int_{\R^d} \jap{k,\eta}^{2\sigma} e^{2\lambda\jap{k,\eta}^s} \abs{\hat{f}_k(\eta)\hat{h}_{\ell}(t)\hat{\tilde{f}}_{k-\ell}(\eta-\ell t)} d\eta  & \nonumber \\ & \hspace{-7cm} \lesssim_{\lambda,\sigma,s,d} \norm{f}_{\cG^{\lambda,\sigma;s}}\norm{\tilde f}_{\cG^{\tilde{c}\lambda,0;s}}\norm{h(t)}_{\cF^{\lambda,\sigma;s}}
 + \norm{f}_{\cG^{\lambda,\sigma;s}}\norm{\tilde f}_{\cG^{\lambda,\sigma;s}} \norm{h(t)}_{\cF^{\tilde{c}\lambda,0;s}}. \label{ineq:GProduct2}
\end{align}
\end{subequations}
Moreover, we have the algebra property  
\begin{align} 
\sum_{k \in \Integer^d} \int_{\R^d}\abs{\sum_{\ell \in \Integer^d} e^{\lambda \jap{k,\eta}^s}\jap{k,\eta}^\sigma \hat{h}_{\ell}(t)\hat{f}_{k-\ell}(\eta - \ell t) }^2 d\eta 
 \lesssim  \norm{h(t)}^2_{\cF^{\lambda,\sigma;s}} \norm{f}^2_{\cG^{\lambda,\sigma;s}}. \label{ineq:GAlg} 
\end{align} 
\end{lemma}
\begin{proof} 
We only prove \eqref{ineq:GProduct2}, which is slightly harder. 
Denote $B_k(\eta) = \jap{k,\eta}^{\sigma} e^{\lambda\jap{k,\eta}^s}$. 
The proof proceeds by decomposing with a paraproduct: 
\begin{align*}
\sum_{k,\ell \in \Integer^d} \int_\eta \abs{B\hat{f}_k(\eta)B_k(\eta)\hat{h}_{\ell}(t)\hat{\tilde{f}}_{k-\ell}(\eta-\ell t)} d\eta &\leq  \sum_{N \geq 8}\sum_{k,\ell \in \Integer_\ast^d} \int_\eta \abs{B\hat{f}_k(\eta)B_k(\eta)\hat{h}_{\ell}(t)_N\hat{\tilde{f}}_{k-\ell}(\eta-\ell t)_{<N/8}} d\eta \\ 
& \hspace{-1cm} + \sum_{N \geq 8}\sum_{k,\ell \in \Integer_\ast^d} \int_\eta \abs{B\hat{f}_k(\eta)B_k(\eta)\hat{h}_{\ell}(t)_{<N/8}\hat{\tilde{f}}_{k-l}(\eta-\ell t)_{N}} d\eta \\ 
& \hspace{-1cm} + \sum_{N \in \mathbb D} \sum_{N/8 \leq N^\prime \leq 8N}\sum_{k,\ell \in \Integer_\ast^d} \int_\eta \abs{B\hat{f}_k(\eta)B_k(\eta)\hat{h}_{\ell}(t)_{N^\prime}\hat{\tilde{f}}_{k-\ell}(\eta-\ell t)_{N}} d\eta \\ 
& = R + T + \mathcal{R}. 
\end{align*}
Note that the $R$ and $T$ term are almost, but not quite, symmetric. 
Consider first the $R$ term. On the support of the integrand we have the frequency localizations 
\begin{subequations} \label{ineq:FreqLocExample}
\begin{align} 
\frac{N}{2} \leq \abs{\ell,\ell t} & \leq \frac{3N}{2}, \\ 
\abs{k-\ell,\eta - \ell t} & \leq \frac{3N}{32}, \\ 
\frac{13}{16} \leq \frac{\abs{k,\eta}}{\abs{\ell,t \ell}} & \leq \frac{19}{16}.  
\end{align}  
\end{subequations}
Therefore, by \eqref{lem:scon2}, there exists some $c = c(s) \in (0,1)$ such that
\begin{align*} 
B_k(\eta) \leq \jap{k,\eta}^{\sigma} e^{\lambda\jap{\ell,\ell t}^s} e^{c\lambda\jap{k-\ell,\eta -\ell t}^s} & \lesssim_\sigma  \jap{\ell,t\ell}^{\sigma} e^{\lambda\jap{\ell,\ell t}^s} e^{c\lambda\jap{k-\ell,\eta -\ell t}^s} \\ 
& \lesssim_{\lambda}  \jap{\ell,t\ell}^{\sigma} e^{\lambda\jap{\ell,\ell t}^s}\jap{k-\ell,\eta -\ell t}^{-\frac{d}{2}-1} e^{\frac{1}{2}(c+1)\lambda\jap{k-\ell,\eta -\ell t}^s}, 
\end{align*}
where in the last inequality we applied \eqref{ineq:SobExp}. 
Adding a frequency localization with \eqref{ineq:FreqLocExample}, denoting $\tilde c = \frac{1}{2}(c + 1)$, and using \eqref{ineq:L2L1L2} we have      
\begin{align*} 
R &  \lesssim \sum_{N \geq 8}\sum_{k,\ell \in \Integer_\ast^d} \int_\eta \abs{B\hat{f}_k(\eta)_{\sim N}B_\ell(\ell t)\hat{h}_{\ell}(t)_Ne^{\tilde c\lambda\jap{k-\ell,\eta -\ell t}^s} \jap{k-\ell,\eta-\ell t}^{-\frac{d}{2}-1} \hat{\tilde{f}}_{k-l}(\eta-lt)_{<N/8}} d\eta \\ 
& \lesssim \sum_{N \geq 8} \norm{f_{\sim N}}_{\cG^{\lambda,\sigma;s}} \norm{h(t)_N}_{\cF^{\lambda,\sigma;s}} \norm{\tilde f}_{\cF^{\tilde c\lambda,0;s}}. 
\end{align*} 
By Cauchy-Schwarz and \eqref{ineq:GeneralOrtho} we have 
\begin{align*} 
R \lesssim \norm{f}_{\cG^{\lambda,\sigma;s}} \norm{h(t)}_{\cF^{\lambda,\sigma;s}} \norm{\tilde f}_{\cF^{\tilde c\lambda,0;s}},  
\end{align*}
which appears on the RHS of \eqref{ineq:GProduct2}. 
Treating the $T$ term is essentially the same except reversing the role of $(\ell,t\ell)$ and $(k-\ell,\eta-t\ell)$ and applying \eqref{ineq:L2L2L1} as opposed to \eqref{ineq:L2L1L2}.  

To treat the $\mathcal{R}$ term we use a simple variant. 
We claim that there exists some $c^\prime = c^\prime(s) \in (0,1)$ such that on the support of the integrand we have
\begin{align} 
B_k(\eta) \lesssim_\sigma e^{c^\prime \lambda \jap{k-\ell,\eta-t\ell}^s}e^{c^\prime \lambda \jap{\ell,t\ell}^s}. \label{ineq:BgainProduct}
\end{align} 
To see this, consider separately the cases (say) $N \geq 128$ and $N < 128$. On the latter, $B_k(\eta)$ is simply bounded by a constant. 
In the case $N \geq 128$ we have the frequency localizations
\begin{subequations} \label{ineq:RRNL_FreqLocI}
\begin{align} 
\frac{N}{2} \leq \abs{k-\ell,\eta - \ell t} & \leq \frac{3N}{2}, \\ 
\frac{N^\prime}{2} \leq \abs{\ell,\ell t} & \leq \frac{3N^\prime}{2}, \\ 
\frac{1}{24} \leq \frac{\abs{k-\ell,\eta - kt}}{\abs{\ell,\ell t}} & \leq 24, 
\end{align}
\end{subequations}
and hence we may apply \eqref{lem:strivial2} since in this case $64 \leq \abs{k-\ell,\eta-\ell\tau} \approx \abs{\ell,\ell\tau}$. 
Further, we can use \eqref{ineq:SobExp} to absorb the Sobolev corrections and, indeed, we have \eqref{ineq:BgainProduct} on the support of the integrand in $\mathcal{R}$.  
Hence, 
\begin{align*}
\mathcal{R} \lesssim \sum_{N \in \mathbb D} \sum_{N^\prime \approx N}\sum_{k,\ell \in \Integer_\ast^d} \int_\eta \abs{B\hat{f}_k(\eta)e^{c^\prime \lambda \jap{\ell,\ell t}^s}\hat{h}_{\ell}(t)_{N^\prime} e^{c^\prime \lambda \jap{k-\ell,\eta - \ell t}^s}\hat{\tilde{f}}_{k-l}(\eta-lt)_{N}} d\eta. 
\end{align*} 
Applying \eqref{ineq:L2L1L2} followed by \eqref{ineq:dyadicDiff} and \eqref{ineq:SobExp} (since $c^\prime < 1$) implies, 
\begin{align*} 
\mathcal{R} & \lesssim \sum_{N \in \mathbb D} \sum_{N^\prime \approx N} \norm{f}_{\cG^{\lambda,\sigma;s}} \norm{h_{N^\prime}}_{\cF^{c^\prime \lambda,\frac{d}{2}+1;s}} \norm{\tilde f_N}_{\cG^{c^\prime \lambda,0;s}} \\ 
& \lesssim \sum_{N \in \mathbb D} \frac{1}{N} \norm{f}_{\cG^{\lambda,\sigma;s}} \norm{h_{\sim N}}_{\cF^{c^\prime \lambda,\frac{d}{2}+2;s}} \norm{\tilde f_N}_{\cG^{c^\prime \lambda,0;s}} \\ 
& \lesssim_{\lambda,\sigma} \norm{f}_{\cG^{\lambda,\sigma;s}} \norm{h}_{\cF^{\lambda,\sigma;s}} \norm{\tilde f}_{\cG^{c^\prime \lambda,0;s}}. 
\end{align*}
Hence (after possibly adjusting $\tilde c$) this term appears on the RHS of \eqref{ineq:GProduct2}. 
\end{proof}
We also need the standard Sobolev space trace lemma, which we will apply on the Fourier side. 
\begin{lemma}[$L^2$ Trace] \label{lem:SobTrace}
Let $f$ be smooth and $C \subset \Real^d$ be an arbitrary straight line. Then for all $\sigma \in \Real_+$ with $\sigma > (d-1)/2$, 
\begin{align*} 
\norm{f}_{L^2(C)} \lesssim \norm{f}_{H^\sigma}. 
\end{align*} 
\end{lemma} 
\begin{proof} 
Follows by induction on co-dimension with the standard $H^{1/2}$ restriction theorem \cite{Adams03}.
\end{proof}

\section{Linear damping in Gevrey regularity} \label{sec:LinDamp}
The first step to proving \eqref{ctrl:Mid} is understanding (forced) linear Landau damping in the $L^2$ Gevrey norms we are using.
This is provided by the following lemma, which also shows that \textbf{(L)} implies linear damping in all Gevrey regularities ($s > 1/3$ is not relevant to the proof). 
It is crucial that the same norm appears on both sides of \eqref{ineq:LinearCtrl} so that we may use it in the nonlinear estimate on $\rho(t)$. 
The main idea of the proof of Lemma \ref{lem:LinearL2Damping} appears in \cite{MouhotVillani11} to treat damping in Gevrey regularity and is based on decomposing the solution into analytic/exponentially decaying sub-components.
We note that Lin and Zeng in \cite{LZ11b} have linear damping results at much lower regularities (similarly, we believe the following proof can be adapted also to Sobolev spaces). 
We will first give a formal proof of Lemma \ref{lem:LinearL2Damping} as an a priori estimate in \S\ref{sec:proof-priori-estim}
and then explain the rigorous justification in \S\ref{sec:rigor-just-priori}. 
In \S\ref{sec:penrose-criterion}, we discuss the proof that the Penrose condition \cite{Penrose} implies \textbf{(L)}.  

\begin{lemma}[Linear integral-in-time
  control] \label{lem:LinearL2Damping} 
  Let $f^0$
  satisfy the condition \textbf{(L)} with constants $M> d/2$ and $C_0,\bar \lambda,\kappa > 0$. 
  Let 
  $A_k(t,\eta)$ 
  be the multiplier defined in \eqref{eq:Akteta} for $s \in (0,1)$,
  and $\lambda=\lambda(t) \in (0,\lambda_0)$ as defined in
  \eqref{def:lambda}.
   Let $F(t)$ and $T^\star > 0$ be given such that, if we denote $I = [0,T^\star)$, 
  \[
  \int_0 ^{T_*} \norm{F(t)}_{\cF^{\lambda(t),\sigma;s}} ^2 \, dt = \sum_{k \in
    \Z^d_*}\norm{A_k(t,kt) F_k(t)}_{L_t^2(I)}^2 < \infty.
\]
 Then there exists a constant $C_{LD} =
C_{LD}(C_0,\bar{\lambda},\kappa,s,d,\lambda_0,\lambda^\prime)$ such
that for all $k \in \Integer_\ast^d$, the solution $\phi_k(t)$ to the
system
\begin{align}\label{eq:volterra}
  \phi_k(t)= F_k(t) + \int_0^tK^0(t-\tau,k) \phi_k(\tau) \, d\tau
\end{align}
in $t \in \R_+$ with $K^0(t,k) := -\tilde f^0\left(kt\right)
\widehat{W}(k)\abs{k}^2t$ satisfies the mode-by-mode estimate
\begin{align}
\forall \, k \in \Integer_\ast^d, \quad
 \int_0 ^{T_*} A_k(t,kt)^2 |\phi_k(t)|^2 \, dt \leq C_{LD} ^2 \int_0 ^{T_*} A_k(t,kt)
|F_k(t))|^2 \, dt \label{ineq:LinearCtrl}
\end{align} 
which is equivalent to $\int_0 ^{T_*}
\norm{\phi}_{\cF^{\lambda(t),\sigma;s}}^2 \, dt \le C_{LD} ^2 \int_0
^{T_*} \norm{F}_{\cF^{\lambda(t),\sigma;s}}^2 \, dt$. 
\end{lemma}

\begin{remark} 
The proof proceeds slightly differently in the case $s = 1$ where the additional requirement $\bar\lambda > \lambda(0)$ occurs naturally (and the constant depends badly in the parameter $\bar\lambda - \lambda(0)$). 
\end{remark}

\subsection{Proof of the a priori estimate}
\label{sec:proof-priori-estim}

We only consider the $s < 1$ case; the analytic case is only a slight variation.  
As the hypothesis on $F$ is known a priori only to hold on $[0,T_*)$, we 
simply extend $F_k(t)$ to be zero for all $t \geq T_*$. 
\mk

\noindent {\em Step 1. Rough Gr\"onwall bound.} First we deduce a
rough bound using Gr\"onwall's inequality with no attempt to be
optimal. This bound shows in particular that the integral equation
\eqref{eq:volterra} is globally well-posed (for each frequency $k \in
\Z^d_*$) in the norm associated with the multiplier $A$. By
\eqref{ineq:JapTri}, the definition of $K^0$, \eqref{ineq:Wbd} and
that $\lambda(t)$ is non-increasing in time,
\begin{align*}
  A_k(t,kt) \abs{\phi_k(t)} & \leq A_k(t,kt)\abs{F_k(t)} + \int_0^t A_k(t,kt)\abs{K^0(t-\tau,k)}\abs{\phi_k(\tau)} \, d\tau \\
  & \lesssim A_k(t,kt)\abs{F_k(t)} + \int_0^t \jap{k(t-\tau)}^\sigma
  e^{\lambda(t)\jap{k(t-\tau)}^s} \abs{\widehat{\grad_v
      f^0}(k(t-\tau))}A_k(\tau,k\tau)\abs{\phi_k(\tau)} \, d\tau.
\end{align*} 
Then by \eqref{ineq:SobExp}, the $H^{d/2+} \hookrightarrow C^0$ embedding theorem and \eqref{ineq:f0Loc}, 
\begin{align*}
  A_k(t,kt) \abs{\phi_k(t)} & \lesssim A_k(t,kt)\abs{F_k(t)} + \left(\sup_{\eta \in \Real^d}\jap{\eta}^\sigma e^{\lambda(0)\jap{\eta}^s} \abs{\widehat{f^0}(\eta)}\right) \int_0^t A_k(\tau,k\tau)\abs{\phi_k(\tau)} \, d\tau \\
  & \lesssim A_k(t,kt)\abs{F_k(t)} + \norm{\jap{\eta}^\sigma
    e^{\lambda(0)\jap{\eta}^s} \widehat{f^0}(\eta)}_{H^M_\eta} \int_0^t
  A_k(\tau,k\tau)\abs{\phi_k(\tau)} \, d\tau.
\end{align*} 
By \eqref{ineq:f0Loc} with \eqref{ineq:SobExp}, it follows by Gr\"onwall's inequality that for some $C > 0$ we have the following (using also Cauchy-Schwarz and \eqref{ineq:SobExp} in the last inequality),  
\begin{align} 
A_k(t,kt) \abs{\phi_k(t)} \lesssim e^{Ct} \int_0^t \abs{A_k(\tau,k\tau)F_k(\tau)} \, d\tau \lesssim e^{2Ct}\norm{AF_k}_{L^2_t(I)}. \label{ineq:sorttimeLinear}
\end{align} 
\mk

\noindent {\em Step 2. Frequency localization.}
We would like to use the Fourier-Laplace transform as in
\cite{MouhotVillani11}, but $F_k(t)$ does not decay exponentially in
time.  Instead, we deduce the estimate by decomposing the problem into
a countable number of exponentially decaying contributions.  
Let $R \geq e$ to be a constant fixed later depending only on
$\bar{\lambda}$, $\lambda(0)^{-1}$ and
$s$, let
$F^n_k(t) = F_k(t) \mathbf{1}_{Rn \leq \abs{kt}^s \leq R(n+1)}$, and
define $\phi_k^n$ as solutions to
\begin{align} 
\phi_k^n(t) = F_k^n(t) + \int_0^t K^0(t-\tau,k) \phi_k^n(\tau) \, d\tau. \label{eq:linbyn}
\end{align}
Then $\phi_k(t) = \sum_{n=0}^\infty \phi_k^n(t)$ by linearity of the
equation, and by the definition of $F^n_k(t)$, $\phi_k^n(t)$ is
supported for $\abs{kt}^{s} \geq Rn$. Moreover, obviously
\eqref{ineq:sorttimeLinear} holds for each $\phi_k^n$.  Now we will
use the Fourier-Laplace transform in time to get an $L_t^2$ estimate
on $\phi_k^n$ as in \cite{MouhotVillani11}, but with a contour which
gets progressively closer to the imaginary axis as $n$ increases.
Define the indices $\mu_n$ as
\begin{align*} 
\mu_0 & = 
\mu_1,  \\ 
\mu_n & = \frac{1}{(Rn)^{1/s}}\left[ \lambda\left(\frac{(Rn)^{1/s}}{\abs{k}}\right)\jap{(Rn)^{1/s}}^s + \sigma \log\jap{(Rn)^{1/s}}\right], \quad\quad n \geq 1. 
\end{align*}
We will have the requirement that $\mu_n < \bar{\lambda}$ so that our
integration contour always lies in the half plane defined in
\textbf{(L)}.  As long as $s < 1$ this only requires us to choose $R$
large relative to $\bar{\lambda}(\lambda(0))^{-1}$, so to fix ideas
suppose that $R$ is large enough so that $\sup\mu_n <
\bar{\lambda}/2$.
Define the amplitude corrections
\begin{align*} 
N_{k,0} & = \jap{k}^{\sigma} \exp \left[\lambda\left(\frac{(R)^{1/s}}{\abs{k}}\right)\jap{k}^\sigma \right], \\ 
N_{k,n}& = \frac{\jap{k,(Rn)^{1/s}}^{\sigma}}{\jap{(Rn)^{1/s}}^{\sigma}} \exp \left[\lambda\left(\frac{(Rn)^{1/s}}{\abs{k}}\right)\left[\jap{k,(Rn)^{1/s}}^s - \jap{(Rn)^{1/s}}^s\right] \right] \quad\quad n \geq 1. 
\end{align*}
The role of this correction is highlighted by the fact that when $\abs{kt} = (Rn)^{1/s}$ and $n \geq 1$, then 
\begin{equation*}
  \begin{cases}\displaystyle
  e^{\mu_n|kt|} = \langle kt \rangle^{\sigma} e^{\lambda(t)
    \langle kt \rangle^s} \vspace{0.2cm} \\
  \displaystyle
  N_{n,k} = \frac{\langle k,kt \rangle^\sigma}{\langle kt \rangle^\sigma} e^{\lambda(t)\left( \langle k, kt \rangle^s - \langle kt \rangle^s \right)} \vspace{0.2cm} \\
  \displaystyle
e^{\mu_n|kt|} N_{n,k} = \langle k, kt \rangle^\sigma e^{\lambda(t)
    \langle k, kt \rangle^s} = A_k(t,kt). 
\end{cases}
\end{equation*}
Hence the index $\mu_n$ is related to the radius of analyticity in the velocity
variable, whereas the correction $N_{k,n}$ measures the ratio of what
is lost by not taking into account the regularity in the space
variable. A further error is introduced
by the fact that these regularity weights only exactly fit the
multiplier $A$ at the left endpoint of the interval of the decomposition. 
\mk

\noindent {\em Step 3. The one-block estimate via the Fourier-Laplace
  transform.} Now multiply
\eqref{eq:linbyn} by $e^{\mu_n\abs{k}t}N_{k,n}$, and denoting
\begin{equation*} 
  \begin{cases}
    R^n_k(t) & = e^{\mu_n\abs{k}t} N_{k,n} F_k^n(t) \vsp \\
    \Phi^n_k(t) & = e^{\mu_n\abs{k}t}N_{k,n} \phi_k^n(t),
  \end{cases}
\end{equation*} 
we have
\begin{align} 
\Phi_k^n(t) = R_k^n(t) + \int_0^t e^{\mu_n\abs{k}(t-\tau)} K^0(t-\tau,k) \Phi_k^n (\tau)\, d\tau. \label{def:PhikVolt}
\end{align} 
Taking the Fourier transform in time $\hat G(\omega) = (1/\sqrt{2\pi})
\int_\R e^{-it\omega} G(t) \, dt$ (extending $R_k^n$ and
$\Phi_k^n$ 
and $K^0$ as zero for negative times)\footnote{The tacit assumption that the Fourier transform of $\Phi_k^n(t)$ exists is where the a priori estimate is not rigorous. We justify that the transform exists and that this computation can be made rigorous in \S\ref{sec:rigor-just-priori} below.}
  we obtain
  \begin{equation*}
    \hat \Phi_k ^n(\omega) = \hat R^n_k(\omega) + \mathcal
    L\left(k,\mu_n + i \frac{\omega}{|k|} \right) \hat \Phi^n
    _k(\omega),  
  \end{equation*}
  (where $\mathcal{L}(\xi,k) := \int_0^{+\infty}
  e^{\bar{\xi}\abs{k}t}K^0(t,k) \, dt$), which is formally solved as
\begin{align*} 
  \Phi_k^n(\omega) = \frac{\hat R^n_k(\omega)}{1 - \mathcal
    L\left(k,\mu_n + i \frac{\omega}{|k|} \right)}.
\end{align*} 
Applying the stability condition \textbf{(L)} and Plancherel's theorem implies 
\begin{align} 
\norm{\Phi_k^n(t)}_{L_t^2(\R)} \lesssim \frac{1}{\kappa}\norm{R_k^n(t)}_{L^2_t(\R)}. \label{ineq:PhiBound}
\end{align} 
\mk 

\noindent {\em Step 4. Coming back to the original multiplier $A$.}
Consider $Rn \leq \abs{kt}^{s} \leq R(n+1)$ for $n \geq 1$. From the definition \eqref{def:lambda} of $\lambda(t)$ we see that if $R$ is chosen sufficiently large then for $t \geq R^{1/s}$ we have 
\begin{align*} 
\frac{d}{dt} \left(\lambda(t) \jap{k,kt}^s \right) = \left(\dot{\lambda}(t) + s\lambda(t) \frac{\abs{k}\abs{k,kt}}{\jap{k,kt}^2} \right)\jap{k,kt}^s  = \left(-\frac{a(\lambda_0-\lambda^\prime)}{4\jap{t}^{1+a}} + s\abs{k}\lambda(t)\frac{\abs{k,kt}}{\jap{k,kt}^{2}}\right)\jap{k,kt}^s > 0, 
\end{align*}
and hence $\lambda(t)\jap{k,kt}^s$ is increasing since $\abs{k} \geq 1$. 
Applying this gives, 
\begin{align} 
N_{k,n} e^{\mu_n\abs{kt}}  & = \exp\left[\lambda\left(\frac{(Rn)^{1/s}}{\abs{k}}\right)\left[\jap{k,(Rn)^{1/s}}^s - \jap{(Rn)^{1/s}}^s\right]\right] \frac{\jap{k,(Rn)^{1/s}}^{\sigma}}{\jap{(Rn)^{1/s}}^{\sigma}} \nonumber \\ & \quad\quad \times \exp \left[\frac{\abs{kt}}{(Rn)^{1/s}}\left(\lambda\left(\frac{(Rn)^{1/s}}{\abs{k}}\right)\jap{(Rn)^{1/s}}^s + \sigma \log\jap{(Rn)^{1/s}} \right) \right] \nonumber \\ 
& \leq \exp\left[\lambda\left(t\right)\jap{k,kt}^s \right]
\jap{k,kt}^{\sigma}\jap{(Rn)^{1/s}}^{\sigma\left(\frac{(n+1)^{1/s}}{n^{1/s}} -1\right)} \nonumber
 \\ & \quad\quad \times \exp\left[\lambda\left(\frac{(Rn)^{1/s}}{\abs{k}}\right)\left(\frac{(n+1)^{1/s}}{(n)^{1/s}} -  1\right)\jap{(Rn)^{1/s}}^s\right] \nonumber \\ 
& \lesssim_{R,\lambda_0,s} A_k(t,kt). \label{ineq:emuNkn}
\end{align} 
A similar result holds for $n=0$ when $\abs{kt}^s \leq R$ using that $\lambda(t)$ is non-increasing, 
\begin{align*} 
N_{k,0} e^{\mu_0\abs{kt}}  & = \exp\left[\lambda\left(\frac{R^{1/s}}{\abs{k}}\right)\jap{k}^s\right] \jap{k}^{\sigma} 
\exp \left[\frac{\abs{kt}}{R^{1/s}}\left(\lambda\left(\frac{R^{1/s}}{\abs{k}}\right)\jap{R^{1/s}}^s + \sigma \log\jap{R^{1/s}} \right) \right] \\
& \lesssim_R A_k(t,kt).  
\end{align*} 
Hence it follows from \eqref{ineq:PhiBound} that 
\begin{align}
\norm{\Phi_k^n(t)}_{L_t^2(\Real)} \lesssim \frac{1}{\kappa}\norm{AF_k(t) \mathbf{1}_{Rn \leq \abs{kt}^s \leq R(n+1)}}_{L_t^2(\Real)}. \label{ineq:LinearBound}
\end{align}

\noindent
{\em Step 5. Summation of the different frequency blocks.} 
Now we want to estimate $A\phi_k$ by summing over $n$ which will
require some almost orthogonality.  This is possible as each
$\phi_k^n$ is exponentially localized near $\abs{kt}^s \approx Rn$.
Computing, noting that by construction $\phi_k^n$ is only supported on
$\abs{kt}^s \geq Rn$,
\begin{align} 
\norm{A\phi_k}_{L_t^2}^2 & = \int_0^{+\infty} \abs{\sum_{n=0}^{+\infty} A_k(t,kt)\phi^n_k(t) }^2 \, dt \nonumber \\ 
  & \lesssim \int_0^{\infty} \abs{A_k(t,kt)\phi^0_k(t)}^2 \, dt +
  \int_0^{\infty} 
  \abs{\sum_{n=1}^{\infty} A_k(t,kt)e^{-\mu_n\abs{kt}} N_{k,n}^{-1} \mathbf{1}_{\abs{kt}^s \geq Rn} \Phi^n_k(t) }^2 \, dt \nonumber \\ 
  & =  \int_0^{\infty} \abs{A_k(t,kt)\phi^0_k(t)}^2 \, dt \nonumber \\
  & \quad + \int_0^{\infty} \sum_{n,n^\prime \geq 1} 
A_k(t,kt)^2 e^{-\mu_n\abs{kt}} N_{k,n}^{-1}e^{-\mu_{n^\prime}\abs{kt}}
N_{k,n^\prime}^{-1} \mathbf{1}_{\abs{kt}^s \geq Rn} 
\mathbf{1}_{\abs{kt}^s \geq Rn^\prime} 
\Phi^n_k(t) \Phi^{n^\prime}_k(t) dt.  \label{ineq:AphiLinearbd}
\end{align}
First we will approach the infinite sum, which is the more challenging
term. For this we use Schur's test. Indeed, if we denote
the interaction kernel
\begin{align*} 
K_{n,n^\prime}(t,k) := A_k(t,kt)^2 e^{-\mu_n\abs{kt}}
N_{k,n}^{-1}e^{-\mu_{n^\prime}\abs{kt}} 
N_{k,n^\prime}^{-1} \mathbf{1}_{\abs{kt}^s \geq Rn} \mathbf{1}_{\abs{kt}^s \geq Rn^\prime},
\end{align*}
then Schur's test (or Cauchy-Schwarz three times) implies 
\begin{align} 
\int_0^{+\infty} \sum_{n,n^\prime \geq 1}K_{n,n^\prime}(t,k)\Phi^{n^\prime}_k(t)\Phi^{n}_k(t) \, dt 
 & \nonumber \\ & \hspace{-5.5cm} \leq  \left(\sup_{t \in
     [0,\infty)}\sup_{n\geq 1} 
   \sum_{n^\prime=1}^{\infty} K_{n,n^\prime}(t,k)\right)^{1/2} 
 \left(\sup_{t \in [0,\infty)}\sup_{n^\prime \geq 1} \sum_{n=1}^{\infty}
   K_{n,n^\prime}(t,k)\right)^{1/2} \sum_{n=1}^{\infty}
 \norm{\Phi_k^n(t)}^2_{L^2_t(\R)}. 
 \label{ineq:AO}
\end{align} 
It remains to see that the row and column sums of the interaction
kernel are uniformly bounded in time.  Since the kernel is symmetric
in $n$ and $n^\prime$ it suffices to consider only one of the sums.
The computations above to deduce \eqref{ineq:emuNkn} can be adapted to show at least that
\begin{align*}
A_k(t,kt)e^{-\mu_{n^\prime}\abs{kt}} N_{k,n^\prime}^{-1}\mathbf{1}_{\abs{kt}^s \geq Rn^\prime} \lesssim_R 1, 
\end{align*}
and hence, since $\lambda(t)$ is decreasing,
\begin{align*} 
  \sum_{n=1}^{+\infty} K_{n,n^\prime}(t,k) 
  & \lesssim_R \sum_{n=1}^{+\infty} A_k(t,kt) e^{-\mu_n\abs{kt}} N_{k,n}^{-1}\mathbf{1}_{\abs{kt}^s \geq Rn}  \\
& = \sum_{n=1}^{+\infty}  \mathbf{1}_{\abs{kt}^s \geq Rn} e^{\lambda(t)\jap{k,kt}^s - \lambda\left(\frac{(Rn)^{1/s}}{\abs{k}}\right)\left(\jap{k,(Rn)^{1/s}}^s - \jap{(Rn)^{1/s}}^s\right)  }\frac{\jap{k,kt}^\sigma \jap{(Rn)^{1/s}}^\sigma}{\jap{k,(Rn)^{1/s}}^\sigma} e^{-\mu_n\abs{kt}} \\
& \lesssim \sum_{n = 1}^{+\infty} \mathbf{1}_{\abs{kt}^s \geq Rn} e^{\lambda(t)\jap{k,kt}^s - \lambda\left(t\right)\left(\jap{k,(Rn)^{1/s}}^s - \jap{(Rn)^{1/s}}^s\right)  }\frac{\jap{k,kt}^\sigma \jap{(Rn)^{1/s}}^\sigma}{\jap{k,(Rn)^{1/s}}^\sigma} \\
& \quad\quad \times e^{-\frac{\abs{kt}}{(Rn)^{1/s}}\left[ \lambda\left(t\right)\jap{(Rn)^{1/s}}^s + \sigma \log\jap{(Rn)^{1/s}}\right]}. 
\end{align*} 
Since $e \leq (Rn)^{1/s} \leq \abs{kt}$, we have 
\begin{align*} 
-\jap{k,(Rn)^{1/s}}^s + \jap{(Rn)^{1/s}}^s \leq \jap{kt}^s - \jap{k,kt}^{s}, 
\end{align*}  
since the LHS is increasing as a function of $n$ and therefore,
\begin{align*} 
\sum_{n=1}^{+\infty} K_{n,n^\prime}(t,k) & 
\lesssim \sum_{n = 1}^{+\infty} \mathbf{1}_{\abs{kt}^s \geq Rn} 
e^{\lambda(t)\jap{kt}^s\left[ 1- \frac{\abs{kt}\jap{(Rn)^{1/s}}^s}{ (Rn)^{1/s} \jap{kt}^s}\right]}
e^{\sigma\log\left[\frac{\jap{k,kt}\jap{(Rn)^{1/s}}}{\jap{k,(Rn)^{1/s}}}\right] - \sigma\frac{\abs{kt}}{(Rn)^{1/s}}\log \jap{(Rn)^{1/s}}}.
\end{align*}  
Since $\jap{x} \jap{k,x}^{-1}$ is increasing in $x$ for $\abs{k} \geq
1$ and $x \ge 0$, we have
\begin{equation*} 
\sum_{n=1}^{+\infty} K_{n,n^\prime}(t,k)  \lesssim \sum_{n = 1}^{+\infty}
\mathbf{1}_{\abs{kt}^s \geq Rn} e^{\lambda(t)\jap{kt}^s\left[ 1-
    \frac{\abs{kt}\jap{(Rn)^{1/s}}^s}{ (Rn)^{1/s} \jap{kt}^s}\right]}
e^{\sigma\log\jap{kt}\left( 1 - \frac{\abs{kt} \log
      \jap{(Rn)^{1/s}}}{(Rn)^{1/s} \log \jap{kt}} \right)}.
\end{equation*}
Finally using that $x / \log \jap{x}$ is increasing for $x \ge e$ we
get 
\begin{equation*}
  \sum_{n=1}^{+\infty} K_{n,n^\prime}(t,k)  
  \lesssim \sum_{n = 1}^{+\infty} \mathbf{1}_{\abs{kt}^s \geq Rn} 
  e^{\lambda(t)\jap{kt}^s\left[ 1- \frac{\abs{kt}\jap{(Rn)^{1/s}}^s}{ (Rn)^{1/s} \jap{kt}^s}\right]}.  
\end{equation*} 
Using that $\jap{x^{1/s}}^s \geq x$, the sum can be bounded by 
\begin{align*} 
\sum_{n=1}^\infty K_{n,n^\prime}(t,k) & \lesssim_{R} e^{\lambda(t)\jap{kt}^s} \int_R^{\abs{kt}^s}e^{-\lambda(t)\abs{kt}\frac{\jap{x^{1/s}}^s}{x^{1/s}}} dx \\
 & \lesssim e^{\lambda(t)\jap{kt}^s} \int_R^{\abs{kt}^s}e^{-\lambda(t)\abs{kt}x^{1-1/s}} dx \\
& \lesssim e^{\lambda(t)\jap{kt}^s} \int_{\abs{kt}^{s-1}}^{R^{1-1/s}} e^{-\lambda(t)\abs{kt}\tau} \tau^{\frac{1}{s-1}} \, d\tau \lesssim 1. 
\end{align*} 
This shows that the row sums of $K_{n,n^\prime}(t,k)$ are uniformly
bounded; by symmetry the column sums are also bounded and by
\eqref{ineq:AO}, this completes the treatment of the summation in
\eqref{ineq:AphiLinearbd}.

Now we turn our attention to the $n = 0$ term in
\eqref{ineq:AphiLinearbd}. By \eqref{ineq:sorttimeLinear},
\begin{align} \nonumber
\int_0^{+\infty} \abs{A_k(t,kt)\phi_k^0(t)}^2 \, dt & =
\int_0^{R^{1/s}} \abs{A_k(t,kt)\phi_k^0(t)}^2 \, dt  + \int_{R^{1/s}}
^{+\infty} \abs{A_k(t,kt)\phi_k^0(t)}^2 \, dt \\ & \lesssim_R
\norm{AF_k}^2_{L^2_t(\R)} 
+ \int_{R^{1/s}}^{+\infty} \abs{A_k(t,kt)\phi_k^0(t)}^2 \, dt. \label{ineq:phik0}
\end{align} 
However, for $\abs{t} \geq R^{1/s}$, we have by \eqref{ineq:JapTri}, that $\lambda(t)$ is non-increasing and \eqref{ineq:IncExp},
\begin{align*} 
A_k(t,kt) \leq e^{\lambda(t)\jap{k}^s + \lambda(t)\jap{kt}^s}\jap{k}^\sigma \jap{kt}^\sigma \lesssim_R  e^{\lambda\left(\frac{R^{1/s}}{\abs{k}}\right)\jap{k}^s}\jap{k}^\sigma e^{\mu_0\abs{kt}} \lesssim N_{k,0} e^{\mu_0\abs{kt}},  
\end{align*} 
which implies with \eqref{ineq:phik0} and \eqref{ineq:LinearBound} that 
\begin{align} 
 \int_0^{\infty} \abs{A_k(t,kt)\phi_k^0(t)}^2 dt \lesssim \left(1 + \frac{1}{\kappa^2}\right)\norm{AF_k}_{L^2_t(\R)}.\label{ineq:phi0kfinal}
\end{align} 

Combining \eqref{ineq:AphiLinearbd}, \eqref{ineq:phi0kfinal},
\eqref{ineq:AO} with \eqref{ineq:LinearBound} we have
\begin{align*} 
\norm{A\phi_k}_{L^2_t(I)}^2 & \lesssim_R \left(1 + \frac{1}{\kappa^2}\right)\norm{AF_k}^2_{L^2_t(\R)} +  \sum_{n=1}^\infty \norm{\Phi^n_k(t)}^2_{L^2_t(I)} \\  & \lesssim \left(1 + \frac{1}{\kappa^2}\right) \sum_{n=0}^\infty \norm{AF_k(t) \mathbf{1}_{Rn \leq \abs{kt}^s \leq R(n+1)}}^2_{L^2_t(\R)}, 
\end{align*}
which completes the proof of the lemma. 

\subsection{Rigorous justification of the a priori estimate}
\label{sec:rigor-just-priori}

The reader may have noticed that in the previous subsection it seems
that we only used the bound from below $|1-\mathcal L(k,\xi)|\ge
\kappa$ with $\xi = \mu_n + i \omega/|k|$, i.e. in the strip $\Re e\,
\xi \in (0,\bar \lambda/2)$. The subtlety is that the Fourier-Laplace
transform of $\Phi_k^n (t)$ is only granted to exist when some $L^2$
integrability as $t \to \infty$ is known. 

To be more more specific, consider \eqref{def:PhikVolt}.  
From the Gr\"onwall bound \eqref{ineq:sorttimeLinear} established in step 1, it is clear that the
Fourier-Laplace transform would exist if one chooses $\mu_n< -2C$, however it is not clear that 
we can perform the computation as $\mu_n$ approaches the imaginary axis. 
In order to avoid a circular argument --establishing some time decay
by assuming the existence of a Fourier-Laplace transform which already
requires some time decay--, we can appeal to several arguments:
\begin{enumerate}
\item We can use as a black box the \emph{Paley-Wiener theory}
  (see~\cite[Chap.~18]{Paley-Wiener}
  or~\cite[Chap.~2]{book-volterra}): for every $f \in L^1_{loc}(\R_+)$
  there exists a unique solution $u \in L^1_{loc}(\R_+)$ to the
  integral equation $u = f + k * u$ with $k \in L^1_t$ and with $k \le
  e^{Ct}$ for some constant $C>0$ (where the convolution over $t \in
  \R_+$ is defined as before by extending functions to zero on $\R_-$),
  given by $u = f - f* r$ where $r \in L^1_{loc}(\R_+)$ is the
  so-called resolvent kernel of $k$. The latter is the unique solution
  to $r = k + r *k$ and the key result of the theory is that $r \in
  L^1(\R_+)$ iff the Fourier-Laplace transform of $k$ satisfies
  $\mathcal L[k](\xi) \not =1$ for any $\Re e \, \xi \le 0$. As soon
  as $r, f \in L^1(\R_+)$ we have $u \in L^1(\R_+)$. Then step 3 of Lemma \ref{lem:LinearL2Damping} can be justified by applying this theory to
  $u(t) := \Phi^n_k(t)$, $f(t) :=R^n_k(t)$ and $k(t) :=
  e^{\mu_n\abs{k}t} K^0(t,k)$. 

\item A second method is to use an approximation argument, in the
  spirit of energy methods in PDEs, which was discussed in
  \cite[Section 3]{villani2010}. Define $\Phi^n_{k,\delta}(t) :=
  \Phi^n_k(t) e^{-\delta t^2/2}$.
Now we have the existence of Fourier-Laplace transform of
$\Phi^n_k(t)$ for any $\mu_n$ thanks to the Gaussian decay in time and \eqref{ineq:sorttimeLinear}. 
Using the Fourier-Laplace transform on \eqref{def:PhikVolt} and \textbf{(L)}, we may deduce that for $\mu_n < -2C$ (where $2C$ comes from the estimate \eqref{ineq:sorttimeLinear}) that the following formula holds:
\[
\hat \Phi^n_{k,\delta}(\omega) = \left( \frac{\hat R^n_k(\omega)}{1 - \mathcal
    L\left(k,\mu_n + i \frac{\omega}{|k|} \right)} \right) \ast \gamma_\delta,
  \quad \gamma_\delta(\omega) := \frac{e^{-\frac{|\omega|^2}{2
        \delta}}}{\sqrt{2 \pi \delta}}.
\]
Since this is an analytic function in $\mu_n$ and $\omega$ as long as we do not approach a singularity, by analytic continuation we may deduce that this formula remains valid for all $\mu_n < \bar\lambda$ by \textbf{(L)}.  
Therefore by Plancherel's theorem
\[
\norm{\Phi_{k,\delta}^n(t)}_{L_t^2(\R)} \lesssim
\frac{1}{\kappa}\norm{R_k^n(t)}_{L^2_t(\R)}
\norm{\gamma_\delta}_{L^1(\R)} \lesssim \frac{1}{\kappa}\norm{R_k^n(t)}_{L^2_t(\R)} 
\]
which is an estimate independent of $\delta>0$. We then let $\delta$
go to zero and deduce by Fatou's lemma the desired bound \eqref{ineq:PhiBound}:
\[
\norm{\Phi_{k}^n(t)}_{L_t^2(\R)} \lesssim \frac{1}{\kappa}\norm{R_k^n(t)}_{L^2_t(\R)}
\]
(which also justifies the existence of the Fourier-Laplace
transform). 
\end{enumerate}
\begin{remark} To finish, let us mention that the present discussion is
  related to the Gerhart-Herbst-Pr\"uss-Greiner theorem
  \cite{MR0461206,MR715559,MR743749,MR839450} (see also
  \cite{MR1721989}) for semigroups in Hilbert spaces. The latter
  asserts that the semigroup decay is given by the spectral bound,
  under a sole pointwise control on the resolvent.  While the
  constants seem to be non-constructive in the first versions of this
  theorem, Engel and Nagel gave a comprehensive and elementary proof
  with constructive constant in \cite[Theorem 1.10; chapter
  V]{MR1721989}. Let us also mention on the same subject subsequent
  more recent works like \cite{Helffer-Sjoestrand}. The main idea in
  the proof of \cite[Theorem 1.10, chapter V]{MR1721989}, which is
  also used in \cite{Helffer-Sjoestrand}, is to use a Plancherel
  identity on the resolvent in Hilbert spaces in order to obtain
  explicit rates of decay on the semigroup in terms of bounds on the
  resolvent. However in the Volterra integral equation we study here,
  there is no semigroup structure on the unknown $\phi_k(t)$, and we
  cannot appeal directly to these classical results. 
\end{remark}

\subsection{The Penrose criterion}
\label{sec:penrose-criterion}

A generalised form of the Penrose criterion \cite{Penrose} was given in
\cite{MouhotVillani11} as follows:
\begin{equation}
  \label{eq:penrose}
  {\bf (P)} \qquad \qquad \forall \, k \in \Z^d_* \ \mbox{ and } \ r
  \in \R
  \ \mbox{ s.t. } \ (f^0_k)'(r) =0,
  \quad \quad \widehat W(k) \left( \mbox{p.v.} \int_\R
    \frac{(f^0_k)'(r)}{r-w} \, dr\right) <1,
\end{equation} 
where $f_k^0$ denotes the marginals of the background $f^0$ along arbitrary wave vector
$k \in \Z^d_*$:
\[
f^0_k(r) := \int_{kr/|k|+k^\bot} f^0(w) \, dw, \quad r \in \R.
\]
The proof that condition {\bf (P)} implies the condition
{\bf (L)} was not quite complete in \cite{MouhotVillani11} as it was proved
only that {\bf (P)} implies the lower bound $|1-\mathcal L(k,\xi)|\ge
\kappa$ in a strip and not in a half-plane. The complete proof due to
Penrose relies on the argument principle. The starting point is to
observe that $\mathcal L(k,\xi) = \int_0 ^{+\infty} e^{\bar \xi|k|t}
K^0(t,k) \, dt$ with $\xi = \lambda + i \zeta$ and $K^0(t,k) :=
-\hat f^0\left(kt\right) \widehat{W}(k)\abs{k}^2t$ is well-defined
for $\lambda < \bar \lambda$ by the analyticity of $f^0$ and is small for large $\zeta$ by
integration by parts.
 We
therefore restrict ourselves to a compact interval $|\zeta|\le C$
in $\zeta$, and we compute by the Plemelj formula (see \cite{MouhotVillani11} for more details) that 
\begin{align} 
\mathcal L(k,i\zeta) = \hat W(k) \left[ \left( \mbox{p.v.} \int_\R
    \frac{(f^0_k)'(r)}{r-\zeta} \, dr \right) - i \pi (f^0_k)'(\zeta) \label{eq:LtoPlemelj}
\right].
\end{align}
Therefore, the condition {\bf (P)} implies that $|1-\mathcal
L(k,\xi)| \ge 2\kappa$ for some $\kappa>0$ at $\lambda=0$ and for
$|\zeta|\le C$. Combined with smallness for large $\zeta$ and continuity,
we deduce the lower bound $|1-\mathcal L(k,\xi)| \ge \kappa$ in a
strip $\Re e \, \xi \in [0,\lambda']$ for some $\lambda'>0$. Since the
function $\xi \mapsto \mathcal L(k,\xi)$ is holomorphic on $\Re e \,
\xi < \bar \lambda$ and the value $1$ is
not taken on $i\R$, by the argument principle, the value $1$ can only be taken on $\Re e \, \xi <0$ if
$\Xi : \zeta \mapsto \mathcal L(k,i\zeta)$ has a positive winding number
around this value. However, this would imply that the curve $\Xi$ crosses
the real axis above $1$, which is prohibited again by \eqref{eq:LtoPlemelj} and 
{\bf (P)}, which concludes the proof. 

\section{Energy estimates on the $(\rho,f)$ system} \label{sec:Energy}
In this section we perform the necessary energy estimates to prove Proposition \ref{lem:Boot}, that is, we deduce the multi-tier controls stated in \eqref{ctrl:BootRes} from \eqref{ctrl:Boot} for suitable $K_i$ and sufficiently small $\epsilon$.

\subsection{$L^2_t(I)$ estimates on the density, \eqref{ctrl:Mid}} \label{sec:L2I}
The most fundamental estimate we need to make is the $L_t^2(I)$ control \eqref{ctrl:Mid}, which requires the linear damping Lemma \ref{lem:LinearL2Damping} and the crucial plasma echo analysis carried out in \S\ref{sec:Proofa} (which we apply as a black box in this section).
The controls \eqref{ctrl:LowCommLoc} and \eqref{ctrl:HiLocalized} were chosen specifically for this. 

To highlight its primary importance, we state the inequality as a separate proposition. 
\begin{proposition}[Nonlinear control of $\rho$] \label{prop:nonlinearmoment}
For suitable $K_i$ and $\epsilon_0$ sufficiently small, the estimate \eqref{ctrl:Mid} holds under the bootstrap hypotheses \eqref{ctrl:Boot}.
\end{proposition} 
\begin{proof} 
Proposition \ref{prop:nonlinearmoment} requires two main controls, which follow from \eqref{ctrl:LowCommLoc} and \eqref{ctrl:HiLocalized} respectively.   
\begin{itemize} 
\item[(a)] Define the time-response kernel $\bar{K}_{k,\ell}(t,\tau)$ for some $c = c(s) \in (0,1)$ (determined by the proof): 
\begin{align} 
\bar{K}_{k,\ell}(t,\tau) & = \frac{1}{\abs{\ell}^{\gamma}}e^{(\lambda(t) - \lambda(\tau))\jap{k,kt}^s}e^{c\lambda(\tau)\jap{k-\ell,kt-\ell\tau}^s}\abs{k(t-\tau)\hat{f}_{k-\ell}(\tau,kt-\ell\tau)} \mathbf{1}_{\ell \neq 0}. \label{def:barK}
\end{align}  
Proposition \ref{prop:nonlinearmoment} will depend on the estimate
\begin{align} 
\left(\sup_{t \geq 0} \sup_{k \in \Integer_\ast^d}\int_0^t\sum_{\ell \in \Integer_\ast^d} \bar{K}_{k,\ell}(t,\tau) d\tau\right)
 \left(\sup_{\tau \geq 0} \sup_{\ell\in \Integer_\ast^d} \sum_{k \in \Integer_\ast^d} \int_{\tau}^\infty\bar{K}_{k,\ell}(t,\tau) dt\right)
 \lesssim K_2\epsilon^2. \label{ineq:momentdef}
\end{align}

\item[(b)] Proposition \ref{prop:nonlinearmoment} will also depend on the estimate
\begin{align} 
\sup_{\tau \geq 0}e^{(c-1)\alpha_0\jap{\tau}^s} \sum_{k \in \Integer^d} \sup_{\omega \in \Integer_\ast^d} \sup_{x\in \Real^d}  \int_{-\infty}^\infty \abs{(A \widehat{\grad f})_{k}\left(\tau,\frac{\omega}{\abs{\omega}}\zeta- x\right)}^2 d\zeta \lesssim K_1\epsilon^2. \label{ineq:controlLines}
\end{align}
\end{itemize}
The condition \eqref{ineq:momentdef} controls \emph{reaction:} the interaction of the density with the lower frequencies of $f$; condition \eqref{ineq:controlLines} controls \emph{transport:} the interaction with higher frequencies of $f$. 
The latter, condition \eqref{ineq:controlLines}, follows from \eqref{ctrl:HiLocalized}: by Lemma \ref{lem:SobTrace} followed by \eqref{ineq:OuterActrl},
\begin{align*} 
\sum_{k \in \Integer^d} \sup_{\omega \in \Integer_\ast^d} \sup_{x \in \Real^d}  \int_{-\infty}^\infty \abs{(A \widehat{\grad f})_{k}\left(\tau,\frac{\omega}{\abs{\omega}}\zeta-x\right)}^2 d\zeta & \lesssim_M \norm{A \widehat{\grad f}(\tau)}^2_{L_k^2 H^M_\eta} \lesssim K_1 \epsilon^2 \jap{\tau}^7, 
\end{align*} 
from which \eqref{ineq:controlLines} follows by \eqref{ineq:SobExp} and $c < 1$. 
Since condition \eqref{ineq:momentdef} is much harder to verify and contains the physical mechanism of the plasma echoes, we prove Proposition \ref{prop:nonlinearmoment} assuming \eqref{ineq:momentdef}. In \S\ref{sec:Proofa} below, we prove that \eqref{ineq:momentdef} follows from \eqref{ctrl:LowCommLoc}. 

Expanding the integral equation \eqref{eq:rhointegral} 
using the paraproduct decomposition: 
\begin{align} 
\hat{\rho}_k(t) & = \hat h_{in}(k,kt) - \int_0^t \hat{\rho}_k(\tau)\abs{k}^2 \widehat{W}(k)(t-\tau)\hat{f}^0(k(t-\tau))  d\tau  \nonumber \\ 
& \quad - \int_0^t\sum_{\ell \in \Integer_\ast^d}\sum_{N \geq 8} \hat{\rho}_\ell(\tau)_{<N/8} \widehat W(\ell) \ell \cdot k(t-\tau) \hat f_{k-\ell}(\tau,kt-\ell\tau)_{N} d\tau \nonumber
\\ 
& \quad - \int_0^t\sum_{\ell \in \Integer_\ast^d} \sum_{N \geq 8} \hat{\rho}_\ell(\tau)_{N} \widehat W(\ell) \ell \cdot k(t-\tau) \hat f_{k-\ell}(\tau,kt-\ell\tau)_{<N/8} d\tau \nonumber \\ 
& \quad - \int_0^t\sum_{\ell \in \Integer_\ast^d}\sum_{N \in \mathbb{D}} \sum_{N/8 \leq N^\prime \leq 8N} \hat{\rho}_\ell(\tau)_{N^\prime} \widehat W(\ell) \ell \cdot k(t-\tau) \hat f_{k-\ell}(\tau,kt-\ell\tau)_{N} d\tau, \nonumber \\ 
& = \hat h_{in}(k,kt) -  \int_0^t \hat{\rho}_k(\tau)\abs{k}^2 \widehat{W}(k)(t-\tau)\hat{f}^0(k(t-\tau))  d\tau  - T_k(t) - R_k(t) - \mathcal{R}_k(t). \label{def:pararho}
\end{align} 
Recall our convention that the Littlewood-Paley projection of $\hat\rho_\ell(\tau)$ treats $\ell\tau$ in place of the $v$ frequency.
We begin by applying Lemma \ref{lem:LinearL2Damping} to \eqref{def:pararho}, which implies for each $k \in \Integer_\ast^d$, 
\begin{align} 
\norm{A\hat{\rho}_k}^2_{L_t^2(I)} & \lesssim C^2_{LD}\norm{A_k(\cdot,k\cdot)\hat h_{in}(k,k\cdot)}^2_{L_t^2(I)} + C^2_{LD}\norm{AT_k}^2_{L_t^2(I)} \nonumber \\ & \quad + C^2_{LD}\norm{AR_k}^2_{L_t^2(I)} + C^2_{LD}\norm{A\mathcal{R}_k}^2_{L_t^2(I)}. \label{ineq:MomentInPf}
\end{align}
First, Lemma \ref{lem:SobTrace} and a version of the argument applied in \eqref{ineq:OuterMctrl} (using also that $\lambda(t)$ is decreasing \eqref{def:lambda}) imply 
\begin{align} 
\sum_{k \in \Integer_\ast^d} \norm{A_k(\cdot,k\cdot)\hat{h}_{in}(k,k\cdot)}^2_{L_t^2(I)} &  = \sum_{k \in \Integer_\ast^d} \int_0^{T^\star} \abs{A_k(t,kt) \hat{h}_{in}(k,kt)}^2 dt \nonumber \\ 
& \lesssim \sum_{k \in \Integer_\ast^d}\norm{A_k(0,\cdot)\hat{h}_{in}(k,\cdot)}_{H^M(\Real^d_\eta)}^2 \nonumber \\ 
& \lesssim \epsilon^2. \label{ineq:rhoL2ID}
\end{align} 

Now we turn to the nonlinear contributions in \eqref{def:pararho}.  

\subsubsection{Reaction}
Our goal is to show
\begin{align} 
\norm{AR}_{L_k^2L^2_t(I)}^2 \lesssim  K_2\epsilon^2\norm{A\hat{\rho}}^2_{L_k^2L^2_t(I)}, \label{ineq:rhoReacConc}
\end{align} 
since for $\epsilon$ chosen sufficiently small, this contribution can then be absorbed on the LHS of \eqref{ineq:MomentInPf}. 

First, by applying \eqref{ineq:Wbd}, 
\begin{align*} 
\norm{AR}_{L_k^2L^2_t(I)}^2  \lesssim  \sum_{k \in \Integer_\ast^d}\int_0^{T^\star} \left[A_k(t,kt)\int_0^t\sum_{\ell \in \Integer_\ast^d} \sum_{N \geq 8}\abs{\hat{f}_{k-\ell}(\tau,kt-\ell\tau)_{<N/8} \frac{\abs{k(t-\tau)}}{\abs{\ell}^{\gamma}} \hat{\rho}_\ell(\tau)_N} d\tau \right]^{2} dt.
\end{align*} 
By definition, the Littlewood-Paley projections imply the frequency localizations (as in \eqref{ineq:FreqLocExample}): 
\begin{subequations} \label{ineq:rhoReactionFreqLoc}
\begin{align} 
\frac{N}{2} \leq \abs{\ell} + \abs{\ell \tau} & \leq \frac{3N}{2}, \\ 
\abs{k-\ell} + \abs{kt-\ell\tau} & \leq \frac{3N}{32}, \\ 
\frac{13}{16} \leq \frac{\abs{k,kt}}{\abs{\ell,\tau \ell}} & \leq \frac{19}{16}. 
\end{align} 
\end{subequations}
From \eqref{ineq:rhoReactionFreqLoc}, on the support of the integrand, \eqref{lem:scon2} implies that  for some $c = c(s) \in (0,1)$:  
\begin{align} 
A_k(t,kt) & = e^{(\lambda(t) - \lambda(\tau))\jap{k,kt}^s}A_k(\tau,kt) \lesssim e^{(\lambda(t) - \lambda(\tau))\jap{k,kt}^s}e^{c\lambda(\tau)\jap{k-\ell,kt-\ell\tau}^s}A_\ell(\tau,\ell\tau). \label{ineq:ReacRhoA}
\end{align} 
Therefore, by definition of $\bar{K}$ we have (dropping the Littlewood-Paley projection on $f$), 
\begin{align*} 
\norm{AR}_{L_k^2L^2_t(I)}^2 & \lesssim \sum_{k\in \Integer_\ast^d}\int_0^{T^\star} \left[\int_0^t\sum_{\ell \in \Integer_\ast^d} \bar{K}_{k,\ell}(t,\tau) A_\ell(\tau,\ell\tau)\sum_{N \geq 8} \abs{\hat{\rho}_\ell(\tau)_N} d\tau \right]^{2} dt.
\end{align*} 
Since the Littlewood-Paley projections define a partition of unity,
\begin{align*} 
\norm{AR}_{L_k^2L^2_t(I)}^2 & \lesssim \sum_{k\in \Integer_\ast^d}\int_0^{T^\star} \left[\int_0^t\sum_{\ell \in \Integer_\ast^d} \bar{K}_{k,\ell}(t,\tau) A_\ell(\tau,\ell\tau)\abs{\hat{\rho}_\ell(\tau)} d\tau \right]^{2} dt.
\end{align*} 
From here we may proceed analogous to \S7 in \cite{MouhotVillani11}, for which we apply Schur's test in $L_k^2L_t^2$. Indeed, 
\begin{align*} 
\norm{AR}_{L_k^2L^2_t(I)}^2& \lesssim \sum_{k\in\Integer_\ast^d}\int_0^{T^\star}\left(\int_0^t\sum_{\ell \in \Integer_\ast^d} \bar{K}_{k,\ell}(t,\tau) d\tau\right)\left(\int_0^t \sum_{\ell \in \Integer_\ast^d} \bar{K}_{k,\ell}(t,\tau) \abs{A_\ell(\tau,\ell\tau)\hat{\rho}_\ell(\tau)}^2 d\tau\right) dt \\
& \leq \left(\sup_{t \geq 0} \sup_{k\in \Integer_\ast^d}\int_0^t\sum_{\ell \in \Integer_\ast^d} \bar{K}_{k,\ell}(t,\tau) d\tau\right) \sum_{k\in \Integer_\ast^d}\int_0^{T^\star}\left(\int_0^t \sum_{\ell \in \Integer_\ast^d} \bar{K}_{k,\ell}(t,\tau) \abs{A_\ell(\tau,\ell\tau)\hat{\rho}_\ell(\tau)}^2 d\tau\right) dt.
\end{align*} 
By Fubini's theorem,  
\begin{align*} 
\norm{AR}_{L_k^2L^2_t(I)}^2 & \lesssim  \left(\sup_{t \geq 0} \sup_{k\in \Integer_\ast^d}\int_0^t\sum_{\ell \in \Integer_\ast^d} \bar{K}_{k,\ell}(t,\tau) d\tau\right) \sum_{\ell \in \Integer_\ast^d} \int_0^{T^\star} \left( \int_{\tau}^{T^\star} \sum_{k\in \Integer_\ast^d} \bar{K}_{k,\ell}(t,\tau) dt\right) \abs{A_\ell(\tau,\ell\tau)\hat{\rho}_\ell(\tau)}^2  d\tau \\ 
& \lesssim \left(\sup_{t \geq 0} \sup_{k\in \Integer_\ast^d}\int_0^t\sum_{\ell \in \Integer_\ast^d} \bar{K}_{k,\ell}(t,\tau) d\tau\right) \left(\sup_{\tau \geq 0} \sup_{\ell\in \Integer_\ast^d} \sum_{k\in\Integer_\ast^d} \int_{\tau}^{T^\star}\bar{K}_{k,\ell}(t,\tau) dt\right) \norm{A\hat{\rho}}^2_{L^2_k L_t^2(I)}. 
\end{align*} 
Hence, condition \eqref{ineq:momentdef} implies \eqref{ineq:rhoReacConc}.

\subsubsection{Transport} 
As above, first apply \eqref{ineq:Wbd}, 
\begin{align*} 
\norm{AT}_{L_k^2L_t^2(I)}^2  \lesssim \sum_{k \in \Integer_\ast^d}\int_0^{T^\star} \left[A_k(t,kt)\int_0^t\sum_{\ell \in \Integer_\ast^d} \sum_{N \geq 8} \abs{\hat{f}_{k-\ell}(\tau,kt-\ell\tau)_{N} \frac{\abs{k(t-\tau)}}{\abs{\ell}^{\gamma}} \hat \rho_\ell(\tau)_{<N/8}} d\tau \right]^{2} dt.
\end{align*} 
By the Littlewood-Paley projections, on the support of the integrand there holds, 
\begin{subequations} \label{ineq:transptRhoLoc}
\begin{align} 
\frac{N}{2} \leq \abs{k-\ell} + \abs{kt - \ell\tau} & \leq \frac{3N}{2} \\ 
\abs{\ell} + \abs{\ell\tau} & \leq \frac{3N}{32} \\
\frac{13}{16} \leq \frac{\abs{k,kt}}{\abs{k-\ell,kt-\tau \ell}} & \leq \frac{19}{16}. 
\end{align} 
\end{subequations}
By \eqref{ineq:transptRhoLoc}, on the support of the integrand, \eqref{lem:scon2} implies that for some $c = c(s) \in (0,1)$:  
\begin{align*} 
A_k(t,kt) & = e^{(\lambda(t) - \lambda(\tau))\jap{k,kt}^s}A_k(\tau,kt) \lesssim e^{(\lambda(t) - \lambda(\tau))\jap{k,kt}^s}e^{c\lambda(\tau)\jap{\ell,\ell\tau}^s}A_{k-\ell}(\tau,kt-\ell\tau). 
\end{align*} 
Using that 
\begin{align} 
\abs{k(t-\tau)} \leq \abs{kt-\ell\tau} + \tau\abs{k-\ell} \leq \jap{\tau}\abs{k-\ell,kt-\ell\tau}, \label{ineq:Echot}
\end{align}
 we have (ignoring the Littlewood-Paley projection on $\rho$ and the $\abs{\ell}^{-\gamma}$ which are not helpful), 
\begin{align*} 
 \norm{AT}_{L_k^2L_t^2(I)} & \lesssim \sum_{k \in \Integer_\ast^d} \int_0^{T^\star}\left[\sum_{\ell\in \Integer_\ast^d} \int_0^t \sum_{N \geq 8} \abs{(A \widehat{\grad f})_{k-\ell}(\tau,kt-\ell\tau)_N} e^{c\lambda(\tau)\jap{\ell,\ell\tau}^s} \jap{\tau} \abs{\hat \rho_\ell(\tau)} d\tau \right]^2 dt.
\end{align*}
Since the Littlewood-Paley projections define a partition of unity (using also Cauchy-Schwarz), 
\begin{align*} 
\norm{AT}_{L_k^2 L_t^2(I)}^2  & \lesssim \sum_{k \in \Integer_\ast^d} \int_0^{T^\star} \left[\sum_{\ell \in \Integer_\ast^d}\int_0^t e^{c\lambda(\tau)\jap{\ell,\ell\tau}^s}\jap{\tau} \abs{\hat{\rho}_\ell(\tau)} d\tau \right] \\ 
& \quad\quad \times 
\left[\sum_{l\in \Integer_\ast^d} \int_0^t\abs{(A \widehat{\grad f})_{k-\ell}(\tau,kt-\ell\tau)}^2 e^{c\lambda(\tau)\jap{\ell,\ell\tau}^s} \jap{\tau}\abs{\hat \rho_\ell(\tau)} d\tau \right] dt. 
\end{align*} 
By Cauchy-Schwarz and $\sigma > \frac{d}{2}+2$, 
\begin{align} 
\sum_{\ell \in \Integer_\ast^d}\int_0^t e^{c\lambda(\tau)\jap{\ell,\ell\tau}^s}\jap{\tau} \abs{\hat{\rho}_\ell(\tau)} d\tau 
 & \leq \left(\int_0^t\sum_{\ell \in \Integer_\ast^d} e^{2(c-1)\lambda(\tau)\jap{\ell,\ell\tau}^s} \jap{\tau}^2\jap{\ell,\ell\tau}^{-2\sigma} d\tau \right)^{1/2} \norm{A\hat{\rho}}_{L_k^2 L^2_t(I)} \nonumber \\ 
& \lesssim \norm{A\hat{\rho}}_{L_k^2 L_t^2(I)}.  \label{ineq:lowrhoctrl}
\end{align}
Then \eqref{ineq:lowrhoctrl} and Fubini's theorem imply
\begin{align*}  
\norm{AT}_{L_k^2 L_t^2(I)} & \lesssim \norm{A\hat{\rho}}_{L^2_kL_t^2(I)} \sum_{k \in \Integer_\ast^d} \int_0^{T^\star} \sum_{\ell \in \Integer_\ast^d} \int_0^t\abs{(A \widehat{\grad f})_{k-\ell}(\tau,kt-\ell\tau)}^2 e^{c\lambda(\tau)\jap{\ell,\ell\tau}^s} \jap{\tau} \abs{\hat \rho_\ell(\tau)} d\tau dt \\ 
& \lesssim \norm{A\hat{\rho}}_{L_k^2L_t^2(I)}\sum_{\ell \in \Integer_\ast^d} \int_0^{T^\star} \left(\sum_{k \in \Integer_\ast^d} \int_\tau^{T^\star}\abs{(A \widehat{\grad f})_{k-\ell}(\tau,kt-\ell\tau)}^2 dt \right) e^{c\lambda(\tau)\jap{\ell,\ell\tau}^s} \jap{\tau} \abs{\hat{\rho}_\ell(\tau)} d\tau \\ 
& \leq \norm{A\hat{\rho}}_{L_k^2L_t^2(I)}\sum_{\ell \in \Integer_\ast^d} \int_0^{T^\star} \left(\sum_{k \in \Integer_\ast^d} \int_{-\infty}^\infty \abs{(A \widehat{\grad f})_{k-\ell}(\tau,kt - \ell\tau)}^2 dt \right) e^{c\lambda(\tau)\jap{\ell,\ell\tau}^s} \jap{\tau} \abs{\hat{\rho}_\ell(\tau)} d\tau \\ 
& \leq \norm{A\hat{\rho}}_{L_k^2L_t^2(I)} \left(\sup_{\tau \geq 0}e^{(c-1)\alpha_0\jap{\tau}^s} \sum_{k \in \Integer^d}\sup_{\omega \in \Integer_\ast^d} \sup_{x \in \Real^d} \int_{-\infty}^\infty \abs{(A \widehat{\grad f})_{k}\left(\tau,\frac{\omega}{\abs{\omega}}\zeta - x\right)}^2 d\zeta \right) \\ & \quad\quad \times \sum_{\ell \neq 0} \int_0^{T^\star} e^{\lambda(\tau)\jap{l,l\tau}^s} \jap{\tau} \abs{\hat{\rho}_\ell(\tau)} d\tau.
\end{align*}  
Proceeding as in \eqref{ineq:lowrhoctrl} then gives 
\begin{align*} 
\norm{AT}_{L_k^2 L_t^2(I)} & \lesssim \norm{A\hat{\rho}}_{L_k^2L_t^2(I)}^2 \left(\sup_{\tau \geq 0}e^{(c-1)\alpha_0\jap{\tau}^s}  \sum_{k \in \Integer^d}\sup_{\omega \in \Integer_\ast^d}\sup_{x \in \Real^d}\int_{-\infty}^\infty \abs{(A \widehat{\grad f})_{k}\left(\tau,\frac{\omega}{\abs{\omega}}\zeta-x\right)}^2 d\zeta \right). 
\end{align*}
Using condition \eqref{ineq:controlLines}, we derive
\begin{align} 
\norm{AT}_{L_k^2 L_t^2(I)} & \lesssim_{\alpha_0} K_1 \epsilon^2\norm{A\hat{\rho}}_{L_k^2L^2_t(I)}^2, \label{ineq:rhoTransConc}
\end{align}
which suffices to treat transport.

\subsubsection{Remainders} 
We treat the remainder with a variant of the method used to treat transport. 
First, by \eqref{ineq:Wbd}:
\begin{align*} 
\norm{A\mathcal{R}}^2_{L_k^2L_t^2(I)}
 & \lesssim \sum_{k \in \Integer_\ast^d} \int_0^{T^\ast}\left[A_k(t,kt)\int_0^t\sum_{\ell \in \Integer_\ast^d}\sum_{N \in \mathbb D} \sum_{N/8 \leq N^\prime \leq 8N} \abs{\hat{\rho}_\ell(\tau)_{N^\prime}} \abs{k(t-\tau)} \abs{f_{k-\ell}(\tau,kt-\ell\tau)_{N}} d\tau \right]^2 dt. 
\end{align*} 
Next we claim that on the integrand there holds for some $c^\prime = c^\prime(s) \in (0,1)$,  
\begin{align}
A_k(t,kt) \lesssim_{\lambda_0,\alpha_0} e^{c^\prime\lambda(\tau)\jap{k-\ell,kt-\ell\tau}^s} e^{c^\prime\lambda(\tau)\jap{\ell,\ell\tau}^s}. \label{ineq:AbdRemainder}
\end{align}  
Indeed, this follows simply by following the argument used to deduce \eqref{ineq:BgainProduct}, with $kt$ replacing $\eta$. 
 Therefore by \eqref{ineq:AbdRemainder}, Cauchy-Schwarz, 
\begin{align*} 
\norm{A\mathcal{R}}^2_{L_k^2L_t^2(I)} & \lesssim \sum_{k \in \Integer_\ast^d} \int_0^{T^\ast}\left[\int_0^t\sum_{\ell \in \Integer_\ast^d} \left(\sum_{N \in \mathbb{D} }e^{2\lambda(\tau)\jap{\ell,\ell\tau}^s}\jap{\tau}^2\abs{\hat{\rho}_\ell(\tau)_{\sim N}}^2\right)d\tau \right] \\ 
& \quad\quad \times \left[\int_0^t\sum_{\ell \in \Integer_\ast^d} \left(\sum_{N \in \mathbb D}e^{2(c^{\prime}-1)\lambda(\tau)\jap{\ell,\ell\tau}^s}e^{2c^{\prime}\lambda(\tau)\jap{k-\ell,kt-\ell\tau}^s}\abs{\widehat{\grad f}_{k-\ell}(\tau,kt-\ell\tau)_{N}}^2 \right) d\tau \right] dt. 
\end{align*}
By the almost orthogonality of the Littlewood-Paley decomposition \eqref{ineq:GeneralOrtho} and $\sigma > 1$, 
\begin{align*} 
\norm{A\mathcal{R}}^2_{L_k^2L_t^2(I)} & \lesssim \\ & \hspace{-1cm} \norm{A\hat{\rho}}^2_{L_k^2L_t^2(I)}\sum_{k \in \Integer_\ast^d} \int_0^{T^\ast} \int_0^t\sum_{\ell \in \Integer_\ast^d} \left(\sum_{N \in \mathbb D}e^{2(c^{\prime}-1)\lambda(\tau)\jap{\ell,\ell\tau}^s}e^{2c^{\prime}\lambda(\tau)\jap{k-\ell,kt-\ell\tau}^s}\abs{f_{k-\ell}(\tau,kt-\ell\tau)_{N}}^2 \right) d\tau dt. 
\end{align*} 
By Fubini's theorem and Lemma \ref{lem:SobTrace}, 
\begin{align*} 
\norm{A\mathcal{R}}^2_{L_k^2L_t^2(I)} & \lesssim \norm{A\hat{\rho}}^2_{L_k^2L_t^2(I)} \int_0^{T^\star}\sum_{\ell \in \Integer_\ast^d} e^{2(c^{\prime}-1)\lambda(\tau)\jap{\ell,\ell\tau}^s} \left(\sum_{N \in \mathbb D} \sum_{k \in \Integer_\ast^d} \int_\tau^{T^\star} \abs{A\widehat{\grad f}_{k-\ell}(\tau,kt-\ell\tau)_{N}}^2  dt\right) d\tau \\ 
& \lesssim \norm{A\hat{\rho}}^2_{L_k^2L_t^2(I)}\int_0^{T^\star}\sum_{\ell \in \Integer_\ast^d} e^{2(c^{\prime}-1)\lambda(\tau)\jap{\ell,\ell\tau}^s} \\ & \quad\quad \times  \left(\sum_{N \in \mathbb D} \sum_{k \in \Integer^d} \sup_{\omega \in \Integer_\ast^d} \sup_{x\in \Real^d} \int_{-\infty}^{\infty} \abs{A \widehat{\grad f}_{k}\left(\tau,\frac{\omega}{\abs{\omega}} \zeta-x\right)_{N}}^2  d\zeta \right) d\tau \\  
& \lesssim \norm{A\hat{\rho}}^2_{L_k^2L_t^2(I)} \int_0^{T^\star}\sum_{\ell \in \Integer_\ast^d} e^{2(c^{\prime}-1)\lambda(\tau)\jap{\ell,\ell\tau}^s} \left(\sum_{N \in \mathbb D} \sum_{k \in \Integer^d} \norm{A^{(1)}\hat{f}_k(\tau)_N}^2_{H^M(\Real^d_\eta)} \right) d\tau. 
\end{align*} 
The Littlewood-Paley projections do not commute with derivatives in frequency space, however since the projections have bounded derivatives we still have (see \S\ref{Apx:LPProduct}), 
\begin{align*} 
\norm{A\mathcal{R}}^2_{L_k^2L_t^2(I)} & \lesssim_M \norm{A\hat{\rho}}^2_{L_k^2L_t^2(I)} \int_0^{T^\star}\sum_{\ell \in \Integer_\ast^d} e^{2(c^{\prime}-1)\lambda(\tau)\jap{\ell,\ell\tau}^s} \left(\sum_{N \in \mathbb D} \sum_{\abs{\alpha} \leq M} \norm{(D_\eta^\alpha A^{(1)} \hat{f}_k)(\tau)_{\sim N}}^2_{L_k^2L^2_\eta} \right) d\tau. 
\end{align*}
Then by the almost orthogonality \eqref{ineq:GeneralOrtho} with \eqref{ineq:OuterActrl}, \eqref{ineq:SobExp} and $c^{\prime} < 1$, we have 
\begin{align} 
\norm{A\mathcal{R}}^2_{L_k^2L_t^2(I)} & \lesssim  K_1\epsilon^2 \norm{A\hat{\rho}}^2_{L_k^2L_t^2(I)} \int_0^{T^\star}\sum_{\ell \in \Integer_\ast^d} e^{2(c^{\prime}-1)\lambda(\tau)\jap{\ell,\ell\tau}^s}\jap{\tau}^7 d\tau \nonumber \\ 
& \lesssim K_1\epsilon^2 \norm{A\hat{\rho}}^2_{L_k^2L_t^2(I)}, \label{ineq:rhoRemainderConc}
\end{align} 
which suffices to treat remainder contributions. 

\subsubsection{Conclusion of $L^2$ bound} 
Putting \eqref{ineq:rhoL2ID}, \eqref{ineq:rhoReacConc}, \eqref{ineq:rhoTransConc} and \eqref{ineq:rhoRemainderConc} together with \eqref{ineq:MomentInPf} we have for some $\tilde K = \tilde K(s,d,M,\lambda_0,\alpha_0)$, 
\begin{align*}
\norm{A\hat{\rho}}^2_{L_k^2L_t^2(I)} \leq \tilde K C^2_{LD}\epsilon^2 + \tilde K C^2_{LD}(K_1 + K_2)\epsilon^2\norm{A\hat{\rho}(t)}^2_{L_k^2L_t^2(I)}.
\end{align*} 
Therefore for $\epsilon^2 < \frac{1}{2}(\tilde K C^2_{LD}(K_1 + K_2))^{-1}$ we have
\begin{align*} 
\norm{A\hat{\rho}(t)}^2_{L_k^2L_t^2(I)} < 2 \tilde K C^2_{LD} \epsilon^2. 
\end{align*}  
Hence, Proposition \ref{prop:nonlinearmoment} follows provided we fix $K_3 = \tilde K C^2_{LD}$.  
\end{proof} 

\subsection{Pointwise-in-time estimate on the density} \label{sec:PtwiseRho}
The constant $K_3$ basically only depends on the linearized Vlasov equation with homogeneous background $f^0$. 
The same is not true of the pointwise-in-time estimate we deduce next.
\begin{lemma}[Pointwise estimate] \label{lem:Ptwise}
For $\epsilon_0$ sufficiently small, under the bootstrap hypotheses \eqref{ctrl:Boot}, there exists some $K_4 = K_4(C_0,\bar{\lambda},\kappa,M,s, d,\lambda_0,\lambda^\prime,K_1,K_2,K_3)$ such that for $t \in [0,T^\star]$, 
\begin{align} 
\norm{A\rho(t)}_2^2 & \leq K_4\jap{t}\epsilon^2. \label{ctrl:MidPt}
\end{align} 
\end{lemma} 
\begin{proof} 
As in \cite{MouhotVillani11}, we use the $L_t^2$ bound together with \eqref{eq:rhointegral}.
Our starting point is again the paraproduct decomposition \eqref{def:pararho}:
\begin{align} 
\norm{A\hat{\rho}(t)}_{L_k^2} & \lesssim  \sum_{k \in \Integer_\ast^d} \abs{A_k(t,kt)\hat{h}_{in}(k,kt)}^2 \nonumber \\ & \quad + \sum_{k \in \Integer_\ast^d} \left(A_k(t,kt)\int_0^t \hat{\rho}_k(\tau)\abs{k}^2 \widehat{W}(k)(t-\tau)\hat{f}^0(k(t-\tau))  d\tau\right)^{2} \nonumber \\ 
& \quad + \sum_{k \in \Integer_\ast^d} \abs{A_k(t,kt)T_k(t)}^2 + \sum_{k \in \Integer_\ast^d} \abs{A_k(t,kt)R_k(t)}^2 + \sum_{k \in \Integer_\ast^d} \abs{A_k(t,kt)\mathcal{R}_k(t)}^2. \label{ineq:RhoPtwiseBd}
\end{align} 
To treat the initial data we use the $H^{d/2+} \hookrightarrow C^0$ embedding and that $\lambda(t)$ is decreasing \eqref{def:lambda}:
\begin{align} 
\sum_{k \in \Integer_\ast^d} \abs{A_k(t,kt)\hat{h}_{in}(k,kt)}^2 & \leq \sum_{k \in \Integer_\ast^d} \sup_{\eta \in \Real^d} \abs{A_k(t,\eta)h_{in}(k,\eta)}^2  \lesssim \norm{A(0)\hat{h}_{in}}^2_{L^2_k H_\eta^M} \lesssim \epsilon^2,  \label{ineq:rhoPtID}
\end{align} 
where we used an argument analogous to \eqref{ineq:OuterMctrl} to deduce the last inequality. 

\subsubsection{Linear contribution}
Next consider the term in \eqref{ineq:RhoPtwiseBd} coming from the homogeneous background. 
By \eqref{ineq:Wbd}, \eqref{ineq:JapTri}, \eqref{ineq:SobExp}, and the $H^{d/2+}\hookrightarrow C^0$ embedding with \eqref{ineq:IncExp}, \eqref{ineq:f0Loc} and \eqref{ctrl:Mid},  
\begin{align} 
\sum_{k \in \Integer_\ast^d} \left(A_k(t,kt)\int_0^t \hat{\rho}_k(\tau)\abs{k}^2 \widehat{W}(k)(t-\tau)\hat{f}^0(k(t-\tau))  d\tau\right)^{2} & 
\nonumber  \\ & \hspace{-7cm} \lesssim \sum_{k \in \Integer_\ast^d} \left(\int_0^t A_k(\tau,k\tau)\abs{\hat{\rho}_k(\tau)} \jap{k(t-\tau)}^{\sigma+1} e^{\lambda(t)\jap{k(t-\tau)}^s} \abs{\hat{f}^0(k(t-\tau))}  d\tau\right)^{2} \nonumber \\
& \hspace{-7cm} \lesssim \left(\sup_\eta e^{\lambda_0\jap{\eta}^s} \abs{\hat{f}^0(\eta)}\right)^2 \sum_{k \in \Integer_\ast^d} \left(\int_0^t A_k(\tau,k\tau)\abs{\hat{\rho}_k(\tau)} e^{\frac{1}{2}(\lambda(0)-\lambda_0)\jap{t - \tau}^s} d\tau \right)^2 \nonumber \\ 
& \hspace{-7cm} \lesssim C_0^2\left(\sum_{k \in \Integer_\ast^d} \int_0^t \abs{A_k(\tau,k\tau)\hat{\rho}_k(\tau)}^2 d\tau\right) \left(\int_0^t e^{(\lambda(0)-\lambda_0)\jap{t - \tau}^s} d\tau \right) \nonumber \\ 
& \hspace{-7cm} \lesssim \sum_{k \in \Integer_\ast^d }\norm{A\hat{\rho}_k}^2_{L^2_t(I)} \nonumber \\ 
& \hspace{-7cm} \lesssim K_3\epsilon^2. \label{ineq:rhoPtBackground}
\end{align} 

\subsubsection{Reaction}
Next we treat the reaction term in \eqref{ineq:RhoPtwiseBd}, which by \eqref{ineq:Wbd}, 
\begin{align*} 
\sum_{k \in \Integer_\ast^d} \abs{A_k(kt)R_k}^2 & \lesssim \sum_{k \in \Integer_\ast^d} \left[A_k(t,kt)\int_0^t\sum_{\ell \in \Integer_\ast^d} \sum_{N \geq 8}\abs{\hat{f}_{k-\ell}(\tau,kt-\ell\tau)_{<N/8}} \frac{\abs{k(t-\tau)}}{\abs{\ell}^{\gamma}} \abs{\hat{\rho}_\ell(\tau)_N} d\tau \right]^{2}. 
\end{align*} 
As in the $L_t^2$ estimate we have by \eqref{ineq:ReacRhoA} and the definition \eqref{def:barK} of $\bar K$: 
\begin{align*} 
\sum_{k \in \Integer_\ast^d} \abs{A_k(kt)R_k(t)}^2
& \lesssim \sum_{k \in \Integer_\ast^d} \left[\int_0^t\sum_{\ell \in \Integer_\ast^d} \bar{K}_{k,\ell}(t,\tau) \abs{A_\ell(\tau,\ell\tau) \hat{\rho}_\ell(\tau)}d\tau \right]^{2}. 
\end{align*} 
By Cauchy-Schwarz and Fubini's theorem, 
\begin{align*} 
\sum_{k \in \Integer_\ast^d} \abs{A_k(kt)R_k}^2 & \lesssim \sum_{k \in \Integer_\ast^d} \left(\int_0^t\sum_{\ell \in \Integer_\ast^d} \bar{K}_{k,\ell}(t,\tau) d\tau\right)  \left( \sum_{\ell \in \Integer_\ast^d} \int_0^t \bar{K}_{k,\ell}(t,\tau)\abs{A_\ell(\tau,\ell\tau) \hat{\rho}_\ell(\tau)}^2 d\tau \right) \\ 
& \lesssim \left(\sup_{t \geq 0}\sup_{k \in \Integer_\ast^d}\int_0^t\sum_{\ell \in \Integer_\ast^d } \bar{K}_{k,\ell}(t,\tau) d\tau\right) \sum_{\ell \in \Integer_\ast^d} \int_0^t \left( \sum_{k \in \Integer_\ast^d}  \bar{K}_{k,\ell}(t,\tau) \right)\abs{A_\ell(\tau,\ell\tau) \hat{\rho}_\ell(\tau)}^2 d\tau  \\ 
& \lesssim \left(\sup_{t \geq 0}\sup_{k \in \Integer_\ast^d}\int_0^t\sum_{\ell \in \Integer_\ast^d} \bar{K}_{k,\ell}(t,\tau) d\tau\right)\left( \sup_{0 \leq \tau \leq t} \sup_{\ell \in \Integer_\ast^d} \sum_{k \in \Integer_\ast^d}  \bar{K}_{k,\ell}(t,\tau) \right) \norm{A\hat{\rho}}_{L^2_k L^2_t(I)}.  
\end{align*} 
The first factor appears in \eqref{ineq:momentdef} and is controlled by Lemma \ref{lem:Moment}. 
The second factor is controlled by Lemma \ref{lem:ptwiseTimeResponse} and results in the power of $\jap{t}$ loss. 
Therefore by \eqref{ctrl:Mid},  
\begin{align} 
\sum_{k \in \Integer_\ast^d} \abs{A_k(t,kt)R_k(t)}^2 \lesssim K_2 K_3 \jap{t} \epsilon^4, \label{ineq:rhoPtReac}
\end{align} 
which suffices to treat the reaction term. 

\subsubsection{Transport} 
By \eqref{ineq:Wbd} and $\abs{k(t-\tau)} \leq \jap{\tau}\abs{k-\ell,kt-\ell\tau}$, the transport term is bounded by
\begin{align*} 
\sum_{k \in \Integer_\ast^d} \abs{A_k(t,kt)T_k(t)}^2 & \lesssim \sum_{k \in \Integer_\ast^d} \left[A_k(t,kt)\int_0^t\sum_{\ell \in \Integer_\ast^d} \sum_{N \geq 8} \abs{\widehat{\grad_{z,v} f}_{k-\ell}(\tau,kt-\ell\tau)_{N}} \jap{\tau}\abs{\hat \rho_\ell(\tau)_{<N/8}} d\tau \right]^{2}. 
\end{align*} 
We begin as in Proposition \ref{prop:nonlinearmoment}. By the frequency localizations \eqref{ineq:transptRhoLoc} (which hold on the support of the integrand), \eqref{lem:scon2} implies that for some $c = c(s) \in (0,1)$ we have (using also that the Littlewood-Paley projections define a partition of unity), 
\begin{align*} 
\sum_{k \in \Integer_\ast^d} \abs{A_k(t,kt)T_k(t)}^2 & \lesssim \sum_{k \in \Integer_\ast^d} \left[\sum_{\ell \in \Integer_\ast^d} \int_0^t \abs{(A \widehat{\grad f})_{k-\ell}(\tau,kt-\ell\tau)} e^{c\lambda(\tau)\jap{\ell,\ell\tau}^s} \jap{\tau} \abs{\hat \rho_\ell(\tau)} d\tau \right]^2. 
\end{align*} 
From Cauchy-Schwarz, \eqref{ineq:lowrhoctrl} and \eqref{ctrl:Mid},  
\begin{align*} 
\sum_{k \in \Integer_\ast^d} \abs{A_k(t,kt)T_k(t)}^2 & \lesssim \sqrt{K_3}\epsilon \sum_{k \in \Integer_\ast^d}\sum_{\ell\in \Integer_\ast^d} \int_0^t\abs{(A \widehat{\grad f})_{k-\ell}(\tau,kt-\ell\tau)}^2 e^{c\lambda(\tau)\jap{\ell,\ell\tau}^s} \jap{\tau} \abs{\hat \rho_\ell(\tau)} d\tau.  
\end{align*} 
By Fubini's theorem, \eqref{ineq:SobExp}, the $H^{d/2+} \hookrightarrow C^0$ embedding theorem and \eqref{ineq:OuterActrl} with \eqref{ineq:SobExp},  
\begin{align*} 
\sum_{k \in \Integer_\ast^d} \abs{A_k(t,kt)T_k(t)}^2 & \lesssim \sqrt{K_3}\epsilon \sum_{\ell \in \Integer_\ast^d} \int_0^t\left(\sum_{k \in \Integer_\ast^d}\abs{(A \widehat{\grad f})_{k-\ell}(\tau,kt-\ell\tau)}^2 e^{\frac{1}{2}(c-1)\alpha_0\jap{\tau}^s}\right) e^{\lambda(\tau)\jap{\ell,\ell\tau}^s} \abs{\hat \rho_\ell(\tau)} d\tau \\ 
& \hspace{-1cm}\lesssim \sqrt{K_3}\epsilon \left(\sup_{\tau \leq t}e^{\frac{1}{2}(c-1)\alpha_0\jap{\tau}^s}\sum_{k \in \Integer^d} \sup_{\eta \in \Real^d} \abs{(A \widehat{\grad f})_{k}(\tau,\eta)}^2\right) \left(\sum_{\ell \in \Integer_\ast^d}  \int_0^t e^{\lambda(\tau)\jap{\ell,\ell\tau}^s}\abs{\hat \rho_\ell(\tau)} d\tau \right) \\ 
& \hspace{-1cm} \lesssim \sqrt{K_3}\epsilon \left(\sup_{\tau \leq t}e^{\frac{1}{2}(c-1)\alpha_0\jap{\tau}^s}\norm{A^{(1)}\hat{f}}^2_{L^2_k H^M_\eta}\right) \left(\sum_{\ell\in \Integer_\ast^d}  \int_0^t e^{\lambda(\tau)\jap{\ell,\ell\tau}^s}\abs{\hat \rho_\ell(\tau)} d\tau \right) \\ 
& \hspace{-1cm} \lesssim K_1\sqrt{K_3}\epsilon^3 \left(\sum_{\ell\in \Integer_\ast^d}  \int_0^t e^{\lambda(\tau)\jap{\ell,\ell\tau}^s}\abs{\hat \rho_\ell(\tau)} d\tau \right). 
\end{align*} 
Proceeding as in \eqref{ineq:lowrhoctrl} and applying \eqref{ctrl:Mid}, we get
\begin{align} 
\sum_{k \in \Integer_\ast^d} \abs{A_k(t,kt)T_k(t)}^2 \lesssim K_3K_1 \epsilon^4. \label{ineq:rhoPtTrans}
\end{align} 

\subsubsection{Remainders}
The remainder follows from a slight variant of the argument used to treat transport. 
By \eqref{ineq:Wbd},
\begin{align*} 
\norm{A\mathcal{R}(t)}_{L_k^2}^2 & \lesssim \sum_{k \in \Integer_\ast^d}\left[A_k(t,kt)\int_0^t\sum_{\ell \in \Integer_\ast^d}\sum_{N \in \mathbb D} \sum_{N/8 \leq N^\prime \leq 8N} \abs{\hat{\rho}_\ell(\tau)_{N^\prime}} \abs{k(t-\tau)} \abs{\hat{f}_{k-\ell}(\tau,kt-\ell\tau)_{N}} d\tau \right]^2.  
\end{align*} 
As in Proposition \ref{prop:nonlinearmoment}, \eqref{ineq:AbdRemainder} holds on the support of the integrand and hence, by \eqref{ineq:Echot} and Cauchy-Schwarz, 
\begin{align*} 
\norm{A\mathcal{R}(t)}_{L_k^2}^2 & \lesssim \sum_{k \in \Integer_\ast^d} \left[\int_0^t\sum_{\ell \in \Integer_\ast^d} \left(\sum_{N \in \mathbb{D} }e^{2\lambda(\tau)\jap{\ell,\ell\tau}^s}\jap{\tau}^2\abs{\hat{\rho}_\ell(\tau)_{\sim N}}^2\right)d\tau \right] \\ 
& \quad\quad \times \left[\int_0^t\sum_{\ell \in \Integer_\ast^d} \left(\sum_{N \in \mathbb D}e^{2(c^{\prime}-1)\lambda(\tau)\jap{\ell,\ell\tau}^s}e^{2c^{\prime}\lambda(\tau)\jap{k-\ell,kt-\ell\tau}^s}\abs{\widehat{\grad f}_{k-\ell}(\tau,kt-\ell\tau)_{N}}^2 \right) d\tau \right]. 
\end{align*} 
By the almost orthogonality of the Littlewood-Paley decomposition \eqref{ineq:GeneralOrtho} and $\sigma > 1$, 
\begin{align*} 
\norm{A\mathcal{R}(t)}^2_{L_k^2} & \lesssim \norm{A\hat{\rho}}^2_{L_k^2L_t^2(I)}\sum_{k \in \Integer_\ast^d} \int_0^t\sum_{\ell \in \Integer_\ast^d} \left(\sum_{N \in \mathbb D}e^{2(c^{\prime}-1)\lambda(\tau)\jap{\ell,\ell\tau}^s}e^{2c^{\prime}\lambda(\tau)\jap{k-\ell,kt-\ell\tau}^s} \abs{\widehat{\grad f}_{k-\ell}(\tau,kt-\ell\tau)_{N}}^2 \right) d\tau. 
\end{align*} 
By Fubini's theorem, the $H^{d/2+} \hookrightarrow C^0$ embedding theorem and \eqref{ineq:SobExp}, 
\begin{align*} 
\norm{A\mathcal{R}(t)}^2_{L_k^2} & \lesssim \norm{A\hat{\rho}}^2_{L_k^2L_t^2(I)}\sum_{\ell \in \Integer_\ast^d} \int_0^t e^{2(c^{\prime}-1)\lambda(\tau)\jap{\ell,\ell\tau}^s}\left( \sum_{k \in \Integer_\ast^d} \sum_{N \in \mathbb D} e^{2c^{\prime}\lambda(\tau)\jap{k-\ell,kt-\ell\tau}^s}\abs{\widehat{\grad f}_{k-\ell}(\tau,kt-\ell\tau)_{N}}^2 \right) d\tau \\ 
& \lesssim \norm{A\hat{\rho}}^2_{L_k^2L_t^2(I)}\sum_{\ell \in \Integer_\ast^d} \int_0^t e^{2(c^{\prime}-1)\lambda(\tau)\jap{\ell,\ell\tau}^s} \left( \sum_{k \in \Integer^d} \sum_{N \in \mathbb D} \norm{(A^{(-\beta)}\hat{f}_{k})(\tau)_{N}}_{H^M_\eta}^2 \right) d\tau. 
\end{align*} 
The Littlewood-Paley projections do not commute with derivatives in frequency space, however since the projections have bounded derivatives we still have (see \S\ref{Apx:LPProduct}), 
\begin{align*} 
\norm{A\mathcal{R}(t)}^2_{L_k^2} & \lesssim_M \norm{A\hat{\rho}}^2_{L_k^2L_t^2(I)} \sum_{\ell \in \Integer_\ast^d} \int_0^t  e^{2(c^{\prime}-1)\lambda(\tau)\jap{\ell,\ell\tau}^s} \left(\sum_{N \in \mathbb D} \sum_{\abs{\alpha} \leq M} \norm{(D_\eta^\alpha A^{(-\beta)} \hat{f}_k)(\tau)_{\sim N}}^2_{L_k^2L^2_\eta} \right) d\tau. 
\end{align*} 
Hence by \eqref{ineq:GeneralOrtho}, \eqref{ineq:OuterMctrl}, \eqref{ctrl:Mid} and \eqref{ineq:SobExp}, 
\begin{align} 
\norm{A\mathcal{R}(t)}^2_{L_k^2} & \lesssim K_3 K_2 \epsilon^4. \label{ineq:rhoPtRemainder}
\end{align} 

Summing \eqref{ineq:rhoPtID}, \eqref{ineq:rhoPtBackground}, \eqref{ineq:rhoPtReac}, \eqref{ineq:rhoPtTrans} and \eqref{ineq:rhoPtRemainder} implies the result with $K_4 \approx 1 + K_3 + K_2K_3 + K_3K_1$ (in fact we are rather suboptimal). 
\end{proof} 

\subsection{Proof of high norm estimate \eqref{ctrl:HiLocalizedB}} \label{subsec:HiNorm}
In this section we derive the high norm estimate on the full distribution, \eqref{ctrl:HiLocalizedB}.  
For some multi-index $\alpha \in \Naturals^{d}$ with $\abs{\alpha} \leq M$, compute the time-derivative 
\begin{align} 
\frac{1}{2}\frac{d}{dt}\norm{A^{(1)}D_\eta^\alpha \hat f}_2^2 & = \sum_{k\in \Integer^d}\int_\eta \dot{\lambda}(t)\jap{k,\eta}^{s} \abs{A^{(1)} D_\eta^\alpha \hat f_k(\eta)}^2 d\eta +   \sum_{k \in \Integer^d} \int_\eta A^{(1)} D_\eta^\alpha \overline{\hat{f}_k(\eta)} A^{(1)} D_\eta^\alpha \partial_t \hat{f}_k(\eta) d\eta \nonumber \\ 
& = CK + E. \label{ineq:HighNormDeriv}
\end{align} 
Like similar terms appearing in \cite{LevermoreOliver97,KukavicaVicol09,BM13}, the $CK$ term (for `Cauchy-Kovalevskaya') is used to absorb the highest order terms coming from $E$. 

Turning to $E$, we separate into the linear and nonlinear contributions 
\begin{align} 
E & = -\sum_{k \in \Integer^d_\ast} \int_\eta A^{(1)}D_\eta^\alpha\overline{\hat{f}_k(\eta)} A^{(1)}_k(t,\eta)D_\eta^\alpha \left[\hat{\rho}_k(t) \widehat{W}(k)k \cdot (\eta - tk) \widehat{f^0}(\eta - kt)\right] d\eta \nonumber \\ 
& \quad -\sum_{k\in\Integer^d} \int_\eta A^{(1)}D_\eta^\alpha\overline{\hat{f}_k(\eta)} A^{(1)}_k(t,\eta)D_\eta^\alpha \left[\sum_{\ell \in \Integer_\ast^d} \hat{\rho}_\ell(t)\widehat{W}(\ell) \ell\cdot (\eta - tk) \hat f_{k-\ell}(t,\eta - t\ell)\right] d\eta \nonumber \\ 
& = -E_{L} - E_{NL}. \label{def:ELENL}
\end{align} 

\subsubsection{Linear contribution}\label{subsec:linearContr}
The linear contribution $E_L$ is easier from a regularity standpoint than $E_{NL}$ since we may lose regularity when estimating $f^0$. 
However, $E_L$ has one less power of $\epsilon$ which requires some care to handle and is the reason we cannot just take $K_1$ in \eqref{ctrl:HiLocalizedB} to be $O(1)$. 
The treatment of $E_L$ begins with the product lemma \eqref{ineq:GProduct2}: 
\begin{align*}
\abs{E_L} & \lesssim \norm{A^{(1)} D_\eta^\alpha \hat{f}}_2\norm{A^{(1)} D_\eta^\alpha(\eta\hat{f}^0(\eta))}_{L^2_\eta} \norm{\grad_x W \ast_x \rho}_{\cF^{\tilde c\lambda(t);s}} 
  \\ & \quad + \norm{A^{(1)} D_\eta^\alpha \hat{f}}_2\norm{v^\alpha(\grad_vf^0)}_{\cG^{\tilde{c} \lambda(t);s}}\norm{A^{(1)}\grad_x W \ast_x \rho(t)}_{2}. 
\end{align*} 
By \eqref{ineq:f0Loc}, \eqref{ineq:SobExp} and \eqref{ineq:IncExp}, 
\begin{align} 
\norm{A^{(1)}D^{\alpha}_\eta(\eta \hat{f}^0(\eta))}_2 & \leq \norm{A^{(1)}\left(\eta D^{\alpha}_\eta\hat{f}^0(\eta)\right)}_2 + \sum_{\abs{j} = 1;j \leq \alpha}\norm{A^{(1)} \left(\eta D^{\alpha-j}_\eta\hat{f}^0(\eta)\right)}_2 \lesssim C_0. \label{ineq:BDalphaf0}
\end{align} 
Next we use $\gamma \geq 1$ to deduce from \eqref{ineq:Wbd}, 
\begin{align*} 
A^{(1)}_k(t,kt) \abs{\widehat{W}(k)k} = \jap{k,kt} \abs{\widehat{W}(k)k} A_k(t,kt) \lesssim \jap{t} A_k(t,kt),   
\end{align*} 
which implies (also using $\tilde c < 1$ and \eqref{ineq:SobExp}), 
\begin{align} 
\abs{E_L} & \lesssim e^{(\tilde c - 1)\alpha_0\jap{t}^s}\norm{A^{(1)}D_\eta^\alpha \hat{f}}_2\norm{A\rho(t)}_2 + \jap{t}\norm{A^{(1)} D_\eta^\alpha \hat{f}}_2\norm{A\rho(t)}_2 \nonumber \\ &  \lesssim \jap{t}\norm{A^{(1)} D_\eta^\alpha \hat{f}}_2\norm{A\rho(t)}_2.  \label{ineq:ELConc}
\end{align}

\subsubsection{Commutator trick for the nonlinear term}\label{sec:NL}
Turn now to the nonlinear term in \eqref{def:ELENL}, $E_{NL}$. 
Here we cannot lose much regularity on any of the factors involved,  however we have additional powers of $\epsilon$ which will eliminate the large constants.  
First, we expand the $D_\eta^\alpha$ derivative
\begin{align*} 
E_{NL} & = \sum_{k \in \Integer^d} \int_\eta A^{(1)}D_\eta^\alpha \overline{\hat{f}_k(\eta)}\left(A^{(1)}_k(t,\eta) \left[\sum_{\ell \in \Integer_\ast^d} \hat{\rho}_\ell(t)\widehat{W}(\ell)\ell\cdot ( \eta - tk ) D_\eta^\alpha \hat f_{k-\ell}(t,\eta - t\ell)\right] \right) d\eta \\
& \quad + \sum_{k\in \Integer^d} \int_\eta A^{(1)}D_\eta^\alpha \overline{\hat{f}_k(\eta)}\left(A^{(1)}_k(t,\eta) \sum_{\abs{j} = 1; j \leq \alpha}\left[\sum_{\ell \in \Integer_\ast^d} \hat{\rho}_\ell(t)\widehat{W}(\ell)\ell_j D_\eta^{\alpha-j} \hat f_{k-\ell}(t,\eta - t\ell)\right] \right) d\eta \\ 
& = E_{NL}^1 + E_{NL}^2. 
\end{align*}
Consider first $E_{NL}^1$, as this contains an extra derivative which results in a loss of regularity that must be balanced by the $CK$ term in \eqref{ineq:HighNormDeriv}.
To gain from the cancellations inherent in transport we follow the commutator trick used in (for example) \cite{FoiasTemam89,LevermoreOliver97,KukavicaVicol09,BM13} by applying the identity, 
\begin{align} 
\frac{1}{2}\int_{\T^d \times \R^d} F(t,z+tv)\cdot (\grad_v - t\grad_z) \left[A^{(1)}(v^\alpha f)\right]^2 dz dv = 0, \label{eq:cancel}
\end{align}
to write, 
\begin{align*} 
E_{NL}^1 = \sum_{k\in\Integer^d} \sum_{\ell \in \Integer_\ast^d} \int_\eta A^{(1)} D_\eta^\alpha \overline{\hat{f}_k(\eta)} \hat{\rho}_\ell(t)\widehat{W}(\ell)\ell \cdot  (\eta - kt) \left[A^{(1)}_{k}(t,\eta) - A^{(1)}_{k-\ell}(t,\eta-t\ell)\right]\left(D_\eta^\alpha \hat f_{k-\ell}(t,\eta - t\ell) \right) d\eta.
\end{align*} 
We divide further via paraproduct:
\begin{align} 
E_{NL}^1 & = \sum_{N \geq 8} T^1_N + \sum_{N \geq 8} R^1_N + \mathcal{R}^1, \label{def:Edecomp}
\end{align} 
where the \emph{transport} term is given by  
\begin{align} 
T^1_N & = \sum_{k\in\Integer^d} \sum_{\ell \in \Integer_\ast^d} \int_\eta A^{(1)}D_\eta^\alpha \overline{\hat{f}_k(\eta)} \hat{\rho}_\ell(t)_{<N/8}\widehat{W}(\ell)\ell \cdot  (\eta - kt) \nonumber \\ & \quad\quad \times \left[A^{(1)}_{k}(t,\eta) - A^{(1)}_{k-\ell}(t,\eta-t\ell)\right]\left(D_\eta^\alpha \hat f_{k-\ell}(t,\eta - t\ell) \right)_{N} d\eta, \label{def:Tn}
\end{align} 
and the \emph{reaction} term by
\begin{align} 
R^1_N & = \sum_{k\in\Integer^d} \sum_{\ell \in \Integer_\ast^d} \int_\eta  A^{(1)}D_\eta^\alpha \overline{\hat{f}_k(\eta)} \hat{\rho}_\ell(t)_{N}\widehat{W}(\ell)\ell \cdot \left(\eta - kt\right) \nonumber \\ & \quad\quad \times \left[A^{(1)}_{k}(t,\eta) - A^{(1)}_{k-\ell}(t,\eta-t\ell)\right] \left(D_\eta^\alpha \hat f_{k-\ell}(t,\eta - t\ell) \right)_{<N/8}d\eta. \label{def:Rn}
\end{align} 
The remainder, $\mathcal{R}^1$, is whatever is left over. 

\subsubsection{Transport} \label{sec:HiNormTrans}
On the support of the integrand in \eqref{def:Tn} we have 
\begin{subequations} \label{ineq:TnFreqLoc}
\begin{align} 
\frac{N}{2} \leq \abs{k-\ell,\eta-t\ell} & \leq \frac{3N}{2}, \\ 
\abs{\ell,\ell t} & \leq \frac{3N}{32}, \\ 
\frac{13}{16} \leq \frac{\abs{k,\eta}}{\abs{k-\ell,\eta-t\ell}} & \leq \frac{19}{16}. 
\end{align}
\end{subequations} 
By \eqref{ineq:TnFreqLoc} we can gain from the multiplier: 
\begin{align} 
\abs{\frac{A^{(1)}_k(\eta)}{A^{(1)}_{k-l}(\eta-\ell t)} - 1} &= \abs{\frac{e^{\lambda \jap{k,\eta}^s} \jap{k,\eta}^{\sigma+1}}{e^{\lambda \jap{k-\ell,\eta-\ell t}^s} \jap{k-\ell,\eta-\ell t}^{\sigma+1}} - 1 } \nonumber \\
& \hspace{-2cm} \leq  \abs{e^{\lambda \jap{k,\eta}^s -\lambda \jap{k-\ell,\eta-\ell t}^s} - 1 } +  e^{\lambda \jap{k,\eta}^s -\lambda \jap{k-\ell,\eta-\ell t}^s}\abs{\frac{\jap{k,\eta}^{\sigma+1}}{\jap{k-\ell,\eta-t\ell}^{\sigma+1}} -1 }. \label{def:BkBkl}
\end{align} 
By $\abs{e^x - 1} \leq xe^x$, \eqref{ineq:TrivDiff2} and \eqref{lem:scon2} (using \eqref{ineq:TnFreqLoc}), there is some $c = c(s) \in (0,1)$ such that: 
\begin{align} 
\abs{e^{\lambda \jap{k,\eta}^s -\lambda \jap{k-\ell,\eta-\ell t}^s} - 1} & \leq \lambda\abs{\jap{k,\eta}^s - \jap{k-\ell,\eta-\ell t}^s}e^{\lambda \jap{k,\eta}^s -\lambda \jap{k-\ell,\eta-\ell t}^s} \nonumber \\ 
& \lesssim \frac{\jap{\ell,\ell t}}{\jap{k,\eta}^{1-s} +\jap{k-\ell,\eta-\ell t}^{1-s}}  e^{c\lambda\jap{\ell,\ell t}^s}. 
\label{ineq:Transgain}
\end{align} 
The other term in \eqref{def:BkBkl} can be treated with the mean-value theorem and \eqref{lem:scon2}, resulting in a bound 
not worse than \eqref{ineq:Transgain}. 
Therefore, applying \eqref{ineq:Wbd}, \eqref{def:BkBkl}, \eqref{ineq:Transgain} and adding a frequency localization by \eqref{ineq:TnFreqLoc} to $T^1_N$ implies 
\begin{align*} 
\abs{T^1_N} & \lesssim \sum_{k \in \Integer^d} \sum_{\ell \in \Integer_\ast^d} \int_\eta \abs{A^{(1)} D_\eta^\alpha \hat{f}_k(\eta)}\abs{\hat{\rho}_\ell(t)_{<N/8}} \frac{\abs{\eta - \ell t - t(k-\ell)}\jap{\ell,\ell t}}{\jap{k,\eta}^{1-s} +\jap{k-\ell,\eta-\ell t}^{1-s}}  e^{c\lambda\jap{\ell,\ell t}^s} \\ & \quad\quad\quad \times A^{(1)}_{k-\ell}(t,\eta-t\ell)\abs{\left(D_\eta^\alpha \hat f_{k-\ell}(t,\eta - t\ell)\right)_{N}}  d\eta \\ 
& \lesssim \jap{t}^2\sum_{k \in \Integer^d} \sum_{\ell \in \Integer_\ast^d} \int_\eta \abs{\left(A^{(1)} D_\eta^\alpha \hat{f}_k(\eta)\right)_{\sim N}}\abs{\hat{\rho}_\ell(t)} \jap{\ell} \abs{k-\ell,\eta - t\ell}^{s/2}\abs{k,\eta}^{s/2} e^{c\lambda\jap{\ell,\ell t}^s} \\ & \quad\quad\quad \times A^{(1)}_{k-\ell}(t,\eta-t\ell)\abs{\left(D_\eta^\alpha \hat f_{k-\ell}(t,\eta - t\ell)\right)_{N}}  d\eta.
\end{align*} 
Applying \eqref{ineq:L2L2L1} implies  
\begin{align*} 
\abs{T^1_N} & \lesssim \jap{t}^2 \norm{\abs{\grad_{z,v}}^{s/2}A^{(1)}(v^\alpha f)_{\sim N}}_2\norm{\abs{\grad_{z,v}}^{s/2}A^{(1)}(v^\alpha f)_{N}}_2 \norm{\rho(t)}_{\cF^{c\lambda(t),\frac{d}{2}+2;s}}. 
\end{align*}
Using the regularity gap provided by $c < 1$ and \eqref{ineq:SobExp} (also $\sigma > \frac{d}{2}+2$),
\begin{align} 
\abs{T^1_N} & \lesssim  \jap{t}^2 e^{(c-1)\lambda(t) \jap{t}^{s}}\norm{A\rho(t)}_2\norm{\jap{\grad_{z,v}}^{s/2}A^{(1)}(v^\alpha f)_{\sim N}}_2\norm{\jap{\grad_{z,v}}^{s/2}A^{(1)}(v^\alpha f)_{N}}_2 \nonumber \\  
& \lesssim_{\alpha_0} e^{\frac{1}{2}(c-1)\alpha_0 \jap{t}^{s}} \norm{A\rho(t)}_2
\norm{\jap{\grad_{z,v}}^{s/2}A^{(1)}(v^\alpha f)_{\sim N}}_2^2. \label{ineq:TNConc}
\end{align} 
We will find that this term is eventually absorbed by the $CK$ term in \eqref{ineq:HighNormDeriv}.

\subsubsection{Reaction} \label{sec:HiReac}
Next we consider the reaction contribution, where the commutator introduced by the identity \eqref{eq:cancel} to deal with transport will not be useful. 
Hence, write $R^1_N = R_N^{1;1} + R_N^{1;2}$ where
\begin{align*} 
R_N^{1;1} & = \sum_{k\in\Integer^d} \sum_{\ell \in \Integer_\ast^d} \int_\eta A^{(1)}D_\eta^\alpha \overline{\hat{f}_k(\eta)} A^{(1)}_{k}(t,\eta) \hat{\rho}_\ell(t)_{N}\widehat{W}(\ell)\ell \cdot \left(\eta - kt\right)  \left(D_\eta^\alpha \hat f_{k-\ell}(t,\eta - t\ell) \right)_{<N/8}d\eta.
\end{align*}
We focus on $R_N^{1;1}$ first; $R_N^{1;2}$ is easier as the norm is landing on the `low frequency' factor.
On the support of the integrand, we have the frequency localizations \eqref{ineq:FreqLocExample},
from which it follows by \eqref{lem:scon2} that there exists some $c = c(s) \in (0,1)$ such that
 \begin{align} 
\abs{R_N^{1;1}} & \lesssim \sum_{k\in\Integer^d} \sum_{\ell \in \Integer_\ast^d} \int_\eta \abs{A^{(1)}D_\eta^\alpha \hat{f}_k(\eta)} A^{(1)}_\ell(t,\ell t)\abs{\widehat{W}(\ell) \ell \hat{\rho}_\ell(t)_{N}} \nonumber \\ 
& \quad\quad \times e^{c\lambda\jap{k-\ell,\eta-t\ell}^s} \abs{\left[ \eta - tk\right] \left(D_\eta^\alpha \hat f_{k-\ell}(t,\eta - t\ell)\right)_{<N/8} }  d\eta. \label{ineq:RNeps11}
\end{align} 
Now again we have the crucial use of the assumption $\gamma \geq 1$ as in $E_L^1$: 
\begin{align*} 
A^{(1)}_\ell(t,\ell t)\abs{\widehat{W}(\ell) \ell} \lesssim  \frac{\jap{\ell,\ell t}}{\abs{\ell}} A_\ell(t,\ell t) & \lesssim \jap{t}A_\ell(t,\ell t). 
\end{align*} 
Therefore (adding a frequency localization by \eqref{ineq:FreqLocExample}), by $\abs{\eta-kt} \leq \jap{t}\abs{k-\ell,\eta-t\ell}$ and \eqref{ineq:SobExp}
\begin{align*} 
\abs{R_N^{1;1}} & \lesssim \jap{t}\sum_{k\in\Integer^d} \sum_{\ell \in \Integer_\ast^d} \int_\eta \abs{ A^{(1)}\left(D_\eta^\alpha \hat{f}_k(\eta)\right)_{\sim N}} \abs{A_\ell(t,t \ell) \hat{\rho}_\ell(t)_{N}}e^{c\lambda\jap{k-\ell,\eta-t\ell}^s} \abs{(\eta - tk) \left(D_\eta^\alpha \hat f_{k-\ell}(t,\eta - t\ell)\right)}  d\eta \\ 
& \lesssim \jap{t}^2\sum_{k\in\Integer^d} \sum_{\ell \in \Integer_\ast^d}\int_\eta \abs{A^{(1)} \left(D_\eta^\alpha \hat{f}_k(\eta)\right)_{\sim N}} \abs{A_\ell(t,t \ell) \hat{\rho}_\ell(t)_{N}} e^{\lambda\jap{k-\ell,\eta-t\ell}^s}\abs{\left(D_\eta^\alpha \hat f_{k-\ell}(t,\eta - t\ell)\right)}  d\eta.
\end{align*} 
Applying \eqref{ineq:L2L1L2} and $\sigma - \beta > \frac{d}{2}+1$, 
\begin{align} 
\abs{R_N^{1;1}} & \lesssim \jap{t}^2 \norm{A^{(1)}\left(v^\alpha f\right)_{\sim N}}_2 \norm{A\rho_N}_2 \norm{A^{(-\beta)}v^\alpha f}_2 \nonumber \\  
& \lesssim \frac{\norm{A^{(-\beta)}v^\alpha f}_2}{\jap{t}^2} \norm{A^{(1)}\left(v^\alpha f\right)_{\sim N}}^2_2  + \norm{A^{(-\beta)}v^\alpha  f}_2 \jap{t}^6 \norm{A\rho_N}^2_2, 
 \label{ineq:RN11Conc}
\end{align} 
which will suffice to treat this term.

Next turn to $R_N^{1;2}$, which is easier. By \eqref{ineq:Wbd},  
\begin{align*} 
\abs{R_N^{1;2}} & \lesssim \sum_{k\in\Integer^d} \sum_{\ell \in \Integer_\ast^d} \int_\eta \abs{A^{(1)} D_\eta^\alpha \hat{f}_k(\eta)} \abs{\hat{\rho}_\ell(t)_{N}}\abs{\eta - kt} \abs{A^{(1)}_{k-\ell}(t,\eta-t\ell) \left(D_\eta^\alpha \hat f_{k-\ell}(t,\eta - t\ell) \right)_{<N/8}}d\eta.
\end{align*}
Since the frequency localizations \eqref{ineq:FreqLocExample} hold also on the support of the integrand of $R_N^{1;2}$ (in particular  $\abs{\eta-kt} \leq \jap{t}\abs{k-\ell,\eta-t\ell} \lesssim \jap{t}\abs{\ell,\ell t}$), 
\begin{align*} 
\abs{R_N^{1;2}} 
& \lesssim \jap{t}\sum_{k\in\Integer^d} \sum_{\ell \in \Integer_\ast^d} \int_\eta\abs{A^{(1)} D_\eta^\alpha \hat{f}_k(\eta)_{\sim N}} \jap{\ell,t\ell}\abs{\hat{\rho}_\ell(t)_{N}} A^{(1)}_{k-\ell}(t,\eta-t\ell)\abs{\left(D_\eta^\alpha \hat f_{k-\ell}(t,\eta - t\ell) \right)_{<N/8}}d\eta. 
\end{align*} 
Therefore, by \eqref{ineq:L2L2L1}, \eqref{ineq:dyadicDiff}, \eqref{ineq:SobExp}, \eqref{ineq:FreqLocExample} and $\sigma > \frac{d}{2} + 3$, 
\begin{align} 
\abs{R_N^{1;2}} & \lesssim \jap{t} \norm{A^{(1)}(v^\alpha f)_{\sim N}}_2 \norm{\rho(t)_N}_{\cF^{0,\frac{d}{2}+2}} 
\norm{A^{(1)}(v^\alpha f)}_2 \nonumber \\
 & \lesssim \frac{\jap{t}}{N} \norm{A^{(1)}(v^\alpha f)_{\sim N}}_2 \norm{\rho(t)_N}_{\cF^{0,\frac{d}{2}+3}} 
\norm{A^{(1)}(v^\alpha f)}_2 \nonumber \\ 
 & \lesssim \frac{e^{-\alpha_0\jap{t}^s}}{N} \norm{A^{(1)}(v^\alpha f)}^2_2 \norm{A\rho(t)}_2, \label{ineq:RN12Conc} 
\end{align} 
which suffices to treat this term. 

\subsubsection{Remainders} \label{sec:HiNormRemainder}
In order to complete the treatment of $E_{NL}^1$ it remains to estimate the remainder $\mathcal{R}$. 
Like $R_N^{1}$, the commutator introduced by \eqref{eq:cancel} is not helpful, so divide into two pieces:
\begin{align*} 
\mathcal{R}^1 & = \sum_{N \in \mathbb D} \sum_{N/8 \leq N^\prime \leq 8N} \sum_{k\in\Integer^d} \sum_{\ell \in \Integer_\ast^d} \int_\eta A^{(1)} D_\eta^\alpha \overline{\hat{f}_k(\eta)} \hat{\rho}_\ell(t)_{N^\prime}\widehat{W}(\ell)\ell \cdot  (\eta - kt) \\ & \quad\quad \times \left[A^{(1)}_{k}(t,\eta) - A^{(1)}_{k-\ell}(t,\eta-t\ell)\right]\left(D_\eta^\alpha \hat f_{k-\ell}(t,\eta - t\ell) \right)_{N} d\eta \\ 
& = \mathcal{R}^{1;1} + \mathcal{R}^{1;2}.
\end{align*} 
Analogous to \eqref{ineq:AbdRemainder}, we claim that on the integrand there holds for some $c^\prime = c^\prime(s) \in (0,1)$,  
\begin{align}
A^{(1)}_k(t,\eta) \lesssim_{\lambda_0,\alpha_0} e^{c^\prime\lambda(t)\jap{k-\ell,\eta-t\ell}^s} e^{c^\prime\lambda(t)\jap{\ell,\ell t}^s}, \label{ineq:RemainGain}
\end{align}  
which again follows by the argument used to deduce \eqref{ineq:BgainProduct}. 
Therefore, \eqref{ineq:RemainGain}  implies (using also \eqref{ineq:Wbd} and $\abs{\eta-kt} \leq \jap{t}\jap{k-\ell,\eta - t\ell}$),
\begin{align*} 
\abs{\mathcal{R}^{1;1}} & \lesssim  \jap{t}\sum_{N  \in \mathbb D} \sum_{N^\prime \approx N} \sum_{k\in\Integer^d} \sum_{\ell \in \Integer_\ast^d} \int_\eta \abs{A^{(1)} D_\eta^\alpha \hat{f}_k(\eta)} e^{c^\prime \lambda(t)\jap{\ell,\ell t}^s}\abs{\hat{\rho}_\ell(t)_{N^\prime}} \\ & \quad\quad\quad \times \jap{k-\ell,\eta-t\ell} e^{c^\prime \lambda(t)\jap{k-\ell,\eta-t\ell}^s}\abs{\left(D_\eta^\alpha \hat f_{k-\ell}(t,\eta - t\ell) \right)_{N}} d\eta.
\end{align*} 
Applying \eqref{ineq:L2L1L2}, $\sigma>\frac{d}{2}+2$, \eqref{ineq:dyadicDiff} and \eqref{ineq:LPOrthoProject}, we have
\begin{align} 
\abs{\mathcal{R}^{1;1}} & \lesssim e^{\frac{1}{2}(c^\prime - 1)\alpha_0\jap{t}^s}\sum_{N  \in \mathbb D} \sum_{N^\prime \approx N}\norm{A^{(1)}(v^\alpha f)}_2 \frac{1}{N^\prime} \norm{A \rho(t)_{N^\prime}}_2 \norm{(v^\alpha f)_{N}}_{\cG^{c^\prime\lambda(t);\frac{d}{2}+2;s}} \nonumber \\   
& \lesssim \norm{A\rho(t)}_2 e^{\frac{1}{2}(c^\prime - 1)\alpha_0\jap{t}^s}\sum_{N  \in \mathbb D}\norm{A^{(1)}(v^\alpha f)}_2 \frac{1}{N} \norm{A^{(1)}(v^\alpha f)_{N}}_2 \nonumber  \\ 
& \lesssim \norm{A\rho(t)}_2\jap{t}^{-1}e^{\frac{1}{4}(c^\prime - 1)\alpha_0\jap{t}^s} \norm{A^{(1)}(v^\alpha f)}^2_2, 
\label{ineq:RRaHConc}
\end{align} 
which will suffice to treat this term. 

Treating $\mathcal{R}^{1;2}$ is very similar to $\mathcal{R}^{1;1}$. 
Indeed, on the support of the integrand $\jap{k-\ell,\eta-t\ell} \lesssim \jap{\ell,t\ell}$ by the same logic used to deduce \eqref{ineq:RemainGain} and hence \eqref{ineq:Wbd} and \eqref{ineq:SobExp} imply 
\begin{align*} 
\abs{\mathcal{R}^{1;2}} & \lesssim  \jap{t}\sum_{N  \in \mathbb D} \sum_{N^\prime \approx N} \sum_{k\in\Integer^d} \sum_{\ell \in \Integer_\ast^d} \int_\eta \abs{A^{(1)} D_\eta^\alpha \hat{f}_k(\eta)} \abs{\hat{\rho}_\ell(t)_{N^\prime}} \\ & \quad\quad\quad \times \jap{k-\ell,\eta-t\ell} A^{(1)}_{k-\ell}(t,\eta-\ell t)\abs{\left(D_\eta^\alpha \hat f_{k-\ell}(t,\eta - t\ell) \right)_{N}} d\eta \\ 
& \lesssim  \jap{t}\sum_{N  \in \mathbb D} \sum_{N^\prime \approx N} \sum_{k\in\Integer^d} \sum_{\ell \in \Integer_\ast^d} \int_\eta \abs{A^{(1)} D_\eta^\alpha \hat{f}_k(\eta)} e^{\frac{1}{2}\lambda(t)\jap{\ell,\ell t}^s}\abs{\hat{\rho}_\ell(t)_{N^\prime}} \\ & \quad\quad\quad \times A^{(1)}_{k-\ell}(t,\eta-\ell t)\abs{\left(D_\eta^\alpha \hat f_{k-\ell}(t,\eta - t\ell) \right)_{N}} d\eta,  
\end{align*}  
which implies by \eqref{ineq:L2L2L1}, \eqref{ineq:SobExp}, \eqref{ineq:LPOrthoProject} and $\sigma > d/2 + 2$, 
\begin{align} 
\abs{\mathcal{R}^{1;2}} & \lesssim \jap{t}e^{-\frac{1}{2}\alpha_0\jap{t}^s}\sum_{N  \in \mathbb D} \sum_{N^\prime \approx N}\norm{A^{(1)}(v^\alpha f)}_2 \frac{1}{N^\prime} \norm{A\hat{\rho}(t)_{N^\prime}}_2 \norm{A^{(1)}(v^\alpha f)_{N}}_2 \nonumber \\ 
& \lesssim \jap{t} e^{-\frac{1}{2}\alpha_0\jap{t}^s}\norm{A\rho(t)}_2 \sum_{N  \in \mathbb D}\norm{A^{(1)}(v^\alpha f)}_2 \frac{1}{N} \norm{A^{(1)}(v^\alpha f)_{N}}_2 \nonumber \\ 
& \lesssim \jap{t}^{-1}\norm{A\rho(t)}_2 e^{-\frac{1}{4}\alpha_0\jap{t}^s} \norm{A^{(1)}(v^\alpha f)}^2_2, \label{ineq:RRbHConc} 
\end{align} 
which suffices to treat $\mathcal{R}^{1;2}$. 

\subsubsection{Treatment of lower moments} \label{sec:HiLowMoments}
Next we turn to the treatment of $E_{NL}^2$. 
First apply \eqref{ineq:Wbd} (using $\gamma \geq 1$), 
\begin{align*}
\abs{E_{NL}^2} &  \lesssim \sum_{\abs{j} = 1; j \leq \alpha} \sum_{k \in \Integer^d} \sum_{\ell \in \Integer_\ast^d} \int_\eta \abs{A^{(1)}D_\eta^\alpha \hat{f}_k(\eta)}A^{(1)}_k(t,\eta) \abs{\jap{\ell}^{-1}\hat{\rho}_\ell(t)D_\eta^{\alpha-j} \hat f_{k-\ell}(t,\eta - t\ell)} d\eta. 
\end{align*} 
Then we apply \eqref{ineq:GProduct2}
\begin{align*} 
\abs{E_{NL}^2} & \lesssim  \sum_{\abs{j} = 1; j \leq \alpha} \norm{A^{(1)} v^\alpha f}_2 \norm{A^{(1)} v^{\alpha-j} f}_2 \norm{\rho(t)}_{\cF^{\tilde c\lambda(t),0;s}} \\  
& \quad +  \sum_{\abs{j} = 1; j \leq \alpha}\norm{A^{(1)} v^\alpha f}_2 \norm{v^{\alpha-j} f}_{\cG^{\tilde{c}\lambda,0;s}}\norm{\jap{k}^{-1} A_k^{(1)} \hat\rho_k(t)}_{L_k^2}. 
\end{align*}
From here, we take advantage of the regularity gap and $\jap{k}^{-1}A^{(1)}_k(k,k t) \lesssim \jap{t}A_k(t,k t)$ to deduce
\begin{align}
\abs{E_{NL}^2} & \lesssim  \sum_{\abs{j} = 1; j \leq \alpha} e^{(\tilde c - 1)\alpha_0\jap{t}^s}\norm{A^{(1)} v^\alpha f}_2 \norm{A^{(1)} v^{\alpha-j} f}_2 \norm{A\rho(t)}_2 \nonumber \\
& \quad +  \sum_{\abs{j} = 1; j \leq \alpha}\jap{t}\norm{A^{(1)} v^\alpha f}_2 \norm{A^{(-\beta)} v^{\alpha-j} f}_2 \norm{A\rho(t)}_{2}, \label{ineq:ENL2Conc}
\end{align}
which suffices to treat this contribution. 

\subsubsection{Conclusion of high norm estimate} 
Denote $\delta = -\frac{1}{4}\min\left(c-1,\tilde c - 1,c^\prime - 1\right)\alpha_0$.
Collecting the contributions of \eqref{ineq:HighNormDeriv}, \eqref{ineq:ELConc} 
\eqref{ineq:TNConc}, \eqref{ineq:RN11Conc}, \eqref{ineq:RN12Conc}, \eqref{ineq:RRaHConc}, \eqref{ineq:RRbHConc} and \eqref{ineq:ENL2Conc} 
then summing in $N$ with \eqref{ineq:GeneralOrtho} (note we used $1 \lesssim \jap{k,\eta}^{s/2}$ to group \eqref{ineq:RN12Conc}, \eqref{ineq:RRaHConc} and \eqref{ineq:RRbHConc} with \eqref{ineq:TNConc}), we have the following for some $\tilde{K} = \tilde{K}(s,M,\sigma,\lambda_0,\lambda^\prime,C_0,d)$, 
\begin{align*} 
\frac{1}{2}\frac{d}{dt} \norm{A^{(1)} v^\alpha f}_2^2  & \leq \left(\tilde K \jap{t}^{-1} e^{-\delta \jap{t}^{s}} \norm{A\rho(t)}_2 + \dot\lambda(t)\right)\norm{\jap{\grad_{z,v}}^{s/2} A^{(1)}(v^\alpha f)}_2^2 \\ 
& 
 \quad  + \tilde{K}\jap{t}\norm{A^{(1)} v^\alpha f}_2 \norm{A\rho(t)}_2 \\ 
& \quad + \tilde{K}\frac{\norm{A^{(-\beta)} v^\alpha f}_2}{\jap{t}^2} \norm{A^{(1)} v^\alpha f}^2_2  + \tilde{K}\jap{t}^6\norm{A^{(-\beta)} v^\alpha f}_2\norm{A\rho}^2_2 \\ 
& \quad + \tilde{K}\sum_{\abs{j} = 1; j \leq \alpha} e^{-\delta\jap{t}^s}\norm{A^{(1)} v^\alpha f}_2 \norm{A^{(1)} v^{\alpha-j} f}_2 \norm{A\rho(t)}_2 \\
& \quad + \tilde{K}\sum_{\abs{j} = 1; j \leq \alpha}\jap{t}\norm{A^{(1)} v^\alpha f}_2 \norm{A^{(-\beta)} v^{\alpha-j} f}\norm{A\rho(t)}_2. 
\end{align*} 
Introducing a small parameter $b$ to be fixed depending only on $\tilde K$ and $\lambda(t)$,
\begin{align} 
\frac{1}{2}\frac{d}{dt} \norm{A^{(1)} v^\alpha f}_2^2  & \leq \left(\tilde K e^{-\delta \jap{t}^{s}} \jap{t}^{-1}\norm{A\rho(t)}_2 + \frac{b\tilde K}{\jap{t}^2} +\dot\lambda(t)\right)\norm{\jap{\grad}^{s/2} A^{(1)}(v^\alpha f)}_2^2 \nonumber \\ 
& \quad + \frac{\tilde{K}}{b}\jap{t}^{4}\norm{A\rho(t)}^2_2  \nonumber \\ 
& \quad + \tilde{K}\frac{\norm{A^{(-\beta)} v^\alpha f}_2}{\jap{t}^2} \norm{A^{(1)} v^\alpha f}^2_2  + \tilde{K}\jap{t}^6\norm{A^{(-\beta)} v^\alpha f}_2\norm{A\rho}^2_2 \nonumber \\ 
& \quad + \tilde{K}\sum_{\abs{j} = 1; j \leq \alpha} e^{-\delta\jap{t}^s}\norm{A^{(1)} v^\alpha f}_2 \norm{A^{(1)} v^{\alpha-j} f}_2 \norm{A\rho(t)}_2 \nonumber \\
& \quad + \tilde{K}\sum_{\abs{j} = 1; j \leq \alpha}\jap{t}\norm{A^{(1)} v^\alpha f}_2 \norm{A^{(-\beta)} v^{\alpha-j} f}\norm{A\rho(t)}_2. \label{ineq:HiNormMidStep}
\end{align} 
By \eqref{ineq:SobExp} and \eqref{ineq:dotlambda} 
 we may fix $b$ and $\epsilon$ small such that  
\begin{align*}
\tilde K e^{-\delta \jap{t}^{s}}\sqrt{K_4}\epsilon + \frac{b\tilde K}{\jap{t}^2} \leq \frac{1}{2}\abs{\dot\lambda(t)}. 
\end{align*} 
Note this requires fixing $\epsilon$ small relative to $K_4$ but $b$ is chosen independently of $K_4$.
Then by \eqref{ctrl:MidPt} we deduce that the first term in \eqref{ineq:HiNormMidStep} is negative.  
Therefore, summing in $\alpha$, integrating and applying the bootstrap hypotheses \eqref{ctrl:Boot} and \eqref{ctrl:MidPt} implies (adjusting $\tilde K$ to $\tilde K^\prime$)
\begin{align*}
\sum_{\abs{\alpha} \leq M}\norm{A^{(1)} v^\alpha f}_2^2 & \leq \epsilon^2 + \tilde{K}^\prime \jap{t}^5 K_3 \epsilon^2 + \tilde K^\prime K_1 \sqrt{K_2} \jap{t}^6 \epsilon^3 + \tilde K^\prime \sqrt{K_2} \jap{t}^7 K_3 \epsilon^3 \\ & \quad + \tilde{K}^\prime K_1 \sqrt{K_4}\epsilon^3 + \tilde{K}^\prime \sqrt{K_1K_2 K_4} \jap{t}^6 \epsilon^3. 
\end{align*}
Hence we may take $K_1 = \tilde{K}^\prime K_3 + 1$ and we have \eqref{ctrl:HiLocalizedB} by choosing
\begin{align*}
\epsilon < K_1\left(4\tilde K^\prime\right)^{-1}\left(K_1\sqrt{K_2} + \sqrt{K_2}K_3 + K_1\sqrt{K_4} + \sqrt{K_1K_2K_4} \right)^{-1}. 
\end{align*}

\subsection{Proof of low norm estimate \eqref{ctrl:LowCommLocB}} \label{subsec:LoNorm}
This proof proceeds analogously to \eqref{ctrl:HiLocalizedB} replacing $A^{(1)}$ with $A^{(-\beta)}$.
First compute the derivative as in \eqref{ineq:HighNormDeriv}, 
\begin{align} 
\frac{1}{2}\frac{d}{dt}\norm{A^{(-\beta)}D_\eta^\alpha \hat f}_2^2 & = \sum_{k\in \Integer^d}\int_\eta \dot{\lambda}(t)\jap{k,\eta}^{s} \abs{A^{(-\beta)} D_\eta^\alpha \hat f_k(\eta)}^2 d\eta +   \sum_{k \in \Integer^d} \int_\eta A^{(-\beta)} D_\eta^\alpha \overline{\hat{f}_k(\eta)} A^{(-\beta)} D_\eta^\alpha \partial_t \hat{f}_k(\eta) d\eta \nonumber \\ 
& = CK_{\cL} + E_{\cL}, \label{ineq:LowNormDeriv}
\end{align}
where
\begin{align} 
E_{\cL} & = -\sum_{k \in \Integer^d} \int_\eta A^{(-\beta)} D_\eta^\alpha\overline{\hat{f}_k(\eta)} A^{(-\beta)}_k(t,\eta)\hat{\rho}_k(t) \widehat{W}(k) D_\eta^\alpha \left[ k \cdot (\eta - tk) \hat f^0(\eta - kt)\right] d\eta \nonumber \\ 
& \quad - \sum_{k\in\Integer^d} \sum_{\ell \in \Integer_\ast^d}\int_\eta A^{(-\beta)} D_\eta^\alpha\overline{\hat{f}_k(\eta)} A^{(-\beta)}_k(t,\eta) \hat{\rho}_\ell(t)\widehat{W}(\ell) D_\eta^\alpha \left[ \ell \cdot (\eta - tk) \hat f_{k-\ell}(t,\eta - t\ell)\right] d\eta \nonumber \\ 
& = -E_{\cL;L} - E_{\cL;NL}. \label{def:ELENL_Lo}
\end{align} 
As in the treatment of $E_L$ in \S\ref{subsec:linearContr}, we may use the product lemma \eqref{ineq:GProduct2} and \eqref{ineq:Wbd} to deduce 
\begin{align*}
\abs{E_{\cL;L}} & \lesssim \norm{A^{(-\beta)} D_\eta^\alpha \hat{f}}_2\norm{A^{(-\beta)} D_\eta^\alpha(\eta\hat{f}^0(\eta))}_{L^2_\eta} \norm{\rho(t)}_{\cF^{\tilde c \lambda(t),0;s}} 
  \\ & \quad + \norm{A^{(-\beta)} D_\eta^\alpha \hat{f}}_2\norm{A^{(-\beta)} D_\eta^\alpha(\eta\hat{f}^0(\eta))}_{L^2_\eta}\norm{A^{(-\beta)} \rho(t)}_2.
\end{align*}  
By the analogue of \eqref{ineq:BDalphaf0}, $\tilde c < 1$ and the regularity gap between $A^{(-\beta)}$ and $A$ and \eqref{ineq:SobExp} 
\begin{align} 
\abs{E_{\cL;L}} & \lesssim e^{(\tilde{c}-1)\alpha_0\jap{t}^s} \norm{A^{(-\beta)} D_\eta^\alpha \hat{f}}_2\norm{A\rho(t)}_2  + \jap{t}^{-\beta} \norm{A^{(-\beta)} D_\eta^\alpha \hat{f}}_2 \norm{A\rho(t)}_2  \nonumber \\ 
& \lesssim \jap{t}^{-\beta} \norm{A^{(-\beta)} v^\alpha f}_2 \norm{A\rho(t)}_2, \label{ineq:ELMConc}
\end{align} 
which suffices to treat this term.

We now turn to the treatment of $E_{\cL;NL}$, which as in \S\ref{sec:NL} is expanded by
\begin{align*}
E_{\cL;NL} & = \sum_{k \in \Integer^d} \int_\eta A^{(-\beta)}D_\eta^\alpha \overline{\hat{f}_k(\eta)}\left(A^{(-\beta)}_k(t,\eta) \left[\sum_{\ell \in \Integer_\ast^d} \hat{\rho}_\ell(t)\widehat{W}(\ell)\ell \cdot ( \eta - tk ) D_\eta^\alpha \hat f_{k-\ell}(t,\eta - t\ell)\right] \right) d\eta \\
& \quad +\sum_{k\in \Integer^d} \int_\eta A^{(-\beta)}D_\eta^\alpha \overline{\hat{f}_k(\eta)}\left(A^{(-\beta)}_k(t,\eta) \sum_{\abs{j} = 1; j \leq \alpha}\left[\sum_{\ell \in \Integer_\ast^d} \hat{\rho}_\ell(t)\widehat{W}(\ell)\ell_j D_\eta^{\alpha-j} \hat f_{k-\ell}(t,\eta - t\ell)\right] \right) d\eta \\ 
& = E_{\cL;NL}^1 + E_{\cL;NL}^2. 
\end{align*} 
First consider $E_{\cL;NL}^1$, to which we apply \eqref{eq:cancel} (with $A^{(-\beta)}$ instead of $A^{(1)}$) and then decompose via paraproduct as in \eqref{def:Edecomp}:   
\begin{align} 
E_{\cL;NL}^1 & = \sum_{N \geq 8} T^1_{\cL;N} + \sum_{N \geq 8} R^1_{\cL;N} + \mathcal{R}_{\cL}^1, \label{def:EMdecomp}
\end{align} 
where the \emph{transport} term is given by  
\begin{align} 
T^1_{\cL;N} & = \sum_{k\in\Integer^d} \sum_{\ell \in \Integer_\ast^d} \int_\eta A^{(-\beta)} D_\eta^\alpha \overline{\hat{f}_k(\eta)} \hat{\rho}_\ell(t)_{<N/8}\widehat{W}(\ell)\ell \cdot  (\eta - kt) \nonumber \\ & \quad\quad \times \left[A^{(-\beta)}_{k}(t,\eta) - A^{(-\beta)}_{k-\ell}(t,\eta-t\ell)\right]\left(D_\eta^\alpha \hat f_{k-\ell}(t,\eta - t\ell) \right)_{N} d\eta, \label{def:TMn}
\end{align} 
and the \emph{reaction} term by 
\begin{align} 
R^1_{\cL;N} & = \sum_{k\in\Integer^d} \sum_{\ell \in \Integer_\ast^d} \int_\eta A^{(-\beta)} D_\eta^\alpha \overline{\hat{f}_k(\eta)} \hat{\rho}_\ell(t)_{N}\widehat{W}(\ell)\ell \cdot \left(\eta - kt\right) \nonumber \\ & \quad\quad \times \left[A^{(-\beta)}_{k}(t,\eta) - A^{(-\beta)}_{k-\ell}(t,\eta-t\ell)\right] \left(D_\eta^\alpha \hat f_{k-\ell}(t,\eta - t\ell) \right)_{<N/8}d\eta, \label{def:RMn}
\end{align}
and as before the remainder is whatever is left over. 
The treatment of the transport term $T_{\cL;N}^1$ and the remainder $\mathcal{R}_{\cL}^1$ is unchanged from the corresponding treatments of $T_{N}^1$ and $\mathcal{R}^1$ in \S\ref{sec:HiNormTrans} and \S\ref{sec:HiNormRemainder} respectively. Hence, we omit it and simply conclude as in \eqref{ineq:TNConc}, \eqref{ineq:RRaHConc} and \eqref{ineq:RRbHConc}: 
\begin{align} 
\abs{T_{\cL;N}^1} & \lesssim \jap{t}^{-1}e^{\frac{1}{2}(c-1)\alpha_0 \jap{t}^{s}} \norm{A\rho(t)}_2 \norm{\jap{\grad_{z,v}}^{s/2}A^{(-\beta)} (v^\alpha f)_{\sim N}}_2^2, \label{ineq:TNcMConc} \\ 
\abs{\mathcal{R}_{\cL}^1} & \lesssim \norm{A\rho(t)}_2\jap{t}^{-1}e^{\frac{1}{4}(c^\prime - 1)\alpha_0\jap{t}^s} \norm{A^{(-\beta)} (v^\alpha f)}^2_2. \label{ineq:RRcMConc}
\end{align}  
The reaction term is slightly altered to gain from the regularity gap and get a uniform bound (as in the linear contribution $E_{\cL;L}$). 
As in the treatment of reaction in $E_{NL}^1$ in \S\ref{sec:HiReac}, we separate into $R_{\cL;N}^{1} =  R_{\cL;N}^{1;1} + R_{\cL;N}^{1;2}$ 
where the leading order reaction term is given by
\begin{align*} 
R_{\cL;N}^{1;1} =  -\sum_{k\in \Integer^d} \sum_{\ell \in \Integer_\ast^d} \int_\eta A^{(-\beta)} D_\eta^\alpha \overline{\hat{f}_k(\eta)} A^{(-\beta)}_k(t,\eta)\hat{\rho}_\ell(t)_{N}\widehat{W}(\ell) \ell\cdot \left[ \eta - tk\right] \left(D_\eta^\alpha \hat f_{k-\ell}(t,\eta - t\ell)\right)_{<N/8}  d\eta. 
\end{align*}  
By the frequency localizations \eqref{ineq:FreqLocExample}, 
 \eqref{lem:scon2} implies for some $c = c(s) \in (0,1)$ (using also \eqref{ineq:Wbd} and $\abs{\eta - kt} \leq \jap{t}\jap{k-\ell,\eta-t\ell}$), 
\begin{align*} 
\abs{R_{\cL;N}^{1;1}} 
&\lesssim \jap{t}\sum_{k \in \Integer_\ast^d} \sum_{\ell \in \Integer_\ast^d} \int_\eta \abs{ A^{(-\beta)} \left(D_\eta^\alpha \hat{f}_k(\eta)\right)_{\sim N}}A^{(-\beta)}_\ell(t,t\ell) \abs{\hat{\rho}_\ell(t)_{N}} \\ & \quad\quad \times e^{c\lambda(t)\jap{k-\ell,\eta-\ell t}^s} \jap{k-\ell,\eta - \ell t} \abs{ \left(D_\eta^\alpha \hat f_{k-\ell}(t,\eta - t\ell)\right)_{<N/8}} d\eta. 
\end{align*} 
Proceeding as in the proof of \eqref{ineq:RN11Conc}, applying \eqref{ineq:L2L1L2} (along with $\sigma>\frac{d}{2} + 2$) and using the regularity gap between $A^{(-\beta)}$ and $A$ implies
\begin{align} 
\abs{R_{\cL;N}^{1;1}} & \lesssim \jap{t}\norm{A^{(-\beta)} \left(v^\alpha f\right)_{\sim N}}_2 \norm{A^{(-\beta)}\rho_N}_2 \norm{A^{(-\beta)} v^\alpha f}_2 \nonumber \\ 
& \lesssim \jap{t}^{1-\beta} \norm{A^{(-\beta)} \left(v^\alpha f\right)_{\sim N}}_2 \norm{A\rho_N}_2\norm{A^{(-\beta)} v^\alpha f}_2 \nonumber \\  
& \lesssim \jap{t}^{1-\beta} \norm{A^{(-\beta)} v^\alpha f}_2 \norm{A^{(-\beta)}\left(v^\alpha f\right)_{\sim N}}^2_2 + \jap{t}^{1-\beta}\norm{A^{(-\beta)} v^\alpha f}_2\norm{A\rho_N}^2_2, \label{ineq:Rcm11Conc}
\end{align}
which will be sufficient for the proof of \eqref{ctrl:LowCommLocB}. 
The term $R_{\cL;N}^{1;2}$ can be treated exactly as $R_{N}^{1;2}$ and hence we omit and simply conclude
\begin{align} 
\abs{R_{\cL;N}^{1;2}} & \lesssim \frac{e^{-\alpha_0\jap{t}^s}}{N} \norm{A^{(-\beta)}(v^\alpha f)}^2_2 \norm{A\rho(t)}_2. \label{ineq:RN12cMConc}
\end{align}

The term $E_{\cL;NL}^2$ is treated as in \S\ref{sec:HiLowMoments}. By \eqref{ineq:Wbd}, \eqref{ineq:GProduct2} and the regularity gap between $A^{(-\beta)}$ and $A$ (also \eqref{ineq:SobExp} in the last line), 
\begin{align}
\abs{E_{\cL;NL}^2} & \lesssim  \sum_{\abs{j} = 1; j \leq \alpha} \norm{A^{(-\beta)} v^\alpha f}_2 \norm{A^{(-\beta)} v^{\alpha-j} f}_2 \norm{\rho(t)}_{\cF^{\tilde c \lambda(t),0;s}} \nonumber \\ 
& \quad +  \sum_{\abs{j} = 1; j \leq \alpha}\norm{A^{(-\beta)} v^\alpha f}_2 \norm{v^{\alpha-j} f}_{\cG^{\tilde{c}\lambda,\sigma;s}} \norm{A^{(-\beta)} \rho(t)}_2 \nonumber \\ 
& \lesssim  e^{(\tilde{c} - 1)\alpha_0\jap{t}^s}\sum_{\abs{j} = 1; j \leq \alpha} \norm{A^{(-\beta)} v^\alpha f}_2 \norm{A^{(-\beta)} v^{\alpha-j} f}_2 \norm{A\rho(t)}_2 \nonumber \\
& \quad +  \jap{t}^{-\beta} \sum_{\abs{j} = 1; j \leq \alpha}\norm{A^{(-\beta)} v^\alpha f}_2 \norm{A^{(-\beta)} v^{\alpha-j} f}_2\norm{A\rho(t)}_2 \nonumber \\ 
& \lesssim \jap{t}^{-\beta} \sum_{\abs{j} = 1; j \leq \alpha}\norm{A^{(-\beta)} v^\alpha f}_2 \norm{A^{(-\beta)} v^{\alpha-j} f}_2\norm{A\rho(t)}_2. \label{ineq:ENLcM2Conc}
\end{align} 
Denote $\delta = -\frac{1}{4}\min\left(c-1,\tilde c - 1,c^\prime - 1\right)\alpha_0$.
Collecting \eqref{ineq:ELMConc} \eqref{ineq:TNcMConc}, \eqref{ineq:RRcMConc}, \eqref{ineq:Rcm11Conc}, \eqref{ineq:RN12cMConc} and \eqref{ineq:ENLcM2Conc} and summing in $N$, splitting the linear terms with a small parameter $b$ and combining \eqref{ineq:RRcMConc} and \eqref{ineq:RN12cMConc} with \eqref{ineq:TNcMConc} as in \eqref{ineq:HiNormMidStep} (using also \eqref{ineq:SobExp}),
we have the following for some $\tilde K = \tilde{K}(s,\sigma,\alpha_0,C_0,d)$ (not the same as the $\tilde K$ in \eqref{ineq:HiNormMidStep} but this is irrelevant), 
\begin{align} 
\frac{1}{2}\frac{d}{dt} \norm{A^{(-\beta)} v^\alpha f}_2^2  & \leq \left(\tilde K \jap{t}^{-1}e^{-\delta \jap{t}^{s}} \norm{A\rho(t)}_2 + \tilde{K}\jap{t}^{-\beta}b + \dot\lambda(t)\right)\norm{\jap{\grad_{z,v}}^{s/2} A^{(-\beta)} (v^\alpha f)}_2^2 \nonumber \\ 
& \quad  + \frac{\tilde K}{b}\jap{t}^{-\beta}\norm{A\rho(t)}_2^2 \nonumber \\ 
& \quad + \tilde{K} \jap{t}^{1-\beta} \norm{A^{(-\beta)} (v^\alpha f)}_2\left(\norm{A^{(-\beta)}(v^\alpha f)}_2^2 + \norm{A\rho(t)}_2^2 \right) \nonumber \\ 
& \quad +  \tilde K \jap{t}^{-\beta} \sum_{\abs{j} = 1; j \leq \alpha}\norm{A^{(-\beta)} v^\alpha f}_2 \norm{A^{(-\beta)} v^{\alpha-j} f}\norm{A\rho(t)}_2.\label{ineq:LoNormMidStep}
\end{align} 
By \eqref{ineq:SobExp} and \eqref{ineq:dotlambda}  we may fix $b$ and $\epsilon$ small such that
\begin{align*}
\tilde K \sqrt{K_4} \epsilon e^{-\delta \jap{t}^{s}} + \tilde{K}\jap{t}^{-\beta}b \leq \frac{1}{2}\abs{\dot\lambda(t)},
\end{align*} 
which by \eqref{ctrl:MidPt}, implies that the first term in \eqref{ineq:LoNormMidStep} is non-positive.
Therefore, summing in $\alpha$, integrating with $\beta > 2$ and applying the bootstrap hypotheses \eqref{ctrl:Boot} and \eqref{ctrl:MidPt} implies (adjusting $\tilde K$ to $\tilde K^\prime$), 
\begin{align*} 
\sum_{\abs{\alpha} \leq M} \norm{A^{(-\beta)} v^\alpha f}_2^2 \leq \epsilon^2 + \tilde{K}^\prime K_3\epsilon^2 + \tilde{K}^\prime K_2^{3/2}\epsilon^3 + \tilde{K}^\prime \sqrt{K_2}K_3 \epsilon^3 + \tilde{K}^\prime K_2\sqrt{K_4} \epsilon^3. 
\end{align*} 
Hence, we take $K_2 = 1 + \tilde{K}^\prime K_3$ and $\epsilon < K_2 \left(3\tilde{K}^\prime\right)^{-1}\left(K_2^{3/2} + \sqrt{K_2}K_3 + K_2\sqrt{K_4}\right)^{-1}$ to deduce \eqref{ctrl:LowCommLocB}. 

\section{Analysis of the plasma echoes} \label{sec:Proofa}
The most important step to pushing linear Landau damping to the nonlinear level is analyzing and controlling the dominant weakly nonlinear effect: \emph{the plasma echo}.
Mathematically, this comes down to verifying condition \eqref{ineq:momentdef} on the time-response kernels, crucial to the proof of \eqref{ctrl:Mid} in \S\ref{sec:L2I}. 
Our choices of $\lambda(t)$ for $t \gg 1$ (in particular the choice of $a$) and $s > 1/(2+\gamma)$ are both determined in this section.
The analysis in this section is similar to the moment estimates carried out on the time-response kernels in \S7 of \cite{MouhotVillani11} except with the regularity loss encoded by our choice of $\lambda(t)$ taking the place of amplitude growth. 
The distinction is arguably minor, but this increased precision allows for a slightly cleaner treatment and highlights more clearly the origin of the regularity requirement.  

\begin{lemma}[Time response estimate I] \label{lem:Moment}
Under the bootstrap hypotheses \eqref{ctrl:Boot}, there holds
\begin{align*} 
\sup_{t \in [0,T^\star]} \sup_{k\in \Integer_\ast^d}\int_0^t\sum_{\ell \in \Integer_\ast^d} \bar{K}_{k,\ell}(t,\tau) d\tau \lesssim_{a,s,d,\lambda_0,\lambda^\prime} \sqrt{K_2} \epsilon.  
\end{align*} 
\end{lemma} 
\begin{proof}
Consider first the effect of $f_0$, the homogeneous part of $f$, which corresponds to $\bar{K}_{k,k}(t,\tau)$: 
\begin{align*} 
\mathcal{I}_{inst}(t) & :=  \int_0^t e^{(\lambda(t) - \lambda(\tau))\jap{k,kt}^s} e^{c\lambda(\tau)\jap{k(t-\tau)}^s} \frac{\abs{k(t-\tau)}}{\abs{k}^{\gamma}}\abs{\widehat{f}_{0}(\tau,k(t-\tau))} d\tau \\ 
& \leq \int_0^te^{c\lambda(\tau)\jap{k(t-\tau)}^s} \frac{\abs{k(t-\tau)}}{\abs{k}^{\gamma}}\abs{\widehat{f}_{0}(\tau,k(t-\tau))} d\tau. 
\end{align*} 
Here $inst$ stands for `instantaneous' as this effect has no time delay (unlike $k \neq \ell$ below); this terminology was borrowed from \cite{MouhotVillani11}. 
Also note that this is only controlling the effect of `low' frequencies in $f_0$.
From the $H^{d/2+} \hookrightarrow C^0$ embedding, $\sigma > \beta + 1$ and \eqref{ineq:OuterMctrl},
\begin{align*} 
\mathcal{I}_{inst}(t) & \leq \int_0^t e^{(c-1)\lambda(\tau)\jap{k(t-\tau)}^s} \left(\sup_{\eta \in \Real^d} e^{\lambda(\tau)\jap{\eta}^s} \abs{\eta}\abs{\widehat{f}_{0}(\tau,\eta)} \right)   d\tau \\ 
& \lesssim_M \int_0^t e^{(c-1)\lambda(\tau)\jap{k(t-\tau)}^s} \norm{A^{(-\beta)}f_0(\tau)}_{H^{M}_\eta} d\tau \\ 
& \lesssim_{\alpha_0} \sqrt{K_2}\epsilon. 
\end{align*} 

Next turn to the contributions from the case $k \neq \ell$, which is the origin of the plasma echoes.   
Using $\abs{k(t-\tau)} \leq \jap{\tau}\abs{k-\ell,kt - \ell\tau}$ and the definition of $\bar K$ in \eqref{def:barK}, 
\begin{align*} 
\mathbf{1}_{k \neq \ell} \bar{K}_{k,\ell}(t,\tau) & \lesssim  e^{(\lambda(t) - \lambda(\tau))\jap{k,kt}^s}e^{c\lambda(\tau)\jap{k-\ell,kt-\ell\tau}^s}\frac{\jap{\tau}}{\abs{\ell}^{\gamma}}\abs{\widehat{\grad f}_{k-\ell}(\tau,kt-\ell\tau)} \mathbf{1}_{k \neq \ell \neq 0}.
\end{align*}
In what follows denote
\begin{align*} 
-\nu(t,\tau) = \lambda(t) - \lambda(\tau). 
\end{align*} 
Then using that $\lambda(t) \geq \alpha_0$ and $c < 1$, if we write $\delta = (1 - c)\alpha_0$ we are left to estimate, 
\begin{align} 
\mathcal{I}(t) & := \int_0^t\sum_{\ell \in \Integer_\ast^d} e^{-\nu(t,\tau)\jap{k,kt}^s}e^{c\lambda(\tau)\jap{k-\ell,kt-\ell\tau}^s}\frac{\jap{\tau}}{\abs{\ell}^\gamma}\abs{\widehat{\grad f}_{k-\ell}(\tau,kt-\ell\tau)} \mathbf{1}_{k \neq \ell} d\tau \nonumber \\ 
& \lesssim \int_0^t\sum_{\ell \in \Integer_\ast^d} e^{-\nu(t,\tau)\jap{k,kt}^s} \frac{\jap{\tau}}{\abs{\ell}^\gamma} e^{-\delta \jap{k-\ell,kt-\ell\tau}^s}\abs{e^{\lambda(\tau)\jap{k-\ell,kt-\ell\tau}^s} \widehat{\grad f}_{k-\ell}(\tau,kt-\ell\tau)} \mathbf{1}_{k\neq \ell} d\tau. \label{ineq:mathcalI}
\end{align} 
By $\sigma \geq \beta+1$, the $H^{d/2+} \hookrightarrow C^0$ embedding and \eqref{ineq:OuterMctrl},
\begin{align*} 
\abs{e^{\lambda(\tau)\jap{k-\ell,kt-\ell\tau}^s} \widehat{\grad f}_{k-\ell}(\tau,kt-\ell\tau)} & \leq  \sup_{\eta \in \Real^d} e^{\lambda(\tau)\jap{k-\ell,\eta}^s} \jap{k-\ell,\eta} \abs{\widehat{f}_{k-\ell}(\tau,\eta)} \\ 
& \leq \left(\sum_{k  \in \Integer^d} \sup_{\eta \in \Real^d} e^{2\lambda(\tau)\jap{k,\eta}^s} \jap{k,\eta}^2 \abs{\widehat{f}_{k}(\tau,\eta)}^2\right)^{1/2} \\ 
&\lesssim \norm{A^{(-\beta)}f(\tau)}_{L_k^2 H_\eta^M} \\ 
&\lesssim_M \sqrt{K_2}\epsilon.
\end{align*} 
Applying this to \eqref{ineq:mathcalI} implies 
\begin{align} 
\mathcal{I}(t) \lesssim \sqrt{K_2}\epsilon\int_0^t\sum_{\ell \in \Integer_\ast^d} e^{-\nu(t,\tau)\jap{k,kt}^s} \frac{\jap{\tau}}{\abs{\ell}^\gamma} e^{-\delta \jap{k-\ell,kt-\ell\tau}^s} \mathbf{1}_{\ell \neq k} d\tau. \label{ineq:mIReduced}
\end{align} 
Following an argument similar to that in \cite{MouhotVillani11} we may reduce to the $d = 1$ case. 
By \eqref{ineq:Comps}, 
\begin{align*} 
\mathcal{I}(t) & \lesssim \sqrt{K_2}\epsilon\int_0^t\sum_{\ell \in \Integer_\ast^d} \sum_{j: \ell_j \neq k_j} e^{-\nu(t,\tau)\jap{k_j,k_jt}^s} \frac{\jap{\tau}}{\abs{\ell}^\gamma} e^{-C_s\delta \jap{k_j-\ell_j,k_jt-\ell_j\tau}^s} \prod_{i \neq j}^d e^{-C^{d-1}_s\delta\jap{k_i - \ell_i}^s} \mathbf{1}_{\ell \neq k} d\tau \\ 
& \lesssim \frac{\sqrt{K_2}\epsilon}{\delta^{\frac{d-1}{s}}} \sum_{1 \leq j \leq d}\int_0^t\sum_{\ell_j \in \Integer} e^{-\nu(t,\tau)\jap{k_j,k_jt}^s} \frac{\jap{\tau}}{\jap{\ell_j}^\gamma} e^{-C_s\delta \jap{k_j-\ell_j,k_jt-\ell_j\tau}^s} \mathbf{1}_{\ell_j \neq k_j} d\tau. 
\end{align*} 
Notice that we may not assert that both $k_j$ and $\ell_j$ are non-zero. 
However, if either $k_j$ or $\ell_j$ is zero we have by \eqref{ineq:Comps}, \eqref{ineq:SobExp} and $\tau \leq t$, 
\begin{align*} 
\jap{\tau}e^{-C_s\delta \jap{k_j-\ell_j,k_jt-\ell_j\tau}^s} \leq \jap{\tau}e^{-C_s^2\delta \jap{k_j-\ell_j}^s - C_s^2 \delta\jap{k_jt-\ell_j\tau}^s}
 \lesssim \delta^{-1/s}e^{-C_s^2\delta \jap{k_j-\ell_j}^s - \frac{1}{2}C_s^2 \delta\jap{\tau}^s}. 
\end{align*}
Hence, we see that such cases cannot contribute anything to the sum in $k$ of $\mathcal{I}(t)$ which is not bounded uniformly in time. 
Therefore, up to adjusting the definition of $\delta$ by a constant, we may focus on the cases such that both $k,\ell\in \Integer_\ast$ and $k \neq \ell$.  
Let us now focus on one such choice: 
\begin{align*} 
\mathcal{I}_{k,\ell}(t) := \int_0^t e^{-\nu(t,\tau)\jap{k,kt}^s} \frac{\jap{\tau}}{\abs{\ell}^\gamma} e^{-\delta \jap{k-\ell,kt-\ell\tau}^s} d\tau. 
\end{align*} 
This term isolates a single possible echo at $\tau = t k/\ell$: notice how the integrand is sharply localized near this time which accounts for the effect $\rho_\ell(\tau \ell)$ has on the behavior of $\rho_k(kt)$. 
Summing them deals with the cumulative effect of all the echoes. See \cite{MouhotVillani11} for more discussion. 
By symmetry we need only consider the case $k \geq 1$.  

Let us first eliminate the irrelevant early times; indeed by \eqref{ineq:Comps}, 
\begin{align} 
\int_0^{\min(1,t)} e^{-\nu(t,\tau)\jap{k,kt}^s} \frac{\jap{\tau}}{\abs{\ell}^\gamma} e^{-\delta \jap{k-\ell,kt-\ell\tau}^s} d\tau & \lesssim \int_0^{\min(1,t)} \frac{1}{\abs{\ell}^\gamma} e^{-\delta \jap{k-\ell,kt-\ell\tau}^s} d\tau  
\lesssim \frac{e^{-C_s\delta\jap{k-\ell}^s}}{\delta^{1/s}\abs{\ell}^{1+\gamma}}. \label{ineq:shorttime}
\end{align} 
Now let us turn to the more interesting contributions of $t \geq \tau \geq 1$. 
Given $t$,$k$ and $\ell$, define the \emph{resonant interval} as 
\begin{align*} 
I_R = \set{\tau \in [1,t] : \abs{kt - \ell\tau} < \frac{t}{2} }  
\end{align*} 
and divide $\mathcal{I}_{k,\ell}$ into three contributions (one from \eqref{ineq:shorttime}):  
\begin{align*} 
\mathcal{I}_{k,\ell}(t) & \lesssim \frac{1}{\delta^{1/s}\abs{\ell}^{1+\gamma}} e^{-C_s\jap{k-\ell}^s} + \left(\int_{[1,t] \cap I_R} + \int_{[1,t] \setminus I_R}\right)\frac{\jap{\tau}}{\abs{\ell}^\gamma}e^{-\delta\jap{k-\ell,kt-\ell\tau}^s} e^{-\nu(t,\tau)\jap{k,kt}^s} d\tau \\ 
 & = \frac{1}{\delta^{1/s}\abs{\ell}^{1+\gamma}} e^{-C_s\delta\jap{k-\ell}^s} + \mathcal I_{R} + \mathcal I_{NR}. 
\end{align*} 
Here `NR' stands for `non-resonant'.
Note that if $\ell \leq k-1$ then in fact $[0,t] \cap I_R = \emptyset$. 

Consider first the easier case of $\mathcal{I}_{NR}$. 
Since $\abs{kt - \ell t} \geq t/2$ on the support of the integrand, by \eqref{ineq:Comps} and \eqref{ineq:SobExp} we have
\begin{align} 
\mathcal{I}_{NR} & \leq  \frac{\jap{t}}{\abs{\ell}^\gamma}\int_{[1,t] \setminus I_R} e^{-C_s\delta\jap{k-\ell}^s-C_s\delta\jap{kt - \ell\tau}^s} e^{- \nu(\tau,t) \jap{k,kt}^s} \dd \tau \nonumber \\
 & \leq \frac{\jap{t}}{\abs{\ell}^\gamma} e^{-C_s\delta\jap{k-\ell}^s - \frac{1}{2}C_s\delta\jap{\frac{t}{2}}^s}\int_0 ^te^{- \frac{1}{2}C_s\delta\jap{kt - \ell\tau}^s} \dd \tau \nonumber \\ 
& \lesssim \frac{\jap{t}}{\delta^{1/s}\abs{\ell}^{1+\gamma}} e^{-C_s\delta\jap{k-\ell}^s - \frac{1}{2}C_s\delta\jap{\frac{t}{2}}^s} \nonumber \\ 
& \lesssim \frac{1}{\delta^{2/s}\abs{\ell}^{1+\gamma}} e^{-C_s\delta\jap{k-\ell}^s}, \label{ineq:INR}
\end{align} 
which suffices to treat all of the non-resonant contributions. 

Now focus on the resonant contribution $\mathcal{I}_{R}$, which as pointed out above, is only present if $\ell \geq k+1$ due to the echo at $\tau = tk/\ell \in (0,t)$. 
Since we are interested in  $t \geq \tau \geq 1$, by the definition of $\lambda(t)$ in \eqref{def:lambda}, there exists some constant $\delta^\prime$ (possibly adjusted by the reduction to one dimension) which is proportional to $\lambda_0 - \lambda^\prime$ such that on $[1,t]\cap I_R$, 
\begin{align*} 
\nu(t,\tau) = \delta^\prime\left(\tau^{-a} - t^{-a} \right) = \delta^\prime \left(\frac{t^a - \tau^a}{\tau^at^a} \right). 
\end{align*}
For $t$ and $\tau$ well separated, this provides a gap of regularity that helps us to control $\mathcal I_R$. 
Hence, we see that the most dangerous echoes occur for $\ell \approx k$ as these echoes are stacking up near $t$ and the regularity gap provided by $\nu$ becomes very small. 
From the formal analysis of \cite{MouhotVillani11} we expect to find the requirement $s > 1/(2+\gamma)$ due precisely to this effect.
Indeed we will see that is the case, in fact, here is the only place in the proof of Theorem \ref{thm:Main} where this requirement is used (also at the analogous step in the proof of Lemma \ref{lem:dual} below). 
By the mean-value theorem and the restriction that $\tau \in I_R$ (also $\tau \leq \frac{3kt}{2\ell}$ and $\ell - k \geq 1$), we have 
\begin{align} 
\nu(t,\tau) \geq a\delta^\prime \frac{t - \tau}{\tau^a t} & = \frac{a\delta^\prime}{\tau^a t}\left[t - \frac{kt}{\ell}\right] -\frac{a\delta^\prime}{\tau^a t}\left[\tau - \frac{kt}{\ell}\right] \nonumber \\ 
& \geq \frac{a\delta^\prime}{\tau^a}\left[1 - \frac{k}{\ell}\right] - \frac{a\delta^\prime}{2 \tau^a \ell} \nonumber \\ 
& \geq \frac{a\delta^\prime}{2\tau^a \ell} \nonumber \\ 
& \geq \frac{a\delta^\prime}{2^{1-a} 3^a (kt)^a \ell^{1-a}}. \label{ineq:LowerBdHardnu}
\end{align} 
Let $\tilde\delta^{\prime} =\frac{a\delta^\prime}{2^{1-a}3^a}$. 
The lower bound \eqref{ineq:LowerBdHardnu} precisely measures the usefulness of $\nu$. 
Indeed, by \eqref{ineq:LowerBdHardnu}, \eqref{ineq:Comps}, \eqref{ineq:SobExp} and $(2+\gamma)(s-a) = 1-a$ we have  
\begin{align} 
\mathcal{I}_R & \lesssim \int_{I_R} \frac{kt}{\ell^{1+\gamma}}e^{-\delta\jap{k-\ell,kt-\ell\tau}^s} e^{-\frac{\tilde\delta^\prime}{\ell^{1-a}}\abs{kt}^{s-a}} d\tau \nonumber \\ 
& \lesssim \frac{kt}{\delta^{1/s} \ell^{2+\gamma}}e^{-\frac{\tilde \delta^\prime}{\ell^{1-a}}\abs{kt}^{s-a}} e^{-C_s \delta \jap{k-\ell}^s} \nonumber \\ 
& \lesssim \frac{kt}{\delta^{1/s} \ell^{2+\gamma}}\left(\frac{\ell^{\frac{1-a}{s-a}}}{(\tilde\delta^{\prime})^{\frac{1}{s-a}}kt} \right) e^{-C_s \delta \jap{k-\ell}^s} \nonumber \\
& \lesssim_{s,a} e^{-C_s \delta \jap{k-\ell}^s} \frac{1}{\delta^{1/s}(a\delta^\prime)^{\frac{1}{s-a}}}. \label{ineq:IRhard}
\end{align} 
The use of $(2+\gamma)(s-a) \geq 1-a$ above is exactly the mathematical origin of the requirement $s > (2+\gamma)^{-1}$. 
Notice also that \eqref{ineq:IRhard} can be summed in either $k$ or $l$, but not in both. 

Assembling \eqref{ineq:shorttime}, \eqref{ineq:INR} and \eqref{ineq:IRhard} implies the lemma after summing in $\ell$ and taking the supremum in $t$ and $k$. 
\end{proof} 

The next estimate is in some sense the `dual' of Lemma \ref{lem:Moment} and is proved in the same way.
\begin{lemma}[Time response estimate II] \label{lem:dual}
Under the bootstrap hypotheses \eqref{ctrl:Boot} there holds
\begin{align*}
\sup_{\tau \in [0,T^\star]} \sup_{\ell \in \Integer_\ast^d} \sum_{k \in \Integer_\ast^d} \int_{\tau}^{T^\star}\bar{K}_{k,\ell}(t,\tau) dt \lesssim_{a,s,d,\lambda_0,\lambda^\prime} \sqrt{K_2} \epsilon. 
\end{align*}
\end{lemma} 
\begin{proof} 
First consider $\bar{K}_{k,k}(t,\tau)$, which corresponds to the homogeneous part of $f$:
\begin{align*} 
\mathcal{I}_{inst}(\tau) & :=  \int_\tau^{T^\star} e^{(\lambda(t) - \lambda(\tau))\jap{k,kt}^s} e^{c\lambda(\tau)\jap{k(t-\tau)}^s} \frac{\abs{k(t-\tau)}}{\abs{k}^{\gamma}}\abs{\widehat{f}_{0}(\tau,k(t-\tau))} dt.
\end{align*} 
By the same argument as used in Lemma \ref{lem:Moment}, it is straightforward to show 
\begin{align*} 
\mathcal{I}_{inst}(\tau) & \lesssim \sqrt{K_2}\epsilon. 
\end{align*}

Next consider the case $k \neq \ell$. 
By following the analysis of Lemma \ref{lem:Moment} the problem reduces to analyzing the analogue of \eqref{ineq:mIReduced}: 
\begin{align*} 
\mathcal{I}(\tau) = \sqrt{K_2}\epsilon\int_\tau^{T^\star} \sum_{k \in \Integer_\ast^d} e^{-\nu(t,\tau)\jap{k,kt}^s}e^{-\delta\jap{k-\ell,kt-\ell\tau}^s}\frac{\jap{\tau}}{\abs{\ell}^\gamma} \mathbf{1}_{\ell \neq k} dt 
\end{align*}
where $\nu(t,\tau) = \lambda(\tau) - \lambda(t)$ and $\delta$ are defined as in Lemma \ref{lem:Moment}. 
As before we may reduce to the one dimensional case at the cost of adjusting the constant and the definition of $\delta$. 
Hence consider the one dimensional integrals with $k,\ell \in \Integer_\ast$, $k \neq \ell$ and $k \geq 1$ (by symmetry):  
\begin{align} 
\mathcal{I}_{k,\ell}(\tau) = \int_\tau^{T^\star} e^{-\nu(t,\tau)\jap{k,kt}^s}e^{-\delta \jap{k-\ell,kt-\ell\tau}^s}\frac{\jap{\tau}}{\abs{\ell}^\gamma} dt. \label{ineq:IklDual}
\end{align} 

As in the proof of Lemma \ref{lem:Moment}, we may eliminate early times; we omit the details and conclude 
\begin{align*} 
\int_\tau^{\max(\tau,\min(1,T^\star))} e^{-\nu(t,\tau)\jap{k,kt}^s}e^{-\delta \jap{k-\ell,kt-\ell\tau}^s}\frac{\jap{\tau}}{\abs{\ell}^\gamma} dt \lesssim \frac{1}{\delta^{1/s}\abs{\ell}^{\gamma}\abs{k}} e^{-C_s\delta\jap{k-\ell}^s}. 
\end{align*} 
For the remainder of the proof, we will henceforth just assume $T^\star > \tau \geq 1$. 
Following the proof of Lemma \ref{lem:Moment}, define the resonant interval as
\begin{align*} 
I_R = \set{t \in [\tau,T^\star] : \abs{kt - \ell\tau} < \frac{\tau}{2} }  
\end{align*} 
and divide the integral into two main contributions: 
\begin{align*} 
\mathcal{I}_{k,\ell}(\tau) & =  \left(\int_{[\tau,T^\star) \cap I_R} + \int_{ [\tau,T^\star) \setminus I_R}\right)\frac{\jap{\tau}}{\abs{\ell}^\gamma}e^{-\delta\jap{k-\ell,kt-\ell\tau}^s} e^{-\nu(t,\tau)\jap{k,kt}^s} d t \\  & = \mathcal I_{R} + \mathcal I_{NR}. 
\end{align*}
The non-resonant contribution follows essentially the same proof as \eqref{ineq:INR} in Lemma \ref{lem:Moment}; we omit the details and conclude
\begin{align} 
\mathcal I_{NR} \lesssim \frac{1}{\delta^{2/s}\abs{\ell}^{\gamma} k} e^{-C_s\delta\jap{k-\ell}^s}. \label{ineq:INRDual}
\end{align} 
Turn now to the resonant integral, 
in which case $\ell \geq k+1$, and there is an echo at $ \ell \tau/k = t \in (\tau,\infty)$.
Since we are interested in  $t \geq \tau \geq 1$, by the definition of $\lambda(t)$ in \eqref{def:lambda}, there exists some constant $\delta^\prime$ (possibly adjusted by the reduction to one dimension) which is proportional to $\lambda_0 - \lambda^\prime$ such that by the mean-value theorem and the restriction that $t \in I_R$ (also since $\frac{kt}{2\ell} \leq \tau$),
\begin{align*} 
\nu(t,\tau) \geq a \delta^\prime \frac{t - \tau}{\tau^a t} & \geq \frac{a\delta^\prime}{\tau^\alpha t}\left[\frac{\ell\tau}{k} - \tau\right] -\frac{a\delta^\prime}{\tau^\alpha t}\left[t - \frac{\ell \tau}{k}\right] \\  
& \geq \frac{a\delta^\prime\tau^{1-a}}{2tk} \\
& \geq \frac{a\delta^\prime}{2^{2-a}\ell^{1-a}(kt)^a}.
\end{align*} 
If we now let $\tilde\delta^{\prime} = \frac{a\delta^\prime}{2^{2-a}}$
and apply \eqref{ineq:Comps}, \eqref{ineq:SobExp} and $(2+\gamma)(s-a) = 1-a$ then we have
\begin{align} 
\mathcal{I}_R & \lesssim \int_{I_R} \frac{kt}{\ell^{1+\gamma}}e^{-\delta\jap{k-\ell,kt-\ell\tau}^s} e^{-\frac{\tilde\delta^{\prime}}{\ell^{1-a}}\abs{kt}^{s-a}} d t \nonumber \\ 
& \lesssim \int_{I_R} \frac{\ell^{\frac{1-a}{s-a}}}{\ell^{1+\gamma} (\tilde\delta^{\prime})^{\frac{1}{s-a}}}e^{-C_s\delta\jap{k-\ell}^s - C_s\delta\jap{kt-\ell\tau}^s} d t \nonumber \\ 
& \lesssim \frac{\ell^{\frac{1-a}{s-a}}}{\ell^{2+\gamma} \delta^{1/s} (\tilde\delta^{\prime})^{\frac{1}{s-a}}} \left(\frac{\ell e^{-C_s\delta\jap{k-\ell}^s}}{k}\right) \nonumber \\ 
& \lesssim \frac{1}{\delta^{2/s} (\tilde\delta^{\prime})^{\frac{1}{s-a}}} e^{-\frac{1}{2}C_s\delta\jap{k-\ell}^s}, \label{ineq:IRHardDual}  
\end{align} 
which is summable in $k$ uniformly in $l$ (the extra power of $k$ in the denominator of the penultimate line came from the time integration). 

Assembling the contributions of \eqref{ineq:INRDual} and \eqref{ineq:IRHardDual}, summing in $k$ and taking the supremum in $\ell$ and $\tau \leq \infty$ completes the proof of Lemma \ref{lem:dual}. 
\end{proof}

The following simple lemma is used in \S\ref{sec:PtwiseRho} above to deduce the pointwise-in-time control on $\rho$.
\begin{lemma} \label{lem:ptwiseTimeResponse} 
Under the bootstrap hypotheses \eqref{ctrl:Boot} we have
\begin{align*} 
\sup_{0 \leq \tau \leq t} \sup_{\ell \in \Integer_\ast^d} \sum_{k \in \Integer_\ast^d} \bar{K}_{k,\ell}(t,\tau) \lesssim \sqrt{K_2} \epsilon\jap{t}.
\end{align*} 
\end{lemma} 
\begin{proof} 
As in the proof of Lemmas \ref{lem:Moment} and \ref{lem:dual}, we may control $f$ by \eqref{ctrl:LowCommLoc} and reduce to dimension one, leaving us to analyze the analogue of \eqref{ineq:IklDual} except without the time integral: 
\begin{align*} 
\mathcal{I}_{k,\ell}(t,\tau) = e^{-\nu(t,\tau)\jap{k,kt}^s}e^{-\delta \jap{k-\ell,kt-\ell\tau}^s}\frac{\jap{\tau}}{\abs{\ell}^\gamma}. 
\end{align*} 
By using \eqref{ineq:Comps} we have,
\begin{align*} 
\mathcal{I}_{k,\ell}(t,\tau) \lesssim e^{-C_s\delta \jap{k-\ell}^s}\jap{\tau}, 
\end{align*} 
which after summing in $k$ and taking the supremum in $\ell$ and $\tau \leq t$ gives the lemma. 
\end{proof}

\section{Final steps of proof} \label{sec:ConcProof}
By Proposition \ref{lem:Boot}, \eqref{ctrl:BootRes} holds uniformly in time. 
By \eqref{ineq:Wbd} and the algebra property \eqref{ineq:GAlg},
\begin{align*} 
\int_0^T \norm{F(t,z+vt) \cdot (\grad_v - t\grad_z)(f^0+f)(t) }_{\cG^{\alpha_0}} dt \lesssim \int_0^T \norm{\rho(t)}_{\cF^{\alpha_0}}\norm{ (\grad_v - t\grad_z)(f^0+f)(t)}_{\cG^{\alpha_0}} dt. 
\end{align*} 
Therefore,  \eqref{ctrl:BootRes}, \eqref{ineq:SobExp}, $\lambda(t) \geq \alpha_0$, \eqref{ineq:f0Loc} and $\sigma > \beta+1$ imply
\begin{align*} 
\int_0^T \norm{F(t,z+vt) \cdot (\grad_v - t\grad_z)(f^0+f)(t) }_{\cG^{\alpha_0}} dt & \lesssim \int_0^T \jap{t}^{-\sigma+1} \norm{A\rho(t)}_{2} \norm{A^{(-\beta)}(f^0 + f)(t)}_{2} dt \\
& \hspace{-2cm} \lesssim \left(\int_0^T \norm{A\rho(t)}^2_{2} dt \right)^{1/2} \left(\int_0^T \jap{t}^{-2\sigma + 2} \norm{A^{(-\beta)}(f^0 + f)(t)}^2_{2} dt\right)^{1/2} \\
& \hspace{-2cm} \lesssim \epsilon.
\end{align*}   
Therefore, we may define $f_\infty$ satisfying $\norm{f_\infty}_{\cG^{\alpha_0}} \lesssim \epsilon$ by the absolutely convergent integral 
\begin{align*} 
f_\infty = h_{in} - \int_0^\infty F(\tau,z+v\tau)\cdot (\grad_v - \tau\grad_z)f(\tau) d\tau.
\end{align*}   
Moreover, again by \eqref{ineq:GAlg}, \eqref{ctrl:BootRes} and \eqref{ineq:SobExp}, 
\begin{align*} 
\norm{f(t) - f_\infty}_{\cG^{\lambda^{\prime}}} 
& \lesssim \int_t^\infty e^{(\lambda^\prime - \alpha_0)\jap{\tau}^s} \jap{\tau}^{-\sigma + 1} \norm{A\rho(\tau)}_{2} \norm{A^{(-\beta)}(f^0+f)(\tau)}_{2} d\tau \\ 
& \lesssim \epsilon e^{\frac{1}{2}(\lambda^\prime - \lambda_0)\jap{t}^s}, 
\end{align*} 
which implies \eqref{ineq:glidingconverg}. 
Then Lemma \ref{lem:Ptwise} implies \eqref{ineq:densitydecay} (since $\sigma > 1/2$), completing the proof of Theorem \ref{thm:Main}. 

We briefly sketch the refinement mentioned in Remark \ref{rmk:LinearTh}. Specifically we verify that the final state predicted by the linear theory is accurate to within $O(\epsilon^2)$.
Indeed, let $f^L$ be the solution to 
\begin{equation}\label{def:VPE_Linear_gliding}
\left\{
\begin{array}{l}
\partial_t f^L + F^L(t,z+vt)\cdot \grad_vf^0 = 0, \\ 
\widehat{F^L(t,z+vt)}(t,k,\eta) = -ik\widehat{W}(k)\widehat{f^L}_k(t,kt) \delta_{\eta = kt}, \\
f^L(t = 0,z,v) = h_{in}(z,v). 
\end{array}
\right.
\end{equation}
By the analysis of \S\ref{sec:LinDamp} we have that $h^L(t,x,v) = f^L(t,x-tv,v)$ satisfies the conclusions of Theorem \ref{thm:Main} 
for $h_\infty^L = h_\infty^L(z,v)$ given by 
\begin{align*} 
h_\infty^L(z,v) = h_{in}(z,v) - \int_0^\infty F^L(t,z+vt) \cdot \grad_v f^0(v) dt.   
\end{align*} 
Consider next the PDE
\begin{align*} 
\partial_t (f-f^L) + (F - F^L)(t,z+vt) \cdot \grad_vf^0 = -F(t,z+vt) \cdot (\grad_v -t\grad_z)f. 
\end{align*}
By treating the right-hand side as a decaying external force, 
the analysis of \S\ref{sec:LinDamp} with $\lambda^\prime$  replaced by $\lambda^{\prime\prime} < \lambda^\prime$, 
then implies 
\begin{align*} 
\norm{f(t)-f^L(t)}_{\lambda^{\prime\prime}} \lesssim_{\lambda^\prime - \lambda^{\prime\prime}} \epsilon^2, 
\end{align*}
which shows that the nonlinearity changes the linear behavior at the expected $O(\epsilon^2)$ order.
Justifying higher order expansions should also be possible, but justifying the convergence of a Newton iteration is significantly more challenging as the constants would need to be quantified carefully. 

\subsubsection*{Acknowledgments}
The authors would like to thank the following people for references and suggestions: Antoine Cerfon, Yan Guo, George Hagstrom, Vladimir \v{S}ver{\'a}k, C\'edric Villani and Walter Strauss. 

\bibliographystyle{plain} \bibliography{eulereqns}

\end{document}